\documentclass[11pt,leqno]{amsart}
\pdfoutput=1
\usepackage{amssymb}
\usepackage{amscd}
\usepackage[all,cmtip,matrix,arrow]{xy}
\usepackage{graphicx}
\usepackage{hyperref}
\xyoption{rotate}
\xyoption{pdf}
\usepackage{xparse}
\usepackage{tikz} 
\usetikzlibrary{arrows,backgrounds,fit,positioning,shapes.symbols,chains}

\renewcommand{\baselinestretch}{1.3}

\newtheorem{theo}{Theorem}[section]
\newtheorem{lemma}[theo]{Lemma}
\newtheorem{defi}[theo]{Definition}
\newtheorem{prop}[theo]{Proposition}

\newtheorem{remark}[theo]{Remark}

\newtheorem{example}[theo]{Example}
\numberwithin{equation}{section}

\mathchardef\mhyphen="2D

\def\CP{\mathbb{CP}}

\def\R{\mathbb{R}}
\def\C{\mathbb{C}}
\def\Z{\mathbb{Z}}

\def\val{\mathrm{val}}
\def\Tw{\mathrm{Tw}}
\def\crit{\mathrm{crit}}
\def\Re{\mathrm{Re}}

\def\F{{\mathcal F}}

\def\O{{\mathcal O}}

\def\W{{\mathcal{W}}}

\def\PP{{\mathbb P}}

\def\pre-tr{\operatorname{pre-tr}}

\def\Hom{\operatorname{Hom}}
\def\End{\operatorname{End}}

\DeclareMathOperator*{\holim}{holim}

\newcommand{\tens}[1]{%
  \mathbin{\mathop{\otimes}\displaylimits_{#1}}%
}

\newcommand{\hhat}{\widehat}
\newcommand{\xto}{\xrightarrow}

\newcommand{\hto}{\hookrightarrow}

\newcommand{\prarr}{\rightrightarrows}

\newcommand{\Coh}{\operatorname{Coh}}

\newcommand{\Sch}{\operatorname{Sch}}

\newcommand{\can}{\operatorname{can}}

\newcommand{\mk}{\mathrm k}

\newcommand{\cF}{{\mathcal F}}

\newcommand{\cO}{{\mathcal O}}

\newcommand{\cL}{{\mathcal L}}
\newcommand{\cM}{{\mathcal M}}
\newcommand{\cN}{{\mathcal N}}

\newcommand{\cE}{{\mathcal E}}
\newcommand{\cW}{{\mathcal W}}
\newcommand{\cU}{{\mathcal U}}

\newcommand{\cK}{{\mathcal K}}

\newcommand{\cX}{{\mathcal X}}

\newcommand{\mX}{{\mathfrak X}}

\newcommand{\Perf}{\operatorname{Perf}}

\newcommand{\supp}{\operatorname{Supp}}

\newcommand{\im}{\operatorname{Im}}

\newcommand{\Ext}{\operatorname{Ext}}

\newcommand{\Spec}{\operatorname{Spec}}
\newcommand{\Spf}{\operatorname{Spf}}

\newcommand{\id}{\operatorname{id}}

\newcommand{\pr}{\operatorname{pr}}

\usepackage{epsf}
\usepackage{amscd}
\newcommand{\m}{\mathfrak{m}}

\title[Lagrangian Floer theory for trivalent graphs and HMS for curves]
{Lagrangian Floer theory for trivalent graphs
and homological mirror symmetry for curves}
\author{Denis Auroux}
\address{Harvard University, Department of Mathematics, 1 Oxford St., Cambridge MA 02138, USA.}
\email{auroux@math.harvard.edu}
\author{Alexander I. Efimov}
\address{The Hebrew University of Jerusalem}
\email{efimov@mccme.ru}
\author{Ludmil Katzarkov}
\address{
College of Arts and Sciences, Department of Mathematics, University of
Miami, Ungar Bldg, 1365
Memorial Dr 515, Coral Gables, FL 33146, USA; and Institute of the
Mathematical Sciences of the
Americas (IMSA), 1365 Memorial Drive, Ungar 515, Coral Gables, FL 33146,
USA.}
\email{l.katzarkov@miami.edu}

\begin{document}

\begin{abstract}
Mirror symmetry for higher genus curves is usually formulated and studied
in terms of Landau-Ginzburg models; however the critical locus of the
superpotential is arguably of greater intrinsic relevance to 
mirror symmetry than the whole
Landau-Ginzburg model.
Accordingly, we propose a new approach to the A-model of the mirror, 
viewed as a trivalent configuration of rational curves together
with some extra data at the nodal points. In this context, we introduce
a version of Lagrangian Floer theory and the Fukaya category for trivalent
graphs, and show that homological mirror symmetry holds,
namely, that the Fukaya category of a trivalent configuration of rational
curves is equivalent to the derived category of a non-Archimedean
generalized Tate curve. To illustrate the concrete nature of this
equivalence, we show how explicit formulas for theta functions 
and for the canonical map of the curve arise naturally under mirror symmetry.
\end{abstract}

\thanks{
DA is partially supported by NSF grants DMS-1937869 and DMS-2202984, and by Simons Foundation grant
\#385573 (Simons Collaboration on Homological Mirror Symmetry).
AE is partially supported by the European Research Council (ERC, CurveArithmetic, 101078157).
LK is partially supported by The Basic Research Program of the National Research
University Higher School of Economics, a Simons Investigator Award,
the Simons Collaboration on Homological Mirror Symmetry, and the
National Science Fund of Bulgaria, National Scientific Program ``Excellent Research and
People for the Development of European Science'' (VIHREN), Project No. KP-06-DV-7.
Part of the research was
conducted while the authors enjoyed the hospitality of IMSA at the University of Miami and the
Laboratory of Mirror Symmetry at HSE Moscow.
The authors thank the referee for useful comments and suggestions about the
previous version of the manuscript.
}

\renewcommand{\baselinestretch}{1.1}

\maketitle
\setcounter{tocdepth}{1}
\tableofcontents

\vfil

\renewcommand{\baselinestretch}{1.3}

\section{Introduction}

Riemann surfaces have been one of the most fruitful sources of examples for
the exploration of homological mirror symmetry, starting with the elliptic curve over
twenty years ago \cite{PZ}, and including some of the
earliest evidence of homological mirror symmetry for varieties of general
type \cite{SeGenus2,EfimovGenusG,AAEKO}. 
Various mirror constructions can be employed to produce mirrors of
Riemann surfaces of arbitrary genus. Most of them rely crucially on the
choice of an embedding into an ambient toric variety, and typically
output a 3-dimensional Landau-Ginzburg model as mirror,
as explained in \cite{AAK} (see also \cite{HoriVafa,Clarke,ChanLauLeung,GKR}).
However there are also some constructions which yield stacky nodal curves 
as mirrors to Riemann surfaces \cite{STZ,GammageShende,LekiliPolishchuk};
the two types of mirrors are in some cases related by a form of Orlov's
generalized Kn\"orrer periodicity \cite{Orlov}.

The various references mentioned above explore the direction of homological mirror
symmetry that compares the Fukaya category of a Riemann surface viewed
as a 2-dimensional symplectic manifold (A-model) with the derived category 
of singularities of the mirror Landau-Ginzburg model (B-model).
Here we study the other direction, comparing the derived category of
coherent sheaves of a smooth curve (B-model) to the Fukaya category of
a mirror Landau-Ginzburg model (A-model). That direction is more challenging, in part
due to the difficulty of defining and working with Fukaya categories of
non-exact Landau-Ginzburg models with non-compact critical loci. 
In the one instance where the
Landau-Ginzburg mirror is exact, namely for pairs of pants, a verification
of the equivalence using the language of microlocal sheaves can be found in \cite{Nadler}.
A comprehensive treatment of this direction of homological mirror symmetry 
for hypersurfaces in $(\C^*)^n$ (the case $n=2$ being of interest here), 
in the language of fiberwise wrapped Fukaya categories of toric 
Landau-Ginzburg models, can be found in \cite{AA}, whereas the example of
a genus 2 curve embedded in an abelian surface (its Jacobian) is treated
using a similar approach (minus the compactness issues)
in Cannizzo's thesis \cite{Cannizzo}. 

The approach pursued in \cite{AA} and \cite{Cannizzo} makes it
clear that the geometry of Landau-Ginzburg mirrors to curves depends
very much on the choice of an embedding: in fact the fiber of the
superpotential is mirror to the ambient space into which the curve is
embedded, with inclusion and restriction functors $i_*,i^*$ on the algebraic side
corresponding under mirror symmetry to a pair of adjoint functors $\cup,\cap$
between the Fukaya category of the Landau-Ginzburg model and that of its
regular fiber. Thus, it should be no surprise that the various
Landau-Ginzburg mirrors to genus 2 curves considered in the papers
\cite{SeGenus2,GKR,AAK,Cannizzo} are actually different: for instance the singular fiber of the mirror in \cite{Cannizzo}
is irreducible, while those of \cite{GKR,AAK} have three irreducible components.
And yet, these mirrors share one common feature, which is that
(after crepant resolution in the case of \cite{SeGenus2}) the critical
loci of the superpotentials always consist of three rational curves meeting
in two triple points. 
Similarly, for a smooth proper curve of genus $g\geq 2$ curve,
the critical locus of a mirror superpotential (possibly after crepant 
resolution of the total space) consists of a configuration of $3g-3$
rational curves meeting in $2g-2$ triple points. 

For the other direction of mirror symmetry, it has been proposed that the algebraic geometry of the Landau-Ginzburg
model can be replaced by direct consideration of this critical locus,
equipped with additional data making it a ``perverse curve''
\cite{GKR,Ruddat}; this is generally sound given the local nature of the
derived category of singularities, which was shown by Orlov to only depend
on the formal neighborhood of the critical locus. 
Our goal in this manuscript is to do the same for the symplectic geometry
(A-model), in order to arrive that a picture of homological mirror symmetry
for curves that allows for explicit computations and is manifestly
independent of a choice of embedding; there is however a price to pay,
due to the non-local nature of Fukaya-Floer theory and
the fact that restriction to the critical locus hides away instanton
corrections that may be present in the global symplectic geometry of the Landau-Ginzburg
model.

The general features of our construction are motivated by considering
the simplest example, which serves as a building block for all others:

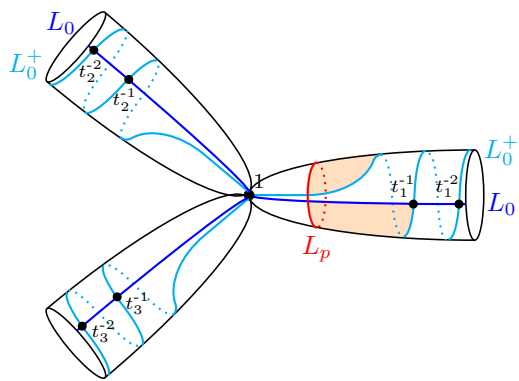
\begin{figure}[t]
\begin{tikzpicture}[scale = 0.6]
\fill[orange!25!white] (1.3,0) arc (270:345:1.5 and 1) arc (165:90:0.22) arc (114:135:5 and 1) arc (90:180:0.2 and 0.72);
\fill[orange!25!white] (1.3,-0.1) arc (190:270:0.2 and 0.72) arc (-135:-109:5 and 1) arc (-90:-20:0.25 and 1.15);
\draw[semithick] (5,1) arc (90:270:5 and 1);
\draw[semithick] (5,0) ellipse (0.2 and 1);
\draw[cyan,thick](0,0)--(1.3,0) arc (270:345:1.5 and 1) arc (165:90:0.22);
\draw[cyan,thick,dotted] (2.92,0.9) arc (90:10:0.25 and 1.1) arc (190:270:0.25 and 1.1);
\draw[cyan,thick] (3.42,-0.92) arc (-90:-10:0.25 and 1.15) arc (170:90:0.25 and 1.15);
\draw[cyan,thick,dotted] (3.98,0.92) arc (90:10:0.2 and 1.16) arc (190:270:0.2 and 1.16);
\draw[cyan,thick] (4.42,-0.98) arc (-90:-10:0.25 and 1.2) arc (170:90:0.25 and 1.2);
\draw[cyan](5.6,1) node {\small $L_0^+$};
\draw[blue,thick](0,0) arc (180:270:4.8 and 0.2);
\draw[blue](5.6,-0.2) node {\small $L_0$};
\draw[red,thick](1.5,0.72) arc (90:270:0.2 and 0.72);
\draw[red,thick,dotted](1.5,0.72) arc (90:-90:0.2 and 0.72);
\draw[red](1.5,-1.2) node {\small $L_p$};
\fill(3.64,-0.2) circle (0.1) (3.42,0.2) node {\tiny $t_1^{\mhyphen1}$};
\fill(4.65,-0.2) circle (0.1) (4.4,0.2) node {\tiny $t_1^{\mhyphen2}$};
\fill(0,0) circle (0.12) (0.22,0.3) node {\SMALL $1$};
\begin{scope}[rotate=140,shift={(0.05,-0.05)}]
\draw[semithick] (5,1) arc (90:270:5 and 1);
\draw[semithick] (5,0) ellipse (0.2 and 1);
\draw[cyan,thick](0,0)--(1.3,0) arc (270:345:1.5 and 1) arc (165:90:0.22);
\draw[cyan,thick,dotted] (2.92,0.9) arc (90:10:0.25 and 1.1) arc (190:270:0.25 and 1.1);
\draw[cyan,thick] (3.42,-0.92) arc (-90:-10:0.25 and 1.15) arc (170:90:0.25 and 1.15);
\draw[cyan,thick,dotted] (3.98,0.92) arc (90:10:0.2 and 1.16) arc (190:270:0.2 and 1.16);
\draw[cyan,thick] (4.42,-0.98) arc (-90:-10:0.25 and 1.2) arc (170:90:0.25 and 1.2);
\draw[cyan](5.6,1) node {\small $L_0^+$};
\draw[blue,thick](0,0) arc (180:270:4.8 and 0.2);
\draw[blue](5.6,-0.2) node {\small $L_0$};
\fill(3.64,-0.2) circle (0.1) (3.42,0.2) node {\tiny $t_2^{\mhyphen1}$};
\fill(4.65,-0.2) circle (0.1) (4.4,0.2) node {\tiny $t_2^{\mhyphen2}$};
\end{scope}
\begin{scope}[rotate=-140,shift={(0.05,0.05)}]
\draw[semithick] (5,1) arc (90:270:5 and 1);
\draw[semithick] (5,0) ellipse (0.2 and 1);
\draw[cyan,thick](0,0)--(1.3,0) arc (270:345:1.5 and 1) arc (165:90:0.22);
\draw[cyan,thick,dotted] (2.92,0.9) arc (90:10:0.25 and 1.1) arc (190:270:0.25 and 1.1);
\draw[cyan,thick] (3.42,-0.92) arc (-90:-10:0.25 and 1.15) arc (170:90:0.25 and 1.15);
\draw[cyan,thick,dotted] (3.98,0.92) arc (90:10:0.2 and 1.16) arc (190:270:0.2 and 1.16);
\draw[cyan,thick] (4.42,-0.98) arc (-90:-10:0.25 and 1.2) arc (170:90:0.25 and 1.2);
\draw[blue,thick](0,0) arc (180:270:4.8 and 0.2);
\fill(3.64,-0.2) circle (0.1) (3.38,0.2) node {\tiny $t_3^{\mhyphen1}$};
\fill(4.65,-0.2) circle (0.1) (4.36,0.2) node {\tiny $t_3^{\mhyphen2}$};
\end{scope}
\end{tikzpicture}
\caption{Wrapped Floer homology in the mirror of the pair of pants,
$M=\bigcup_{i=1}^3 (\C,0)$}\label{fig:pantsmirror}
\end{figure}

\begin{example} \label{ex:pants}
Let $X$ be the pair of pants, i.e.\ $\PP^1$ minus three points. The 
mirror Landau-Ginzburg model is $(\C^3,-xyz)$, with critical locus the
union of the three coordinate axes in $\C^3$, i.e.\ the mirror we
consider is a configuration $M=\bigcup_{i=1}^3 (\C,0)$ consisting of three copies of the complex plane $\C$ meeting in 
a triple point at the origin. The mirror to the structure sheaf $\O_X$
is a Lagrangian graph $L_0=\bigcup_{i=1}^3 \R_{\geq 0}$ consisting of the real positive
axis in each component of~$M$. The wrapped Floer cohomology of $L_0$ inside
$M$ has an additive basis consisting of one generator at the origin, and
three infinite sequences of generators in each of the ends of $M$ (see
Figure~\ref{fig:pantsmirror}); these
correspond respectively to the constant function~1 and to successive
powers of the inverses of coordinates~$t_i$ near the three punctures of $X$.
Considering the multiplicative structure on $HW^*(L_0,L_0)$, however, it is 
clear that the structure maps of Lagrangian Floer theory in $M$ must 
include holomorphic discs that ``propagate'' from one component to another
through the origin, as we explain further 
\hbox{in \S\S \ref{s:pantsmotivate}--\ref{s:Amodel}}.
\end{example}

In order to pass from the pair of pants to the
general case, recall first that mirror symmetry is expected to hold near
the ``large complex structure limit'', i.e., in a non-Archimedean setting.
Lee's thesis \cite{Lee} illustrates the general expectation that mirror 
symmetry for curves is compatible with 
pair-of-pants decompositions induced by maximal degenerations.
Namely, the construction in \cite{AAK} produces a toric Landau-Ginzburg 
model from a maximally degenerating family of complex curves in $(\C^*)^2$
near the tropical limit; this mirror is built out of standard affine charts 
$(\C^3,-xyz)$ glued to each other by toric coordinate changes in a
manner that reflects the combinatorial pair-of-pants decomposition of the
curve induced by the tropical limit.  Lee constructs a
version of the wrapped Fukaya category of the curve that can be viewed as
a \v Cech model for this pair-of-pants decomposition, and uses it to prove an
equivalence with the derived category of singularities of the mirror
\cite{Lee}. 

While the language of degenerating families of complex curves
is convenient when the curve lives on the symplectic
side of mirror symmetry, in our setting it is more fruitful to consider
a curve $X$ defined over a non-Archimedean field $K$, the
Novikov field of power series with real exponents in a formal
variable $T$, which is the natural field of definition of Fukaya
categories in the non-exact setting. We consider non-Archimedean
curves obtained by smoothing a maximally degenerate nodal configuration
$X^0$, given by a union of rational curves with three marked points, 
identified pairwise across components according to a trivalent graph.

\begin{defi}\label{def:combidata}
The {\bf combinatorial data} for our construction is the following.
Let $G$ be a finite (unoriented) graph, with set of vertices $V$ and 
set of edges $E,$ such that each vertex $v\in V$ has degree $3,$
and without loops (edges from a vertex to itself). We write $e/v$ when $e\in E$ is incident to $v\in V.$

For each $v\in V,$ we take $X^0_{v}$ to be a copy of $\PP^1_{\Z},$ and for each $e/v,$ we fix a $\Z$-point $x_{e/v}\in X^0_{v},$ so that $x_{e/v},$ and $x_{e'/v}$ are disjoint for $e\ne e'.$

For each $e/v,$ we choose a coordinate $t_{e/v}$ on $X_{v}^0,$ such that $t_{e/v}(x_{e/v})=0$ and $t_{e/v}$ takes values $1,\infty$ at the other two marked points. 

We also introduce formal variables $\{q_{e}\}_{e\in E},$ which will be set
to elements of the Novikov field with valuation $\val(q_e)=A_e>0$.
\end{defi}

We explain in Section \ref{s:Bmodel} how to produce generalized
Tate curves by smoothing the nodal curve
$X^0=\left(\bigsqcup_{v \in V} X^0_v\right)/(x_{e/v}\sim x_{e/v'}\ \forall
e\in E, v\neq v')$. In terms of rigid analytic geometry, the construction amounts to replacing 
each node of $X^0$ by its smoothing defined in terms of local coordinates
by $t_{e/v}t_{e/v'}=q_e$, producing a curve $X_K$ on which the valuations of the 
coordinates $t_{e/v}$ naturally take values in a metric graph modelled on
$G$ and with edge lengths $A_e=\val(q_e)$.

The A-side is a trivalent configuration $M$ of $2$-spheres, where the
components are in bijection with $E,$ and the nodes are in bijection with 
$V.$ (Thus each component of $M$ passes through two triple points).
We denote by $\{A_e\}_{e\in E}$ the symplectic areas of the components.
The Fukaya category $\F(M)$ is defined in Section \ref{s:Amodel}. 
Besides simple closed curves in the complement of the nodes, this
category also includes objects which are embedded trivalent
graphs in $M$, consisting of one arc joining the two nodes inside each 
component; the Floer theory of these objects involves configurations of
holomorphic discs which propagate through the vertices, according to rules
determined by the coordinates $t_{e/v}$ chosen as part of the combinatorial
data (see \S \ref{s:Amodel}).

Our main result is then:

\begin{theo}\label{thm:main}
Given combinatorial data as above, and setting $q_{e}=T^{A_{e}},$ the 
Fukaya category $\F(M)$ is equivalent to $\Perf(X_K)$.
\end{theo}

\begin{remark}\label{rmk:main}
Equipping $M$ with a B-field or {\em bulk deformation} of the
Fukaya category gives an extension of this result to arbitrary values of
$q_e\in K$ with $\val(q_e)=A_e>0$.
Also, the requirement that $G$ has no loops is purely for convenience
of notation, so that the half-edges of $G$ can be labelled unambiguously;
apart from the notation issues, the result extends immediately to the case
with loops, with the same proof.
\end{remark}

\begin{remark}\label{rmk:opencase}
On the A-side we can also allow some components of $M$ to be $S^2\setminus
\{pt\}$, i.e.\ the complex plane $\C$, with a single triple point on each
such component. These noncompact components are equipped with a symplectic 
form of infinite area, and the Fukaya category can be defined either with wrapping at infinity or with a stop at infinity.
Combinatorially this amounts to allowing $G$ to have ``external
edges'' (so that each vertex still has three edges attached to it, but
external edges do not connect to another vertex; we do not associate
a formal parameter $q_e$ to the external edge). On the B-side, we do not
attach any other component to $X^0_v$ at the marked point $x_{e/v}$
corresponding to an external edge, but in the wrapped case we delete the point $x_{e/v}$ 
from $X^0$ and $X$; in the stopped case we do not do anything at $x_{e/v}$.
For instance, the pair of pants (Example \ref{ex:pants}) corresponds to
the case of a single vertex, with three external edges.
The analogue of Theorem \ref{thm:main} in this setting follows readily from
our proof of the theorem.
\end{remark}

\begin{remark}We mention that one can verify explicitly that the product structure on the ring of regular
functions of an affine elliptic curve matches the structure constants of
the Floer product on the A-model (which in this case has one component of the form $S^2\setminus\{pt\},$ with wrapping at infinity, and one component of the form $S^2/(p\sim q)$).\end{remark}

Another extension of Theorem \ref{thm:main} is to consider curves near
a non-maximal degeneration, i.e.\ graphs whose vertices may
have valency greater than 3. On the B-side, this amounts to considering
curves obtained by smoothing nodal configurations where each
$\PP^1$ may carry more than three nodes (we accordingly relax the requirements 
on the local coordinates $t_{e/v}$ used to construct $X$). On the A-side, 
this amounts to allowing $M$ to have nodes where more than three 
components attach to each other; objects are still supported on graphs
consisting of one arc joining the two nodes in each component of $M$.
Our proof of Theorem \ref{thm:main} can be adapted to this setting to
establish homological mirror symmetry over the entire moduli space of rigid analytic
curves.

The rest of this paper is organized as follows. Section \ref{s:pantsmotivate} 
discusses the case of the pair of pants and
the symplectic geometry of the Landau-Ginzburg model $(\C^3,-xyz)$ in order
to motivate some of the key features of our A-model construction; we also
highlight some key differences between our construction and other approaches.
Section \ref{s:Amodel} is devoted to the definition of our A-model (the Fukaya category 
of a trivalent configuration of spheres).
In Section \ref{s:Bmodel} we describe the construction of the B-model (the curve
$X$) from the combinatorial data, and Theorem \ref{thm:main} is then proved in
Section \ref{s:equivalence}; the argument involves a version of
the Fukaya category $\F(M)$ with Hamiltonian perturbations
(similar to the construction in \cite{Lee}),
homological perturbation theory, and a restriction diagram 
for decompositions of $X_K$ and $M$ into pairs of pants and their mirrors.  
Sections \ref{ss:theta} and \ref{s:canonicalmap} illustrate the very concrete nature of the
equivalence of A- and B-models in our setup (in sharp contrast with Fukaya
categories of Landau-Ginzburg models): we show how theta functions arise
from the construction of the mirror functor, and we
determine explicitly the canonical map of the curve $X_K$ and its A-model 
counterpart for a general trivalent graph. 
Finally, in Section \ref{s:higherdim} we give a tentative (and highly
speculative) description of how our construction and results ought to generalize
to the higher-dimensional setting.

\section{Motivation and comparison with other approaches} \label{s:pantsmotivate}

In this section we discuss some features of the symplectic geometry of
Landau-Ginzburg mirrors to plane curves, focusing in particular on the
case of $(\C^3,-xyz)$ (mirror to the pair of pants).
This material is useful to understand the rationale
for the construction described in Section \ref{s:Amodel}, and some of its
key differences with other approaches, but it is not
part of the main argument; the reader who wishes to get straight to the
precise formulation of our construction and the proof of Theorem
\ref{thm:main} can skip this section altogether.

\subsection{Motivation: the mirror of the pair of pants} 

We first turn our attention to the symplectic geometry of 
the Landau-Ginzburg model $(\C^3,-xyz)$ and the manner in which it is
reflected in our A-model construction in the case of the pair of pants (Example
\ref{ex:pants}). 

The general philosophy of trying to reduce the symplectic geometry of
a Landau-Ginzburg model to that of its critical locus
is motivated by the well-understood case of Lefschetz fibrations and, less
well understood but closer to our setting, Morse-Bott fibrations.
For instance, the construction in \cite{AAK} associates to
a smooth elliptic curve $X$ (embedded into a toric surface) a 
3-dimensional Landau-Ginzburg model $(Y,W)$ whose singularities are
Morse-Bott along a smooth elliptic curve $M=\crit(W)\subset Y$, which 
is in fact the ``usual'' mirror of $X$. We can then upgrade an object
of the Fukaya category of $M$ (i.e., a simple closed curve with a local
system) to a Lagrangian {\em thimble} in $Y$, obtained by parallel transport
over an arc connecting the critical value of $W$ (the origin) to $+\infty$:
to $L\in \F(M)$ we associate $\mathcal{T}(L)\in \F(Y,W)$, the admissible
Lagrangian consisting of those points of $Y$ where the negative gradient flow of 
$\Re(W)$ with respect to a K\"ahler metric converges to a point of $L$
(together with the pullback local system). In this example the construction
gives rise to a functor $\mathcal{T}:\F(M)\to \F(Y,W)$, which is in fact an
equivalence; we note however that for a general Morse-Bott fibration the
situation can be slightly more complicated (see e.g.\ \cite[Corollary 7.8]{AAK}).

The case of interest to us falls outside of the Morse-Bott setting:
we consider the Landau-Ginzburg model $(\C^3,-xyz)$ and its fiberwise 
wrapped Fukaya category. The objects of $\F(\C^3,-xyz)$ are admissible Lagrangian 
submanifolds of $\C^3$, whose image under the projection $W=-xyz:\C^3\to \C$
consists, near infinity, of one or more rays pointing towards $\Re(W)\to
+\infty$, while morphisms involve Hamiltonian perturbations that act on
Lagrangians by wrapping at infinity within the fibers of $W$ and by pushing
rays in the base of the fibration slightly in the counterclockwise direction
\cite{AA}. 

The Fukaya category of a Landau-Ginzburg model is related to that of the
regular fiber (in this case, the wrapped Fukaya category of $(\C^*)^2$) 
by a pair of spherical functors \cite{AG,AS},
often denoted $\cup$ and $\cap$, which we briefly describe.
On objects, the cup functor  (also called Orlov functor)
$$\cup:\W((\C^*)^2)\to \F(\C^3,-xyz)$$
takes a Lagrangian submanifold $\ell$ of $(\C^*)^2\simeq \{xyz=1\}=W^{-1}(-1)$ and 
considers its parallel transport in the fibers of $W=-xyz$ over a U-shaped
arc to produce an admissible Lagrangian submanifold $\cup \ell\subset \C^3$.
The cap functor $$\cap:\F(\C^3,-xyz)\to \Tw\,\W((\C^*)^2)$$
restricts an admissible Lagrangian $\mathbf{L}\subset \C^3$ to the fiberwise
Lagrangians in its ends at $\Re(W)\to\infty$; if there is only one such
end this produces an object of $\W((\C^*)^2)$, otherwise one obtains a
twisted complex built from the objects in the various ends of $\mathbf{L}$, with
connecting differentials given by counts of holomorphic discs in
$\C^3$ with boundary in $\mathbf{L}$ (with one outgoing strip-like end towards
$\Re(W)\to \infty$).
The argument in \cite{AA} proves homological mirror symmetry for the pair
of pants (and for other very affine hypersurfaces) in a manner compatible with
these functors, namely:

\begin{theo}[\cite{AA}]
$\F(\C^3,-xyz)$ is equivalent to the derived category of the pair of
pants $X=\{(x_1,x_2)\in (K^*)^2\,|\,1+x_1+x_2=0\}$, and we have
a commutative diagram
$$\xymatrix{
\!\!\!\!\Tw\,\F(\C^3,-xyz) \ar@<.5ex>[r]^{\cap} \ar[d]^{\simeq} 
& \Tw\,\W((\C^*)^2) \ar@<.5ex>[l]^{\cup}\ar[d]^{\simeq} \\
\Perf(X)\ar@<.5ex>[r]^{i_*} & \Perf((K^*)^2) \ar@<.5ex>[l]^{i^*}
}
$$
i.e.\ the functors $\cap$ and $\cup$ correspond under mirror symmetry to the 
inclusion and restriction functors $i_*$ and $i^*$ between the derived categories
of $X$ and of the ambient space $(K^*)^2$.
\end{theo}

The critical locus $M=\crit(W)$ is the union of the coordinate axes in $\C^3$,
hence not smooth, but the
singularities of $W$ are Morse-Bott away from the origin; given an embedded 
Lagrangian submanifold $L_p$ in the smooth part of $M$, we can 
build a thimble $\mathcal{T}(L_p)\subset \C^3$ by parallel transport over
the real positive axis. For example, if we use the standard K\"ahler form of
$\C^3$, and start from $L_p=\{(x,0,0)\,|\,|x|=r\}\subset M$, we obtain the
Lagrangian $$\mathcal{T}(L_p)=\bigl\{(x,y,z)\in\C^3\,\big|\,|x|^2-r^2=|y|^2=|z|^2,\,
-xyz\in \R_{\geq 0}\bigr\}$$
(in \cite{AA} a different toric K\"ahler form is used for technical reasons, but
this is immaterial to our discussion).
$\mathcal{T}(L_p)$ can be equipped with a (unitary, i.e.\ valuation-preserving)
local system of rank 1 over the Novikov field $K$, and should also be endowed
with a bounding cochain to cancel out the Floer-theoretic obstruction arising from the holomorphic discs bounded by
$\mathcal{T}(L_p)$ (namely, the disc of radius $r$ in the $x$-axis, whose symplectic
area we denote by $A$, and its multiple covers); this yields a
so-called Aganagic-Vafa Lagrangian brane in $\F(\C^3,-xyz)$, which is mirror
to the skyscraper sheaf $\O_p$ of a point $p$ of the pair of pants $X=\{1+x_1+x_2=0\}$
with $\val(x_1(p))=A$; the values of the coordinates $(x_1,x_2)$ depend on
the choice of local system and bounding cochain. 
The vanishing cycle, i.e.\ the boundary at infinity $\Lambda_p=\cap\mathcal{T}(L_p)$,
is a Lagrangian torus in $(\C^*)^2$ equipped with a rank~1 local system
(whose holonomy is nontrivial even along the $S^1$-factor that bounds a
disc inside $\mathcal{T}(L_p)$, due to the obstruction-cancelling bounding cochain);
it is in fact mirror to the skyscraper sheaf of the point $p$ in $(K^*)^2$,
as expected given that $\cap$ corresponds to $i_*$ under mirror symmetry.

Since the object which corresponds to the structure sheaf of $X$ should intersect each of
the point objects once, it is natural to consider 
the singular Lagrangian $L_0=\bigcup_{i=1}^3 \R_{\geq 0}$ consisting of
the union of the real positive axes in the three components of $M$. 
Parallel transport can be used to produce a piecewise linear Lagrangian
cycle in $(\C^3,-xyz)$ out of $L_0$, whose intersection $\Lambda_0^{PL}$ 
with a smooth fiber $\{-xyz=c\gg 0\}$ near infinity (the ``PL vanishing cycle'')
is the union of the semi-infinite cylinders $\{|x|\geq |y|=|z|,\ \arg(x)=0\}$,
$\{|y|\geq |x|=|z|,\ \arg(y)=0\}$, $\{|z|\geq |x|=|y|,\ \arg(z)=0\}$ and two
triangular portions of the torus $\{|x|=|y|=|z|\}$.
However it is not clear how one could modify this construction to produce a smooth admissible Lagrangian
in $\C^3$. 

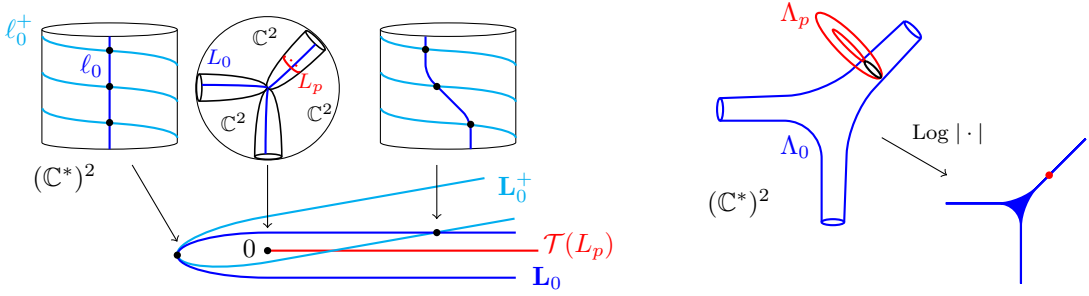
\begin{figure}[t]
\begin{center}
\begin{tikzpicture}[scale=0.6]
\draw(-2.5,1.7) node {\small $(\C^*)^2$};
\draw[cyan](-3.5,5) node {\small $\ell_0^+$};
\draw[blue](-1.9,4.2) node {\small $\ell_0$};
\draw[blue](8.2,-0.5) node {\small $\mathbf{L}_0$};
\draw[red](8.9,0.2) node {\small $\mathcal{T}(L_p)$};
\draw[cyan](7.5,1.5) node {\small $\mathbf{L}_0^+$};
\draw[red,thick](2,0.1)--(8,0.1);
\draw[cyan,thick,rotate=10](7,0.5)--(2,0.5) arc (90:270:2 and 0.5)--(7.5,-0.5);
\draw[blue,thick](7.5,0.5)--(2,0.5) arc (90:270:2 and 0.5)--(7.5,-0.5);
\fill(2,0.1) circle (0.08) (1.6,0.15) node {\small 0};
\fill(0,0) circle (0.08);
\fill(5.75,0.5) circle (0.08);
\draw[->](-1,2)--(0,0.3);
\draw[->](2,2)--(2,0.6);
\draw[->](5.75,2)--(5.75,0.8);
\draw(-1.5,5) ellipse (1.5 and 0.15);
\draw(-3,5)--(-3,2.5) arc (180:360:1.5 and 0.15)--(0,5);
\draw[cyan,thick](0,2.6) arc (0:70:2.27 and 0.35) arc (-110:-180:2.27 and 0.35);
\draw[cyan,thick](0,3.4) arc (0:70:2.27 and 0.35) arc (-110:-180:2.27 and 0.35);
\draw[cyan,thick](0,4.2) arc (0:70:2.27 and 0.35) arc (-110:-180:2.27 and 0.35);
\draw[blue,thick](-1.5,2.35)--(-1.5,4.85);
\fill(-1.5,2.94) circle (0.08);
\fill(-1.5,3.74) circle (0.08);
\fill(-1.5,4.54) circle (0.08);
\draw(6,5) ellipse (1.5 and 0.15);
\draw(4.5,5)--(4.5,2.5) arc (180:360:1.5 and 0.15)--(7.5,5);
\draw[cyan,thick](7.5,2.6) arc (0:70:2.27 and 0.35) arc (-110:-180:2.27 and 0.35);
\draw[cyan,thick](7.5,3.4) arc (0:70:2.27 and 0.35) arc (-110:-180:2.27 and 0.35);
\draw[cyan,thick](7.5,4.2) arc (0:70:2.27 and 0.35) arc (-110:-180:2.27 and 0.35);
\draw[blue,thick,rounded corners](6.5,2.35)--(6.5,3)--(5.5,4)--(5.5,4.85);
\fill(6.5,2.9) circle (0.08);
\fill(5.75,3.74) circle (0.08);
\fill (5.5,4.56) circle (0.08);
\begin{scope}[shift={(2,3.7)},scale=0.3]
\begin{scope}[rotate=45]
\draw(0,0) circle (5.3);
\draw(2.6,2.8) node {\SMALL $\C^2$};
\draw(1.5,-4) node {\SMALL $\C^2$};
\draw(-3.5,0) node {\SMALL $\C^2$};
\draw[semithick] (5,1) arc (90:270:5 and 1);
\draw[semithick] (5,0) ellipse (0.2 and 1);
\draw[blue,thick](0,0) arc (180:270:4.8 and 0.2);
\draw[red,thick](2.5,0.8) arc (90:270:0.2 and 0.8);
\draw[red,thick,dotted](2.5,0.8) arc (90:-90:0.2 and 0.8);
\draw[red](2.5,-2) node {\SMALL $L_p$};
\end{scope}
\begin{scope}[rotate=180,shift={(0.05,-0.05)}]
\draw[semithick] (5,1) arc (90:270:5 and 1);
\draw[semithick] (5,0) ellipse (0.2 and 1);
\draw[blue,thick](0,0) arc (180:270:4.8 and 0.2);
\draw[blue](3.5,-2.2) node {\SMALL $L_0$};
\end{scope}
\begin{scope}[rotate=-90,shift={(0.05,0.05)}]
\draw[semithick] (5,1) arc (90:270:5 and 1);
\draw[semithick] (5,0) ellipse (0.2 and 1);
\draw[blue,thick](0,0) arc (180:270:4.8 and 0.2);
\end{scope}
\end{scope}
\end{tikzpicture}
\qquad
\begin{tikzpicture}[scale=0.5]
\draw[blue,thick](-3,0) ellipse (0.08 and 0.3);
\draw[blue,thick](0,-3) ellipse (0.3 and 0.08);
\draw[blue,thick,rounded corners=5mm](-3,0.3)--(-0.2,0.4)--(1.9,2.35);
\draw[blue,thick,rounded corners=5mm](0.3,-3)--(0.4,-0.2)--(2.35,1.9);
\draw[blue,thick,rounded corners=5mm](-0.3,-3)--(-0.3,-0.3)--(-3,-0.3);
\draw[blue,thick,rotate=45](3,0) ellipse (0.08 and 0.3);
\begin{scope}[rotate=45]
\draw[red,thick](1.5,0.9) ellipse (0.08 and 0.6);
\draw[red,thick](1.5,0.9) ellipse (0.3 and 1.2);
\draw[thick](1.5,0) ellipse (0.08 and 0.3); 
\end{scope}
\begin{scope}[rotate=-30]
\draw[->](1.5,0)--(3.5,0); \draw (3,1) node {\SMALL Log\,$|\cdot|$};
\end{scope}
\filldraw[blue,thick,rounded corners=3mm](3,-2.5)--(5,-2.5)--(7,-0.5)--(5,-2.5)--(5,-5)--(5,-2.5)--(3,-2.5);
\fill[red](5.75,-1.75) circle (0.1);
\draw[red](-1,2.5) node {\small $\Lambda_p$};
\draw[blue](-1,-1) node {\small $\Lambda_0$};
\draw(-2.5,-2.5) node {\small $(\C^*)^2$};
\end{tikzpicture}
\end{center}
\caption{Left: $\mathbf{L}_0=\cup \ell_0\in \F(\C^3,-xyz)$
and the thimble $\mathcal{T}(L_p)$. Right: 
the tropical Lagrangian pair of pants $\Lambda_0\simeq \cap\mathbf{L}_0\subset (\C^*)^2,$
and $\Lambda_p=\cap\mathcal{T}(L_p)$.}
\label{fig:pantsLGmirror}
\end{figure}

Thus, the argument in \cite{AA} bypasses attempts
to construct a thimble and instead considers
the object $\mathbf{L}_0=\cup \ell_0\in \F(\C^3,-xyz)$ obtained by parallel transport of
$\ell_0=(\R_+)^2\subset (\C^*)^2$ over a U-shaped arc in the complex plane;
see Figure \ref{fig:pantsLGmirror}.

The proof of homological mirror symmetry in \cite{AA} relies on a direct calculation to show
that the fiberwise wrapped Floer complex of $\mathbf{L}_0$ is given by
$$\End(\mathbf{L}_0)\simeq
\left\{CW^*(\ell_0,\ell_0)[1] \stackrel{\partial}{\longrightarrow}
CW^*(\ell_0,\ell_0)\right\}\simeq 
\left\{K[x_1^{\pm 1},x_2^{\pm 1}][1]\stackrel{1+x_1+x_2}{\xrightarrow{\hspace{1cm}}}
 K[x_1^{\pm 1},x_2^{\pm 1}]\right\},$$
and that the cohomology algebra agrees with the ring of functions of the pair of pants.
(The two terms in the complex correspond to intersections between $\mathbf{L}_0$ and
its positive perturbation $\mathbf{L}_0^+$ inside the two fibers of $W$ depicted
on Figure \ref{fig:pantsLGmirror} left; each of these amounts to the wrapped
Floer cohomology of $\ell_0$ in $(\C^*)^2$, and the connecting differential is a count
of holomorphic sections over the bigon visible in the base of the fibration.)
While this calculation leads to a proof of homological mirror symmetry 
for the pair of pants $X$ and the Landau-Ginzburg model $(\C^3,-xyz)$, it does
not shed light on how the endomorphisms of $\mathbf{L}_0$
might relate to a version of wrapped Floer homology for $L_0=\bigcup_{i=1}^3 \R_{\geq 0}$
inside $M$ (cf.\ Figure~\ref{fig:pantsmirror}): indeed, $HW^0_M(L_0,L_0)$ comes with a distinguished basis (up to scaling)
corresponding to Floer generators, while $H^0\End(\mathbf{L}_0)$ arises as a
quotient of $HW^0_{(\C^*)^2}(\ell_0,\ell_0)\simeq K[x_1^{\pm 1},x_2^{\pm 1}]$
by the ideal generated by $1+x_1+x_2$, and does not have a preferred basis.

A more promising approach stems from the observation that, even though
$\mathbf{L}_0$ has two ends at $\Re(W)\to +\infty$ and hence maps under
the cap functor to a twisted complex rather than a single
Lagrangian, specifically the mapping cone 
$\bigl\{\ell_0[1]\stackrel{1+x_1+x_2}{\xrightarrow{\hspace{1cm}}}\ell_0\bigr\}
\in \Tw\,\W((\C^*)^2)$, in fact this twisted complex can be represented
geometrically by an embedded Lagrangian $\Lambda_0\subset (\C^*)^2$, the
{\em tropical Lagrangian pair of pants} introduced independently by
Hicks, Matessi and Mikhalkin \cite{Hicks,Matessi,Mikhalkin}; not coincidentally,
$\Lambda_0$ is in fact a smoothing of the PL vanishing cycle $\Lambda_0^{PL}$.
We note that the construction given by Hicks explicitly realizes the tropical Lagrangian
pair of pants as a mapping cone between $\ell_0$ and its image under
the monodromy of the fibration $W,$ making it apparent that
$\cap \mathbf{L}_0\simeq \Lambda_0$ \cite{Hicks}. This is relevant because
the map $\Hom(\mathbf{L}_0,\mathbf{L}_0)\to \Hom(\cap
\mathbf{L}_0,\cap \mathbf{L}_0)$ induced by the cap functor is injective
(in fact this holds for every object of $\F(\C^3,-xyz)$, because the 
exact triangle of functors involving the 
counit of the adjunction $\cup\cap\to \id$ is split). Therefore
$H^*\End(\mathbf{L}_0)$ naturally arises as a summand in the wrapped Floer
cohomology $HW^*(\Lambda_0,\Lambda_0)$ in $(\C^*)^2$, specifically
it is the degree zero part $HW^0(\Lambda_0,\Lambda_0)$. 
This corresponds under mirror symmetry to the fact
that $\Hom^0(i_*\O_X,i_*\O_X)\simeq \End(\O_X)$.
Summarizing, we have:

\begin{prop}
The degree zero wrapped Floer cohomology $HW^0(\Lambda_0,\Lambda_0)$ of
the tropical Lagrangian pair of pants $\Lambda_0$ inside $(\C^*)^2$ is
isomorphic (as a ring) to $\End(\O_X)$, i.e.\ the ring of functions of the pair of pants $X$.
\end{prop}

Thus, our definition of the wrapped Floer cohomology of $L_0$ inside $M$
is motivated by an analogy with the degree 0 wrapped Floer cohomology
of $\Lambda_0$ in $(\C^*)^2$. 
$HW^0(\Lambda_0,\Lambda_0)$ has one generator $e$ corresponding 
to the minimum of the wrapping Hamiltonian, representing the
identity element for the Floer product, and one infinite sequence 
of generators $\theta_{i,k}$, $k\geq 1$, $1\leq i\leq 3$ in each of the three legs of 
$\Lambda_0$ (corresponding to trajectories of the Hamiltonian flow which
wrap $k$ times in the $\arg(x)$ (resp.\ $\arg(y)$, $\arg(z)$) direction).

\begin{lemma} \label{l:floergens_as_rationalfunctions}
Under the isomorphism $HW^0(\Lambda_0,\Lambda_0)\simeq \End(\O_X)$, 
the Floer generator $\theta_{i,k}$ corresponds to a regular function on $X$ which,
as a rational function on $\PP^1$, has a pole of order $k$ at the
$i^\text{th}$ puncture, and no other poles.
\end{lemma}

\proof Recall that the wrapped Floer
complex of $\Lambda_0$ is the direct limit of the Floer complexes
$CF^*(\Lambda_{n},\Lambda_0)$, where $\Lambda_{n}$ is the image of
$\Lambda_0$ under a Hamiltonian diffeomorphism which wraps each of the
three legs $n$ times at infinity.
The direct limit is taken with respect to the continuation maps 
$CF^*(\Lambda_n,\Lambda_0)\to CF^*(\Lambda_{n+1},\Lambda_0)$ associated
to positive Hamiltonian isotopies from $\Lambda_n$ to $\Lambda_{n+1}$
(``wrapping once''); it is not hard to check that the image of
$CF^0(\Lambda_n,\Lambda_0)$ inside $CW^0(\Lambda_0,\Lambda_0)$ is the span
of $e$ and $\theta_{i,k}$, $k\leq n$.

These Floer complexes describe morphisms in the Fukaya category 
$\F((\C^*)^2,x+y+z)$, which is equivalent to $D^b(\PP^2)$, with $\Lambda_0$
(resp.\ $\Lambda_n$) 
corresponding to $\O_{\bar{X}}$ (resp.\ $\O_{\bar{X}}(-3n)$), where
$\bar{X}=\{(x_1\!:\!x_2\!:\!x_3)\,|\,x_1+x_2+x_3=0\}\subset \PP^2$,
while the continuation map for wrapping once amounts to multiplication by the monomial
$x_1x_2x_3$ \cite{Hanlon}. The direct limit thus corresponds to rational
functions on $\bar{X}$ which are allowed to have arbitrary pole orders at
the three points where one of the homogeneous coordinates vanishes, i.e.\
regular functions on $X$, while the image of $HF^0(\Lambda_n,\Lambda_0)$ in
$HW^0(\Lambda_0,\Lambda_0)$ corresponds to rational functions with poles
of order at most $n$ at the punctures of $X$.  

In fact, $\theta_{i,k}$ is the image under continuation of
a generator of the Floer complex of $\Lambda_0$ with its image under
wrapping just the $i^\text{th}$ leg $k$ times. The continuation map for this
Hamiltonian isotopy amounts to multiplication by $x_i^k$ (again by
\cite{Hanlon}), and thus we conclude that $\theta_{i,k}$ corresponds to a
rational function which has only a pole of order at most $k$ at the
$i^\text{th}$ 
puncture of $X$; and the pole order has to be exactly $k$ since there is
no generator corresponding to $\theta_{i,k}$ when we wrap $k-1$ times.
\endproof

As a sanity check, we note that 
any collection of rational functions as in the lemma gives an additive basis of
$H^0(X,\O_X)$ (as follows e.g.\ from partial fraction decomposition).


The multiplicative structure on $HW^0(\Lambda_0,\Lambda_0)$ is surprisingly
difficult to calculate explicitly, and so is the Floer product 
\begin{equation}\label{eq:evalpantsatpoint}
\mu^2:HW^0(\Lambda_0,\Lambda_p)\otimes HW^0(\Lambda_0,\Lambda_0)\to
HW^0(\Lambda_0,\Lambda_p)
\end{equation}
where $\Lambda_p$ is a Lagrangian torus in $(\C^*)^2$ with a rank one
local system, corresponding to the skyscraper sheaf of a point $p\in X$,
and under our dictionary, to a circle $L_p$ inside the smooth part of $M$,
equipped with a rank one local system. The leading order terms
of these products, corresponding to the holomorphic discs with the lowest
geometric energy, can be determined readily; when considering generators
which lie within a single end, the projections
of these holomorphic discs from $(\C^*)^2$ onto the appropriate coordinate
axis in $M$ look precisely like the configurations depicted in 
Figure \ref{fig:pantsmirror}, and in fact they replicate the geometry of
wrapped Floer homology in (one half of) the infinite cylinder.

The geometric reason for this similarity is that, in the open subset
$U_x\subset (\C^*)^2$ where $|x|>\max(|y|,|z|)+C$ for a suitable constant $C>0$, we
can treat the geometry as the product of a factor $\C^*$ with coordinate $x$,
inside which $\Lambda_0$ corresponds to the real positive axis $\arg(x)=0$
while $\Lambda_p$ corresponds to a circle $|x|=$\,constant, and another
factor inside which $\Lambda_0$ and $\Lambda_p$ both correspond to the
circle $|y|=|z|$ (whose self-Floer homology is responsible for the presence
of generators in two different degrees, even though only degree 0 is of
interest to us). Thus, among the holomorphic discs contributing to the product
structure on $HW^0(\Lambda_0,\Lambda_0)$ and to \eqref{eq:evalpantsatpoint},
those which remain within $U_x$ can be determined explicitly, and agree with
the corresponding products in wrapped Floer cohomology for the real axis and a 
circle inside (one half of) $\C^*$. (Similarly for the two other ends of
$\Lambda_0$.)

If these discs were the only ones contributing to Floer products, then it
would follow that $\theta_{i,k}=(\theta_{i,1})^{k}$, so that
$\theta_{i,k}$ corresponds to the $k^\text{th}$ power of a rational function of 
degree~1 with a single pole at the $i^\text{th}$ puncture of $X$ (i.e., the 
inverse of a local coordinate $t_i$), and the product
\eqref{eq:evalpantsatpoint} corresponds to evaluation at a point $p$
where the value of the coordinate $t_i$ is directly
determined by the position of 
$\Lambda_p$ and the holonomy of its local system around the $x$ factor.
However, there is no obvious reason why every holomorphic disc contributing
to the Floer product should be entirely contained in $U_x$, even if its inputs
and output all lie near $|x|\to\infty$; for example, the Floer differential
on $CF^0(\Lambda_0,\Lambda_p)$ is known to involve not only holomorphic
discs within $U_x$ but also some whose
image under the logarithm map propagates all the way to the vertex of the 
tropical pants \cite{Hicks}. The model we construct in \S \ref{s:Amodel}
below ignores the contributions from any such discs, and instead
chooses the correspondence between the wrapped Floer cohomology of
$L_0$ in $M$ and the ring of functions of $X$ to be the simplest possible one,
{\em even though this means that the identification between $\End(L_0)$ and
$HW^0(\Lambda_0,\Lambda_0)$ may differ from the expected one by instanton
corrections.}

Additionally, there is a well-known ambiguity in the manner in which a local system
on a simple closed curve $L_p\subset M$ determines one on $\Lambda_p=
\partial\mathcal{T}(L_p)$ in such a way that the point $p$ lies on the pair of 
pants $X$. This is because equipping the thimble $\mathcal{T}(L_p)$
with a bounding cochain requires the choice of a splitting of the
map $H_1(\Lambda_p)\twoheadrightarrow H_1(\mathcal{T}(L_p))\simeq H_1(L_p)$;
in the literature on open Gromov-Witten theory this is called a 
{\em framing} for each leg of $M$. It is not hard to see that the choice of
framing amounts to a choice of local coordinate on $X$; the most natural
choices for each puncture are those given by ratios of homogeneous coordinates
on the compactification $\bar{X}\subset \PP^2$, which take the values $-1$ and
$\infty$ at the other two punctures (compare with Definition \ref{def:combidata}), but 
from a Floer-theoretic perspective there is no particular reason to
restrict oneself to these. In fact, considerations about equivariance with
respect to permuting the coordinates $(x,y,z)$ suggest that the zeroes of the
rational functions $t_i^{-1}$ associated to the generators $\theta_{i,1}$
are not at the punctures of $X$ but rather at the points with homogeneous
coordinates $(-\frac12\!:\!-\frac12\!:\!1)$, 
$(-\frac12\!:\!1\!:\!-\frac12)$, and $(1\!:\!-\frac12\!:\!-\frac12)$.

Regardless of the above issues, the most important unexpected feature of
wrapped Floer theory in $M$ that emerges from our geometric considerations
is that {\em holomorphic discs in $M$ must be allowed to propagate through
the vertex at the origin}. By using mirror symmetry and calculating the
product in the ring of functions $H^0(X,\O_X)$, the following is a direct
consequence of Lemma \ref{l:floergens_as_rationalfunctions}:

\begin{lemma}
For $i\neq j$ and $k,\ell\geq 1$, the Floer product $\mu^2(\theta_{i,k},
\theta_{j,\ell})\in HW^0(\Lambda_0,\Lambda_0)$ is a nontrivial linear combination of
the generators $e$, $\theta_{i,k'}$ $(k'\leq k)$ and $\theta_{j,\ell'}$
$(\ell'\leq \ell)$.
Moreover, for any given generator $\theta_{i,k}$, the Floer product 
\eqref{eq:evalpantsatpoint} is nonzero for all but finitely many tori with
local systems $\Lambda_p$ corresponding to skyscraper sheaves $\O_p$, $p\in X$.
\end{lemma}

Therefore, irrespective of the exact manner in which
we transcribe the wrapped Floer cohomologies $HW^0(\Lambda_0,\Lambda_0)$ and 
$HW^0(\Lambda_0,\Lambda_p)$ into Lagrangian Floer theory for $L_0$ and $L_p$
inside $M$ and the instanton corrections that may be packaged into this
dictionary, Floer products in $M$ must include not
only holomorphic discs which lie inside one of the three components of $M$,
but also nodal configurations of discs which lie in different components and
are attached to each other through the origin. That such a construction can
be carried out in a way that accurately reflects the geometry of 
homological mirror symmetry is {\em a priori} not clear; thus, instead of
relying on the above intuition, in Sections \ref{s:Amodel} and
\ref{s:equivalence} we describe our 
A-model construction from scratch, verify that its product operations 
satisfy the $A_\infty$-relations, and verify homological mirror symmetry.

\subsection{Beyond the pair of pants}

While our main focus is on mirrors of closed curves, our construction 
also applies in the punctured case, where the more usual approaches to
mirror symmetry stem from Hori-Vafa and Aganagic-Vafa's
construction of mirror curves from toric Calabi-Yau 3-folds \cite{HoriVafa,AV}.

The next simplest example after the pair of pants is the 4-punctured genus
0 curve $X=\PP^1\setminus \{0,1,q,\infty\}$. It was first shown by
Hori-Vafa \cite{HoriVafa} and Aganagic-Vafa \cite{AV} that the curve $X$ (or rather, a specific embedding
of $X$ as a very affine plane curve) arises out of mirror symmetry
for a toric Calabi-Yau 3-fold $Y,$ the so-called ``resolved conifold'',
i.e.\ the total space of the bundle $\O(-1)\oplus\O(-1)\to \CP^1$;
or more precisely, a toric Landau-Ginzburg model $(Y,W)$
consisting of $Y$ equipped with a suitable superpotential \cite{ChanLauLeung,AAK}.
Choosing a toric K\"ahler form on $Y$ for which the area of the zero section $\CP^1\subset Y$
is $A$, its Hori-Vafa mirror curve is 
$C=\{(x,y)\,|\,1+x+y+qxy=0\}$ in $(\C^*)^2$ or,
working over the Novikov field, $(K^*)^2$; here $q=T^A$.

The critical locus $M=\mathrm{crit}(W)$ is the union of the 1-dimensional toric strata of $Y$, i.e.\
the zero section $\CP^1$ and four copies of $\C$ (the fibers of
the two distinguished line subbundles over $0$ and $\infty$). 
This configuration is decribed by a graph $G$ with two trivalent vertices
$v,v'$ (corresponding to the two triple points of $M$) connected by a single
edge $e$ corresponding to the $\CP^1$ component of $M$, while the other
(external) edges $e_1,\dots,e_4$ correspond to the four $\C$ components.
The construction in Section \ref{s:Bmodel} then produces $X$
out of two punctured rational curves $W^0_v,W^0_{v'}$, equipped with
coordinates $t_{e/v}$, $t_{e/v'}$ identifying them with 
$\PP^1\setminus \{0,1,\infty\}$, by setting
$t_{e/v}t_{e/v'}=q_e=T^A$.

In this case, Hori-Vafa mirror symmetry and our construction agree (with
the same mirror parameter $q=q_e=T^A$), as we can identify $X$ with $C$ by
considering the distinguished coordinates $t_{e/v}=-x^{-1}$ on $
X\simeq C\simeq \PP^1\setminus
\{0,1,q,\infty\}$. However, the Hori-Vafa curve $C$ arises naturally 
as a plane curve, with distinguished ambient coordinates $x,y$ (up to
monomial coordinate changes) arising from the toric geometry of $Y$; while
our mirror curve $X$ does not come with an embedding, though it does have a
preferred collection of coordinates after combinatorial data is chosen for
the external edges as
in Definition \ref{def:combidata}. The most natural choice is to have
the coordinate for each
external edge take the value $\infty$ at the puncture corresponding
to the edge $e$ and $1$ at the other puncture, which yields: for the
vertex $v$, $t_{e/v}=-x^{-1}$, $t_{e_1/v}=-x$, and $t_{e_2/v}=1+x$;
for the vertex $v'$, $t_{e/v'}=-qx$, $t_{e_3/v'}=-(qx)^{-1}$, $t_{e_4/v'}=
1+(qx)^{-1}$. In toric mirror symmetry, the most natural
coordinate on the end of $C$ where $y\to 0$ is $y^{-1}$ (or $y^{-1}$ times
some power of $x$); while in our setting it is
$t_{e_2/v}^{-1}=(1+x)^{-1}=(q-1)^{-1} (y^{-1}+q)$.

This discrepancy is due to the different manners in which
the mirror curve is built from the local pairs of pants in the two approaches.
The Hori-Vafa mirror curve $C$ is assembled by
patching together the pairs of pants $1+x+y=0$ (corresponding to the vertex
$v$) and $q+x^{-1}+y^{-1}=0$ (or equivalently $qxy+x+y=0$, corresponding to the vertex $v'$), combining
the various monomials into a single Laurent polynomial (here, $1+x+y+qxy$).
There is no preferred mapping from the pairs of pants to the assembled
curve but, if we choose to match the pieces to each other via the $x$ 
coordinate, then the gluing of the two pairs of pants involves a correction
of the $y$ coordinate, which can be understood in terms of the enumerative
geometry of ``outer Aganagic-Vafa branes'' on the resolved conifold
\cite{AV}. By contrast, in our construction the gluing between pairs of
pants to build $X$ does not involve a deformation of the 
local coordinates $t_{e_i/v}$ (meaning, a posteriori, that the expression
for $t_{e_2/v}$ in terms of $y$ already incorporates the needed
corrections).

The differences between the two approaches become starker if we consider
higher genus examples. Consider e.g.\ the case where $Y$
is the total space of $\O(-3)\to\CP^2$. Choosing a toric K\"ahler form for which
the symplectic area of $\CP^1$ is $A$, the Hori-Vafa mirror curve is
$C=\{(x,y)\,|\,f(q)+x+y+q/xy=0\}$ in $(K^*)^2$, where
$q=T^A$ and
$f(q)=1-2q+5q^2-32q^3+\dots$ is a generating series for certain open Gromov-Witten
invariants in $Y$ (see e.g.\ \cite[\S 5.3.3]{ChanLauLeung} and references therein). 
The thrice-punctured elliptic curve $C$ is isomorphic to
$X=(K^*- q_e^\Z)/q_e^{3\Z},$ where $q_e=-q(1-9q+108q^2-1461q^3+\dots)$.

In our setup, $X$ can be obtained by gluing together three pairs of pants
$W^0_{v_i}$ ($i\in \Z/3$) with coordinates
$t_{e_i/v_i}=t_i\in \PP^1\setminus \{0,1,\infty\}$, $t_{e_{i+1}/v_i}=t_i^{-1}$, and $t_{e'_i/v_i}=1-t_i$,
via $t_{e_{i+1}/v_i}t_{e_{i+1}/v_{i+1}}=q_e$ for all $i\in \Z/3$,
or equivalently, $t_{i+1}=q_et_i$. Thus, our construction produces $X$ as
mirror curve if we take $M$ to be the union of the 1-dimensional strata of
$Y$, assigning the corrected K\"ahler parameter $q_e$ to each of the three
$\PP^1$ components and using the above sets of coordinates on the pair of
pants as combinatorial data for the three triple points.

We note that the functions $t_{e_i/v_i}$ do not determine actual coordinates on
the punctured elliptic curve $X$, but rather on its $\Z$-cover
$\tilde{X}=K^*- q_e^\Z$. 
(This is related to the
construction of our curves as quotients of punctured $\PP^1$'s via Schottky
uniformization, cf.\ Section \ref{s:Bmodel}.)
Thus, they behave very differently from the ambient coordinates $(x,y)$ for the
Hori-Vafa curve $C$. A posteriori, those can be recovered from the coordinate $t$ on
$\tilde{X}=K^*-q_e^\Z$ as ratios of suitable theta functions:
$$x=\frac{\vartheta_{-\frac16,\frac12}(q_e^3,t)^2}
{\vartheta_{\frac16,\frac12}(q_e^3,t)\vartheta_{\frac12,\frac12}(q_e^3,t)}=
\frac{-q_e^{-1/3}t^{-1}\Bigl(\sum_{n\in\Z}(-1)^n q_e^{(3n^2-n)/2}t^n\Bigr)^2}
{\Bigl(\sum_{n\in\Z}(-1)^n q_e^{(3n^2+n)/2}t^n\Bigr)
\Bigl(\sum_{n\in\Z}(-1)^n q_e^{(3n^2+3n)/2}t^n\Bigr)}$$
and
$$y=\frac{\vartheta_{\frac16,\frac12}(q_e^3,t)^2}
{\vartheta_{-\frac16,\frac12}(q_e^3,t)\vartheta_{\frac12,\frac12}(q_e^3,t)}=
\frac{q_e^{-1/3}\Bigl(\sum_{n\in\Z}(-1)^n q_e^{(3n^2+n)/2}t^n\Bigr)^2}
{\Bigl(\sum_{n\in\Z}(-1)^n q_e^{(3n^2-n)/2}t^n\Bigr)
\Bigl(\sum_{n\in\Z}(-1)^n q_e^{(3n^2+3n)/2}t^n\Bigr)}$$
are invariant under $t\mapsto q_e^{3} t$ and satisfy $$x+y+\frac{1}{xy}=q_e^{-1/3}(1+5q_e-7q_e^2+3q_e^3+\dots)=
-q^{-1/3}(1-2q+5q^2-32q^3+\dots),$$ which is the same as the Hori-Vafa curve
$C$ after rescaling $x$ and $y$ by a factor of $q^{1/3}$.
We do not know whether the coefficients in the
Laurent series expansions of these formulas for $x$ and $y$
should be expected to admit an enumerative/combinatorial interpretation,
nor even whether the above choices of parameters for our $A$-model on $M$ are
the most natural ones in this respect.

\section{The A-model: Lagrangian Floer theory in trivalent configurations} \label{s:Amodel}

\subsection{Objects and morphisms in $\F(M)$}

Let $G$ be a graph with finite set of vertices $V$ and edges $E$, such that
each vertex $v\in V$ has degree 3. As noted in Remark \ref{rmk:opencase}, we
allow ``external edges'' which only connect to one vertex.
We denote by $E^0$ the
set of external edges and by $E^i$ the internal edges. We also fix for
each internal edge $e\in E^i$ an area parameter $A_e>0$, and an element $q_e\in K$
with $\val(q_e)=A_e$ (we will mostly focus on the case $q_e=T^{A_e}$); for external edges there are no 
area parameters but we consider either {\em wrapped} or {\em stopped}
Lagrangian Floer theory.

For each internal edge $e\in E^i$,
we consider $M_e=S^2=\CP^1$, equipped with a symplectic form $\omega$ of
total area $A_e$ (eg.\ a multiple of the standard symplectic form), and optionally a bulk deformation class $\mathfrak{b}\in H^2(M_e,\cO_K)$ such
that $T^{A_e} \exp(\int_{M_e} \mathfrak{b})=q_e$. We also fix two marked points on
$M_e$, which we think of as $0$ and $\infty$ in $\CP^1$,
and assign them to the vertices $v,v'\in V$ joined by the edge $e$: 
$\{p_{e/v},p_{e/v'}\}=\{0,\infty\}\subset M_e$.  
For each external edge $e\in E^0$, we set $M_e=\C$, with the standard
symplectic form (of infinite area) and a single marked point $p_{e/v}=0\in M_e$.

Let $M$ be the space obtained by attaching the surfaces $M_e$, $e\in E$ to
each other at the triples of marked points which correspond to the same vertex 
of the graph $G$:
$$M=\textstyle \Bigl(\bigsqcup\limits_{e\in E} M_e\Bigr)/(p_{e/v}\sim p_{e'/v}\sim p_{e''/v}\ \forall v\in V).$$
We denote by $p_v$ the resulting nodal point of $M$. This gluing is purely cosmetic, as the
actual symplectic geometry will take place on the individual components
$M_e$. On the other hand, one important piece of data associated to each
vertex $v\in V$ is that of local
coordinates $t_{e/v}$ on the abstract curve $X^0_v = \PP^1$ which vanish
at the respective marked points $x_{e/v}\in X^0_v$ (cf.\ Definition
\ref{def:combidata}).

We fix an asymptotic direction near $0$ and $\infty$ on each component 
$M_e\subset M$, for example the real positive axis; all
Lagrangians we consider will be required to approach the nodes of $M$ and
its infinite ends along this prescribed direction. 

\begin{defi} The {\em objects} of $\F(M)$ are pairs $(L,\cE)$, where
$L\subset M$ is a properly embedded (trivalent) graph whose vertices lie at the
nodes of $M$ and whose edges lie in the smooth part of $M$, in such a way
that:
\begin{itemize} \item the arc components of $L_e=L\cap M_e$ approach $0$
and $\infty$ along the prescribed asymptotic directions; 
\item the closed curve components of $L_e$ are homotopically non-trivial
in the complement of the marked points;
\item a node $p_v\in M$ lies on $L$ if and only if it is an end point of an arc in
each of the three components of $M$ which meet at $p_v$;
\end{itemize}
and $\cE$ is a unitary local system, i.e.\ a local system of free
finite rank $\cO_K$-modules over $L$.
\end{defi}

Because each component of $M$ is either $(\CP^1,\{0,\infty\})$ or
$(\C,\{0\})$, this definition only allows for two types of indecomposable objects.
\begin{enumerate}
\item {\em Point-type objects:} $L$ is a 
simple closed curve in the smooth part of a component $M_e$, separating $0$
from $\infty$. When $\cE$ has rank 1, the object $(L,\cE)$ corresponds under
mirror symmetry to the skyscraper sheaf of a point of $X_K$ where the
valuation of the coordinate $t_{e/v}$ equals the symplectic area enclosed
by $L$ around the marked point $p_{e/v}$.
\item {\em Vector bundle (v.b.) type objects:} 
$L$ is a trivalent graph with the same sets of edges and vertices as $G$,
consisting of an arc $L_e$ connecting $0$ to $\infty$ in each component $M_e$,
and passing through all the nodes. When $\cE$ has rank 1, the object
$(L,\cE)$ corresponds to a line bundle over $X_K$, as described in \S
\ref{ss:assignbundles} below.
\end{enumerate}

We also specify a class of
smooth Hamiltonian perturbations to be used for defining Floer complexes
between objects of $\F(M)$.

\begin{defi}
A {\em positive Hamiltonian} is a smooth function $h:M\to\R$ which, on each compact component $M_e\simeq
\CP^1$, $e\in E^i$, has local minima at the two marked points $0$ and $\infty$,
$h(0)=h(\infty)=0$, and on each non-compact component $M_e\simeq \C$, $e\in
E^0$, has a minimum at the origin $h(0)=0$, and linear, resp.\
quadratic growth at infinity (in terms of the coordinate $r=|z|^2$) when the 
non-compact end does, resp.\ doesn't carry a stop.
\end{defi}

The flow of such a
Hamiltonian rotates the asymptotic directions near the
marked points in the positive direction, and pushes or wraps the infinite
ends in the customary manner for (partially) wrapped Floer theory.

For each pair $(L,L')$ we choose a positive Hamiltonian $h$ and a small
$\varepsilon>0$ such that
$L^+=\phi^1_{\varepsilon h}(L)$ is transverse to $L'$, and define the generators
of the Floer complex to be time~1 trajectories of the Hamiltonian vector field
generated by $\varepsilon h$ which start at $L$ and end at $L'$,
\begin{equation}\label{eq:floergenerators}
\cX(L,L')=\{\gamma:[0,1]\to M\,|\,\gamma(0)\in L,\ \gamma(1)\in L',\ \dot\gamma(t)=X_{\varepsilon h}(\gamma(t))\}
\end{equation}
or equivalently, pairs of points in $L$ and $L'$ which match under the flow:
$$\cX(L,L')\simeq
\{(p,p')\in L\times L'\,|\,\phi^1_{\varepsilon h}(p)=p'\},$$
or even simpler, points of $L^+\cap L'$. (Abusing notation we think of elements of $\cX(L,L')$ interchangeably
as points, pairs of points, or trajectories of $X_{\varepsilon h}$.)
Note that, when $L$ and $L'$ are of vector bundle type, $\cX(L,L')$ always
includes one generator at each node of $M.$ 

We define morphism spaces by
\begin{equation} \label{eq:morphisms}
\hom_{\F(M)}((L,\cE),(L',\cE'))=CF^*((L,\cE),(L',\cE');\varepsilon h)=\bigoplus_{(p,p')\in \cX(L,L')} \cE_p^*\otimes
\cE'_{p'}
\end{equation}
(Another option would be to define $\F(M)$ by considering a directed category whose
objects are images of $(L,\cE)$ under positive Hamiltonian flows, and
localizing with respect to continuation elements $e_{(L,\cE),\varepsilon}\in
CF^*(\phi^1_{\varepsilon h}(L,\cE),(L,\cE))$; while this is more consistent
with some of the recent literature \cite{AA,AS}, there is no benefit
to doing so in our setting.)

The choice of a
trivialization of the tangent bundle $TM$ outside of the nodes determines
a $\Z$-grading on $\F(M)$;
the preferred choice in our case is the trivialization determined by the radial line
field on the open stratum $\C^*\subset M_e$ of each component.  Objects
should then be graded by choosing a
real-valued lift of the angle between $TL$ and the chosen line field outside
of the nodes. Here again there is a preferred choice: for
v.b.-type objects we declare the angle between $TL$ and the outward radial
line field to be zero near both ends (at $0$ and $\infty$) in each
component, and for point-type objects where $L$ is a circle centered at the
origin in $M_e$ we declare the angle to be $-\pi/2$. With this convention,
all Floer cohomology groups are concentrated in degrees 0 and 1, and for
pairs of v.b.-type objects the
generators which lie at the nodes of $M$ are in degree 0.

\begin{remark}
Because of the positive Hamiltonian perturbations involved in defining
morphism spaces, when $M$ is compact (all components are $\PP^1$\!'s) the category $\F(M)$ is never Calabi-Yau.
The study of open-closed and closed-open maps for $\F(M)$
is beyond the scope of this paper, but we point out that the
Hochschild cohomology of $\F(M)$ is expected to be isomorphic to the
fixed point Floer cohomology of a small positive Hamiltonian, via the
closed-open map $$\mathcal{CO}:HF^*(\phi^1_{\varepsilon h})\to HH^*(\F(M)).$$
For instance, when $M$ consists of $3g-3$
$\PP^1$'s meeting in $2g-2$ triple points, there is a positive Hamiltonian
with $2g-2$ minima (at the nodes), $3g-3$ saddle points, and $3g-3$
maxima. The Floer differential on $CF^*(\phi^1_{\varepsilon h})$
agrees with the Morse
differential within each component of $M$, so each minimum maps to the
sum of three saddle points, and
$$\dim HF^0(\phi^1_{\varepsilon h})=1,\quad 
\dim HF^1(\phi^1_{\varepsilon h})=g,\quad \text{and}\ 
\dim HF^2(\phi^1_{\varepsilon h})=3g-3,$$
in agreement with the Hochschild cohomology of 
the derived category of a genus $g$ curve.
\end{remark}

\subsection{$A_\infty$-operations: propagating discs}

The $A_\infty$-operations in $\F(M)$ are determined by weighted counts of
``propagating'' configurations of (perturbed) holomorphic discs for some
choice of complex structure $J$ on $M$ (the choice is 
immaterial). To define $$\mu^k:\hom((L_{k-1},\cE_{k-1}),(L_k,\cE_k))\otimes
\dots\otimes \hom((L_0,\cE_0),(L_1,\cE_1))\to 
\hom((L_0,\cE_0),(L_k,\cE_k))[2-k]$$ we consider maps whose
domain $S$ is a nodal union of discs, modelled on a planar rooted 
tree $T$ with $k+1$ external edges (one root and $k$ leaves).
For each internal vertex $v_j$ of $T$ we consider a disc $D_j$ with $|v_j|$ boundary marked points,
and define $S=\bigsqcup D_j/\sim$, where for each internal edge of $T$
connecting vertices $v_j,v_{j'}$, we glue $D_j$ to $D_{j'}$ by identifying 
the two boundary marked points that correspond to the end points of the edge. 
The resulting nodal configuration still carries $k+1$ marked
points $z_0$ (corresponding to the root of $T$), $z_1,\dots,z_k$ (corresponding to the leaves), 
in that order along the boundary of $S$. We label each portion of
$\partial S$ from $z_i$ to $z_{i+1}$ (or $z_k$ to $z_0$, for $i=k$) by 
the Lagrangian $L_i$. 
Orienting the tree $T$ from the leaves to the root, 
each component of $S$
has one output marked point (towards the root) and one or more input marked
points (towards the leaves). We choose strip-like ends near each of these,
i.e.\ local coordinates $s+it$ such that the input ends are modelled on 
$\R_+\times [0,1]$ and the output end on $\R_-\times [0,1]$. 
We also choose a 1-form $\beta$ on $S$, such that $\beta_{|\partial S}=0$
and $\beta$ is a small positive multiple of $dt$ on each strip-like end.

\begin{defi}\label{def:proptree}
Given $L_0,\dots,L_k$, generators $p_i\in \cX(L_{i-1},L_i)$ for $1\leq i\leq k$
and $p_0\in \cX(L_0,L_k)$, and a planar tree $T$\!,
a {\bf propagating holomorphic disc} modelled on $T$ is a map
$u:(S,\partial S)\to (M,L_0\cup\dots\cup L_k)$, where the domain $S$
is modelled on $T$, such that 
\begin{enumerate}
\item each component of $S$ maps to a single
component of $M$;
\item $u$ satisfies the perturbed Cauchy-Riemann equation
\begin{equation}\label{eq:cauchyriemann}
(du-X_h\otimes \beta)^{0,1}=0
\end{equation} 
on each component of $S$, where
$h$ is the positive Hamiltonian used to define morphism spaces, and $\beta$ is the chosen
1-form on $S$;
\item the nodes of $S$ map to nodes of $M$;
\item the map $u$ converges at each input marked point $z_i$, resp.\ the output $z_0$, to the flowline
of $X_{\varepsilon h}$ which defines the generator $p_i\in \cX(L_{i-1},L_i)$, resp.\
$p_0\in \cX(L_0,L_k)$;
\item the components of $u$ are not allowed to pass through the nodes of
$M$ except at the nodes of $S$, at input marked points 
$z_i\in S$, or at a constant component carrying the output marked point $z_0\in S$;
\item  when an input marked point $z_i\in S$ maps to a node of $M$,
the restriction of $u$ to the strip-like end near $z_i$
 does not surject onto a neighborhood of the node in the appropriate component of $M$;
\item if the output marked point $z_0\in S$ maps to a node of $M$ then the restriction
of $u$ to the component of $S$ carrying $z_0$ is a constant map.
\end{enumerate}
The moduli space of such propagating discs $u$ 
in a fixed homotopy class $[u]$, up to reparametrization, 
is denoted by $\cM(p_0,\dots,p_k,[u])$.
\end{defi}

\begin{figure}[t]
\begin{tikzpicture}[scale = 0.6]
\begin{scope}[shift={(-10,1.3)},scale=1.3]
\fill[orange!25!white] (-1.2,-0.2) circle (0.65);
\fill[orange!25!white] (1.2,-0.2) circle (0.65);
\fill[orange!25!white] (0,0.3) circle (0.65);
\draw[green!50!gray,semithick](1.85,-0.2) arc (0:157.4:0.65) arc (-22.6:90:0.65);
\draw[red,semithick](0,0.95) arc (90:202.6:0.65) arc (22.6:180:0.65);
\draw[blue,semithick](-1.85,-0.2) arc (-180:22.6:0.65) arc (202.6:337.4:0.65) arc (157.4:360:0.65);
\fill(0,0.95) circle (0.08) (0,1.2) node {\SMALL $z_1$};
\fill(-1.85,-0.2) circle (0.08) (-2.15,-0.2) node {\SMALL $z_2$};
\fill(1.85,-0.2) circle (0.08) (2.15,-0.2) node {\SMALL $z_0$};
\draw[red](-1.2,1) node {\SMALL $L_1$};
\draw[green!50!gray](1.2,1) node {\SMALL $L_0$};
\draw[blue](0,-0.9) node {\SMALL $L_2$};
\draw[semithick,->] (3,-0.4)--(4.5,-0.4);
\draw(-2.5,1) node {$S$};
\draw(-2.5,-2.5) node {$T$};
\draw[thick] (-0.8,-3) circle (0.08) (0,-2.5) circle (0.08) (1,-2.5) circle (0.08);
\draw[thick,-latex] (-1.5,-3)--(-0.85,-3);
\draw[thick,-latex] (-0.72,-3)--(-0.05,-2.5);
\draw[thick,-latex] (-0.7,-2.1)--(-0.05,-2.5);
\draw[thick,-latex] (0.08,-2.5)--(0.95,-2.5);
\draw[thick,-latex] (1.08,-2.5)--(1.8,-2.5);
\end{scope}
\fill[orange!25!white] (0,0) arc (180:360:3 and 0.2) arc (0:103.5:3 and 1) arc (90:160:0.2 and 0.3) arc (345:270:1.2 and 1.02);
\fill[orange!25!white] (2.95,-0.98) arc (-90:-26:0.35 and 1.17) arc (-108:-90:5) arc (-90:-45:1.77 and 2) arc (0:-90:3 and 1);
\draw[green!50!gray,thick](0,0) arc (110:65:3.5) arc (-115:-90:5) arc (-90:-45:1.77 and 2);
\draw[red,thick](0,0)--(1,0) arc (270:345:1.2 and 1.02) arc (160:90:0.2 and 0.3);
\draw[red,thick,dotted] (2.4,0.92) arc (90:10:0.25 and 1.1) arc (190:270:0.25 and 1.15);
\draw[red,thick] (2.95,-0.98) arc (-90:-10:0.35 and 1.17) arc (170:90:0.35 and 1.15);
\draw[red,thick,dotted] (3.7,0.9) arc (90:10:0.25 and 1.1) arc (190:270:0.25 and 1.1);
\draw[red,thick] (4.15,-0.89) arc (-90:-30:0.2 and 0.3) arc (170:100:1.5 and 0.6)
arc (-70:-20:0.73 and 0.33);
\draw[blue,thick](0,0) arc (180:360:3 and 0.2);
\draw[semithick] (3,0) ellipse (3 and 1);
\fill(3.25,-0.35) circle (0.1) (2.95,-0.55) node {\SMALL $p_1$};
\draw[green!50!gray](5.35,-1) node {\SMALL $L_0^{+2}$};
\draw[blue](1.8,-0.5) node {\SMALL $L_2$};
\draw[red](2.5,1.35) node {\SMALL $L_1^+$};
\begin{scope}[rotate=140,shift={(0.05,-0.05)}]
\draw[green!50!gray,thick](0,0) arc (110:65:3.5) arc (-115:-101.5:5);
\draw[red,thick](0,0)--(1,0) arc (270:345:1.2 and 0.98) arc (160:90:0.2 and 0.26);
\draw[red,thick,dotted] (2.4,0.88) arc (90:10:0.25 and 1.06) arc (190:270:0.25 and 1.12);
\draw[red,thick] (2.95,-0.96) arc (-90:-10:0.35 and 1.17) arc (170:90:0.35 and 1.15);
\draw[blue,thick](0,0) arc (180:270:3.8 and 0.2);
\draw[semithick] (4,1) arc (90:270:4 and 1);
\draw[semithick] (4,0) ellipse (0.2 and 1);
\end{scope}
\begin{scope}[rotate=-140,shift={(0.05,0.05)}]
\fill[orange!25!white](0,0) arc (180:254.5:4 and 1) arc (-90:-20:0.35 and 1.17) arc (260:180:3.8 and 0.2);
\fill[orange!25!white](0,0) arc (180:114:4 and 1) arc (90:160:0.2 and 0.26) arc (345:270:1.2 and 0.98);
\draw[green!50!gray,thick](0,0) arc (110:65:3.5) arc (-115:-101.5:5);
\draw[red,thick](0,0)--(1,0) arc (270:345:1.2 and 0.98) arc (160:90:0.2 and 0.26);
\draw[red,thick,dotted] (2.4,0.88) arc (90:10:0.25 and 1.06) arc (190:270:0.25 and 1.12);
\draw[red,thick] (2.95,-0.96) arc (-90:-10:0.35 and 1.17) arc (170:90:0.35 and 1.15);
\draw[blue,thick](0,0) arc (180:270:3.8 and 0.2);
\draw[semithick] (4,1) arc (90:270:4 and 1);
\draw[semithick] (4,0) ellipse (0.2 and 1);
\fill(3.3,-0.18) circle (0.1) (2.95,0.15) node {\SMALL $p_2$};
\end{scope}
\begin{scope}[shift={(6.02,0.07)},rotate=40]
\fill[orange!25!white](0,0) arc (110:63:3.5) arc (255:180:3.8 and 0.2);
\draw[green!50!gray,thick](0,0) arc (110:65:3.5) arc (-115:-101.5:5);
\draw[red,thick](0,0)--(1,0) arc (270:345:1.2 and 0.98) arc (160:90:0.2 and 0.26);
\draw[red,thick,dotted] (2.4,0.88) arc (90:10:0.25 and 1.06) arc (190:270:0.25 and 1.12);
\draw[red,thick] (2.95,-0.96) arc (-90:-10:0.35 and 1.17) arc (170:90:0.35 and 1.15);
\draw[blue,thick](0,0) arc (180:270:3.8 and 0.2);
\draw[semithick] (4,1) arc (90:270:4 and 1);
\draw[semithick] (4,0) ellipse (0.2 and 1);
\fill(2.85,-0.18) circle (0.1) (3,0.1) node {\SMALL $p_0$};
\end{scope}
\begin{scope}[shift={(6.02,-0.07)},rotate=-40]
\draw[green!50!gray,thick](0,0) arc (110:65:3.5) arc (-115:-101.5:5);
\draw[red,thick](0,0)--(1,0) arc (270:345:1.2 and 0.98) arc (160:90:0.2 and 0.26);
\draw[red,thick,dotted] (2.4,0.88) arc (90:10:0.25 and 1.06) arc (190:270:0.25 and 1.12);
\draw[red,thick] (2.95,-0.96) arc (-90:-10:0.35 and 1.17) arc (170:90:0.35 and 1.15);
\draw[blue,thick](0,0) arc (180:270:3.8 and 0.2);
\draw[semithick] (4,1) arc (90:270:4 and 1);
\draw[semithick] (4,0) ellipse (0.2 and 1);
\end{scope}
\end{tikzpicture}
\caption{A propagating disc contributing to
the Floer product $\mu^2$}\label{fig:proptree_ex}
\end{figure}

\noindent (The gluing behavior and consistency needed to establish the $A_\infty$-relations are most easily proved if
$\beta=\varepsilon\,dt$ at all strip-like ends, however this may not be possible
on the non-compact components of $M$, where energy estimates require $d\beta\leq 0$; 
the easiest way around this is to use 
Abouzaid's rescaling trick \cite{AbGen}. Another approach, which we shall not pursue, would be to consider Floer complexes
constructed using arbitrary small multiples of the positive Hamiltonian $h$
and localize at quasi-isomorphisms induced by continuation.)

By a standard trick, when the 1-form $\beta$ is closed
we can recast perturbed holomorphic curves $u:S\to M$ (solutions of \eqref{eq:cauchyriemann}) with boundary on
$L_0,\dots,L_k$ as genuine holomorphic curves $v:S\to M$ (solutions of $(dv)^{0,1}=0$ for a suitable, possibly domain-dependent complex
structure) with boundary on $L_0^{+k}=\phi_{\varepsilon h}^k(L_0),\dots,L_{k-1}^+=\phi^1_{\varepsilon
h}(L_{k-1}), L_k$, by setting $v(z)=\phi_h^{\tau(z)}(u(z))$, where $\tau:S\to\R$ satisfies $d\tau=-\beta$. 
The holomorphic curves $v:S\to M$ are easier to visualize and enumerate, as they are simply
polygons drawn on $M$, so we always use this viewpoint for graphical representations,
as in Figure \ref{fig:proptree_ex}.

The operations $\mu^k$ count {\em rigid} propagating holomorphic discs,
i.e., those which occur in zero-dimensional moduli spaces. This happens
precisely when each component taken separately is rigid, i.e.\ an immersed
polygon with locally convex corners. (For a constant component carrying the
output marked point $z_0$
and mapping to a node of $M$, rigidity amounts to the component having 
exactly two inputs). Rigidity implies that the degrees
of the Floer generators satisfy $\deg(p_0)=\sum \deg(p_i)+2-k$.
Each rigid propagating disc contributing to $\mu^k$ is counted with a weight, 
which is determined by multiplying
several quantities associated to the homotopy class $[u]$: area and holonomy
weights of the components of $u$, as is customary when defining Fukaya categories
over Novikov fields, as well as {\em propagation coefficients} at the nodes of $S$, which 
are unique to our setting.

Consider a node $z_\bullet\in S$, at which the output vertex of a component $D_{in}$ 
is attached to an input vertex of another component $D_{out}$ (recall that we
orient the tree $T$ from the inputs of the operation, i.e.\ the leaves, to the 
output, i.e.\ the root). Under $u:S\to M$, $z_\bullet$ maps to a
node $p_v\in M$ corresponding to some vertex $v$ of the graph $G$, where
the components $M_{e_{in}}$ and $M_{e_{out}}$ which contain
$u(D_{in})$ and $u(D_{out})$ are attached to each other; here $e_{in}$ and
$e_{out}$ are two of the three edges of $G$ which meet at the vertex $v$.
Because the Lagrangian graphs in $M$ which serve as boundary conditions
for $u$ on $D_{in}$ and $D_{out}$ approach the node $p_v$ from fixed directions,
the restrictions of $u$ to the strip-like ends of $D_{in}$ and $D_{out}$
near $z_\bullet$ have well-defined integer {\em degrees} $k_{in}$ and $k_{out}$,
namely the total multiplicities with which the images of the strip-like ends
cover neighborhoods of $p_v$ inside $M_{e_{in}}$ and $M_{e_{out}}$.
For example, the two nodes of the configuration in Figure \ref{fig:proptree_ex}
both have $k_{in}=1$ and $k_{out}=0$. In general, because our Hamiltonian perturbation is
a small positive multiple of $h$, with a local minimum at the node,
for non-constant maps we always have $k_{in}\geq 1$ and $k_{out}\geq 0$.

Recall that the combinatorial data of Definition \ref{def:combidata} includes the choice
of coordinate functions $t_{e/v}$ vanishing at the points $x_{e/v}\in X^0_v\simeq \PP^1$ 
for each of the three edges $e/v$ in the graph $G$. The function $t_{e_{in}/v}^{-k_{in}}$,
with a pole of order $k_{in}$ at $x_{e_{in}/v}$, can be expanded as a power series
in $t_{e_{out}/v}$ in a neighborhood of $x_{e_{out}/v}$.

\begin{defi}\label{defi:propcoeff} 
For given edges $e_{in}/v$, $e_{out}/v$ and degrees $k_{in}\geq 1$, $k_{out}\geq 0$, 
we define the {\bf propagation coefficient}
$C^{v;e_{in},e_{out}}_{k_{in},k_{out}}$ to be the coefficient of
$t_{e_{out}/v}^{k_{out}}$ in the expansion of $t_{e_{in}/v}^{-k_{in}}$ as a
power series in $t_{e_{out}/v}$.
Given a rigid propagating disc $u:S\to M$ whose output does not lie at a node of $M$,
the {\bf propagation multiplicity} $\Pi C([u])$ is defined to be the product
of the propagation coefficients $C^{v;e_{in},e_{out}}_{k_{in},k_{out}}$
at all the nodes of $S$.
\end{defi}

\begin{example}
Recall our preferred choices of coordinates on $X^0_v=\PP^1$ are those which
take values $0,1,\infty$ at the three marked points: for example one might
take $t_0=z$, $t_1=(z-1)/z$, $t_\infty=(1-z)^{-1}$ as coordinates near
the marked points $0$, $1$ and $\infty$. In this case, $t_0^{-1}=1-t_1=-(t_\infty+t_\infty^2+\dots)$, 
and similarly for the other pairs of coordinates,
so the propagation coefficients are
$$C^{v;e_{in},e_{out}}_{k_{in},k_{out}}=\begin{cases}
\displaystyle (-1)^{k_{out}}\binom{k_{in}}{k_{out}} & \text{for }(x_{e_{in}/v},x_{e_{out}/v})\in\{(0,1),(1,\infty),(\infty,0)\},
\\[10pt]
\displaystyle (-1)^{k_{in}}\binom{k_{out}-1}{k_{in}-1} & \text{for }(x_{e_{in}/v},x_{e_{out}/v})\in\{(0,\infty),(1,0),(\infty,1)\}.
\end{cases}$$
\end{example}
\smallskip

\noindent {\bf Output mapping to a node.}
The case where the output marked point $z_0\in S$ maps to a node $p_v\in M$ has a different flavor.
Recall that the whole component $D_0$ of $S$ carrying $z_0$ is required to map to $p_v$, and
rigidity implies that $D_0$ carries exactly two inputs.
If an input of $D_0$ is a node of $S$, we denote by $e_i$ the edge of $G$ such that the component of $S$ attached to $D_0$
at this node maps to $M_{e_i}$, and by $k_i\geq 1$ the degree of its output strip-like end (the incoming degree into the node),
and we associate to it the function $t_{e_i/v}^{-k_i}$ on $X^0_v\simeq \PP^1$.
If an input of $D_0$ is an input marked point of $S$, we instead consider the constant function 1 (this amounts to setting $k_i=0$).
The contribution of the nodes adjacent to the constant component $D_0$ to the
propagation multiplicity is then defined to be the constant term in the
expression of $t_{e_1/v}^{-k_1}t_{e_2/v}^{-k_2}$ as a linear combination
of $\{1, t_{e_1/v}^{-j}, t_{e_2/v}^{-j}\,|\, j\geq 1\}$. We denote this
coefficient by $K^{v;e_1,e_2}_{k_1,k_2}$. (Of note, this can only be nonzero when either
$S=D_0$, for a constant curve contributing to $\mu^2$, or both
inputs of $D_0$ are nodes of $S$ and $e_1\neq e_2$). The propagation 
multiplicity $\Pi C([u])$ is then defined to be the product of
$K^{v;e_1,e_2}_{k_1,k_2}$ (for the two nodes adjacent to the constant component
$D_0$) and the propagation
coefficients $C^{v;e_{in},e_{out}}_{k_{in},k_{out}}$ at all the other nodes
of $S$.

\begin{defi} The {\bf area weight} of a 
propagating holomorphic disc $u:S\to M$ with boundary on $L_0,\dots,L_k$,
inputs $(p_i,p'_i)\in \cX(L_{i-1},L_i)$ and output $(p_0,p'_0)\in \cX(L_0,L_k)$
is $$W([u]):=T^{A([u])}\,\int_S u^*\mathfrak{b}\in K, \quad \text{where}\ A([u])=\int_S u^*\omega.$$
When the $L_i$ are equipped with local systems $\cE_i$, the {\bf holonomy weight} of $u$ is the map
\begin{eqnarray*}
\mathrm{hol}([\partial u]):\bigotimes_{i=1}^k \hom(\cE_{i-1|p_i},\cE_{i|p'_i})&\to& \hom(\cE_{0|p_0},\cE_{k|p'_0})\\[-6pt]
(\rho_1,\dots,\rho_k)&\mapsto& \gamma_k\cdot \rho_k\cdot \dots \cdot \gamma_1\cdot \rho_1\cdot \gamma_0,
\end{eqnarray*}
where for $i=0,\dots,k$ we denote by $\gamma_i\in \hom(\cE_{i|p'_i},\cE_{i|p_{i+1}})$ the isomorphism
defined by parallel transport in the fibers of $\cE_i$ along the
portion of $u(\partial S)$ that lies on $L_i$.
\end{defi}

For simplicity, and since our main focus is not on the wrapped setting, our weights are defined in terms of symplectic area, rather than the more commonly used
{\em topological energy}
$$E([u])=\int_S u^*\omega-d((u^*h)\,\beta).$$
The two notions are equivalent up to rescaling each generator $p$ by $T^{\varepsilon h(p)}$, or
by simply taking the limit $\varepsilon\to 0$ in our choices of Hamiltonian perturbations, except for generators in
wrapped noncompact ends of $M.$ In this latter case, it is more advantageous 
to use action rescaling to eliminate the area contributions of wrapped components
of propagating discs (involving only v.b.-type objects) altogether.

The final ingredient for the definition of $\mu^k$ is the orientation of
the zero-dimensional moduli spaces of rigid propagating discs; this works
just as in ordinary Floer theory on Riemann surfaces, following a recipe
due to Seidel \cite[\S 13]{SeBook}. First we fix orientations for our objects
in a manner consistent with the choices made above for grading, namely
objects of point type loop clockwise around the origin
in each component of $M$, and v.b.-type objects to run from 0 to $\infty$ in each component of $M$.
Given a propagating disc $u:S\to M$ with inputs $(p_i,p'_i)\in \cX(L_{i-1},L_i)$
and output $(p_0,p'_0)\in \cX(L_0,L_k)$, for each $i=0,\dots,k$, if $\deg p_i$ is even
then we set $(-1)^{\sigma_i}=+1$, whereas if $\deg p_i$ is odd
we assign $(-1)^{\sigma_i}=+1$ if the orientation of $L_i$ ($L_k$ in the case of $i=0$)
at $p'_i$ agrees with that of the oriented curve $u(\partial S)$, and $-1$ otherwise.
The overall sign is then $(-1)^{\sigma(u)}=\prod_{i=0}^k (-1)^{\sigma_i}$.  Finally:

\begin{defi} Given $(p_i,p'_i)\in \cX(L_{i-1},L_i)$ and
$\rho_i\in \cE_{i-1|p_i}^*\otimes \cE_{i|p'_i}$ for $1\leq i\leq k$, we set
$$\mu^k(\rho_k,\dots,\rho_1)=\sum_{\substack{
(p_0,p'_0)\in \cX(L_0,L_k)\\[2pt] [u]\ \text{rigid}\\[1pt]
u\in \cM(p_0,\dots,p_k,[u])}} 
(-1)^{\sigma(u)}\,\Pi C([u])\,W([u])\,\mathrm{hol}([\partial u])(\rho_1,\dots,\rho_k).$$
\end{defi}

\subsection{The $A_\infty$-relations} We now state and prove
\begin{theo}\label{thm:ainfty}
The operations $\mu^k$ defined above satisfy the 
$A_\infty$-relations 
\begin{equation}\label{eq:ainfty}
\sum_{\ell=1}^{k}\sum_{j=0}^{k-\ell} (-1)^* \mu^{k+1-\ell}(\rho_k,\dots,
\rho_{j+\ell+1},\mu^\ell(\rho_{j+\ell},\dots,\rho_{j+1}),\rho_j,\dots,\rho_1)=0
\end{equation}
where $*=j+\deg(p_1)+\dots+\deg(p_j)$.
\end{theo}
The proof relies on the same geometric idea as
in the usual case, namely showing that 1-dimensional moduli spaces define
cobordisms between the pairs of rigid configurations which appear in the 
left-hand side of \eqref{eq:ainfty}, but the argument requires substantial modifications
to account for propagation through the nodes of $M$.
We start with the following lemma, which is a direct consequence of the analogous statement for ordinary discs:

\begin{lemma}\label{l:localslit}
The one-dimensional strata of moduli spaces of propagating holomorphic discs 
correspond to configurations where one of the disc components has a single
boundary branch point (and is otherwise immersed), or is a constant map carrying four marked points, one
of which is the output, while all the other components are rigid
(i.e., immersed polygons with locally convex corners).

Near a boundary branch point, the boundary of the disc
doubles back onto itself along a ``slit'',
and the deformation proceeds by moving the
branch point along the boundary, either extending the
slit further into the disc or shrinking it. The ends of such a one-dimensional
stratum arise when either:

$(i)$ the slit shrinks all the way into a marked point, or 

$(ii)$ the slit extends all the way across the disc to break it into a pair of
discs attached to each other at a Floer generator or at a node.

\end{lemma}

\begin{proof}
The statement follows from automatic regularity for 
holomorphic discs on Riemann surfaces (see e.g.\ \cite[\S 13]{SeBook})
and from the fact that the dimension of each moduli space of (non-constant)
discs is equal to the number of boundary branch points plus twice the number 
of interior branch points, with local coordinates given by the positions of the
images of the branch points (see e.g.\ \cite[Proposition 7.8]{ENS}).
\end{proof}

We will not discuss constant components with four marked points any further, 
as the only instance where one explicitly needs to deal
with them is the case of a pair of constant discs at a node, which does not involve any propagation.


Looking at the entirety of a propagating holomorphic disc with a boundary
branch point, the portion of the boundary that backtracks onto itself (the
slit) may extend beyond the component that carries the branch point. Namely,
assume a propagating holomorphic disc $u:S\to M$ has a boundary branch point along an
edge mapping to $L_i$. Then the portion of $\partial S$ that runs from
$z_{i-1}$ to $z_i$ maps to an arc that
travels from $p_{i-1}$ to $p_i$ along $L_i$, possibly via the nodes of $M,$
doubling back on itself exactly once at
the branch point. The portion of the boundary that gets traced twice (the {\em slit}) can run through nodes of $M,$ in which case we say that the
slit {\em propagates} through these nodes.
The slit ends either at one of
$p_{i-1}$ or $p_i$, in which case the image of the map has a {\em concave
corner} (unless $p_{i-1}=p_i$, a {\em repeated corner}), or at a node of $M$ where the arcs running from the branch point to
$p_{i-1}$ and to $p_i$ bifurcate into different components. In the latter case,
we say that the image of the map has a {\em bifurcated node}.
These possibilities are mutually exclusive.

To prove Theorem \ref{thm:ainfty}, we consider the union of the various
1-dimensional moduli spaces of propagating holomorphic discs with boundary on given Lagrangians
$L_0,\dots,L_k$ and corners at fixed input and output generators
$p_1,\dots,p_k,p_0$, whose image $u(S)$ covers a fixed domain in $M$ with fixed
multiplicities; as noted above, this domain must have either one concave (or
repeated) corner or one bifurcated node.
The different 1-dimensional strata correspond to the different ways
in which the slit can extend from the concave (or repeated) corner or from
the bifurcated node, propagate through
nodes along the way, and reach a branch point inside one of the components.

We glue these moduli spaces to each other at their end points whenever those
correspond to propagating discs whose slits trace exactly the same path (ending 
at a node of $M$), to form a one-dimensional cell complex whose ends
correspond to broken configurations; 
see e.g.\ Figure \ref{fig:ainfty}. The difficulty is that in general this cell complex 
has the structure of a {\em tree}, with edges branching off in multiple
directions, rather than a 1-dimensional manifold.
The root of the tree corresponds to the case where the slit has shrunk to a
point; the different paths away from the root correspond to the different
paths along which the slit can grow.


%


Thus, the key step in the proof of Theorem \ref{thm:ainfty} is a property which we call {\em invariance
of the total propagation multiplicity} under extending the slit past a node:
\begin{lemma}\label{l:tot_mult_invariance}
The propagation
multiplicity of a configuration whose slit stops just before reaching a node of $M$
is equal to the sum of the propagation multiplicities of the various
configurations that can be obtained by extending the slit slightly past the
node, 
plus those of
any broken configurations exhibiting a slit that ends
exactly at the node.
\end{lemma}

\noindent
We will show in \S \ref{sss:concave} that this lemma implies the desired cancellation between ends of the union of
moduli spaces for configurations whose image has a concave corner or 
a repeated corner. We then prove the lemma in \S \ref{sss:proof_lemma_tot}.
The case of bifurcated nodes requires an additional ingredient, which we
discuss in \S \ref{sss:bifurcated}.

\subsubsection{Propagating discs with one concave corner}\label{sss:concave}

\begin{figure}[t]
\begin{tikzpicture}[scale = 0.5]
\begin{scope}[shift={(-9,2.5)}] 
\fill[orange!25!white](-3.5,0) arc (-135:-45:1.414 and 2) arc (-135:-45:1.06 and 2) arc (-135:-45:1.06 and 2) arc (-135:-45:1.414 and 2) arc (45:135:1.414 and 2) arc (45:135:2.12 and 3) arc (45:135:1.414 and 2);
\fill[gray!15!white](-3.5,0) arc (-135:-45:1.414 and 2) arc (-135:-45:1.06 and 2) arc (45:135:1.06 and 1.15) arc (-90:-132.5:2) arc (105:135:1.414 and 2);
\draw[green!50!gray,thick](-3.5,0) arc (135:45:1.414 and 2) arc (135:45:2.12 and 3) arc (135:45:1.414 and 2);
\draw[red,thick](-3.5,0) arc (-135:-45:1.414 and 2) arc (-135:-45:1.06 and 2);
\draw[blue,thick,double](-1.5,0) arc (-90:-132.5:2);
\draw[blue,thick,double](-1.5,0) arc (135:45:1.06 and 1.15); 
\draw[blue,thick](0,0) arc (-135:-45:1.06 and 2) arc (-135:-45:1.414 and 2);
\draw[green!50!gray](-1.5,0.7)node {\SMALL $L_0$};
\draw[red](-1.5,-0.7)node {\SMALL $L_1$};
\draw[blue](1.5,-0.7)node {\SMALL $L_2$};
\fill (-3.5,0) circle (0.08) (-3.6,-0.4) node {\SMALL $p_1$};
\fill (3.5,0) circle (0.08) (3.6,-0.5) node {\SMALL $p_0$};
\fill (-1.5,0) circle (0.08) (1.5,0) circle (0.08);
\fill (0,0) circle (0.08) (0,-0.7) node {\SMALL $p_2$};
\fill (-2.82,0.53) circle (0.08);
\end{scope}
\begin{scope}[shift={(-11,-2.5)}] 
\fill[orange!25!white](-3.5,0) arc (-135:-45:1.414 and 2) arc (-135:-45:1.06 and 2) arc (-135:-45:1.06 and 2) arc (-135:-45:1.414 and 2) arc (45:135:1.414 and 2) arc (45:135:2.12 and 3) arc (45:135:1.414 and 2);
\fill[gray!15!white](-1.5,0) arc (-135:-45:1.06 and 2) arc (45:135:1.06 and 1.15) arc (90:132.5:2) arc (-105:-45:1.414 and 2);
\draw[green!50!gray,thick](-3.5,0) arc (135:45:1.414 and 2) arc (135:45:2.12 and 3) arc (135:45:1.414 and 2);
\draw[red,thick](-3.5,0) arc (-135:-45:1.414 and 2) arc (-135:-45:1.06 and 2);
\draw[blue,thick,double](-1.5,0) arc (90:132.5:2);
\draw[blue,thick,double](-1.5,0) arc (135:45:1.06 and 1.15); 
\draw[blue,thick](0,0) arc (-135:-45:1.06 and 2) arc (-135:-45:1.414 and 2);
\draw[green!50!gray](-1.5,0.7)node {\SMALL $L_0$};
\draw[red](-1.5,-0.7)node {\SMALL $L_1$};
\draw[blue](1.5,-0.7)node {\SMALL $L_2$};
\fill (-3.5,0) circle (0.08) (-3.6,-0.4) node {\SMALL $p_1$};
\fill (3.5,0) circle (0.08) (3.6,-0.5) node {\SMALL $p_0$};
\fill (-1.5,0) circle (0.08) (1.5,0) circle (0.08);
\fill (0,0) circle (0.08) (0,-0.7) node {\SMALL $p_2$};
\fill (-2.82,-0.53) circle (0.08);
\end{scope}
\begin{scope}[shift={(-2.5,-2.5)}] 
\fill[orange!25!white](-3.5,0) arc (-135:-45:1.414 and 2) arc (-135:-45:1.06 and 2) arc (-135:-45:1.06 and 2) arc (-135:-45:1.414 and 2) arc (45:135:1.414 and 2) arc (45:135:2.12 and 3) arc (45:135:1.414 and 2);
\fill[gray!15!white](-1.5,0) arc (-135:-45:1.06 and 2) arc (45:135:1.06 and 1.15) arc (45:135:1.414 and 2) arc (-135:-45:1.414 and 2);
\draw[green!50!gray,thick](-3.5,0) arc (135:45:1.414 and 2) arc (135:45:2.12 and 3) arc (135:45:1.414 and 2);
\draw[red,thick](-3.5,0) arc (-135:-45:1.414 and 2) arc (-135:-45:1.06 and 2);
\draw[blue,thick,double](-1.5,0) arc (135:45:1.06 and 1.15); 
\draw[blue,thick](0,0) arc (-135:-45:1.06 and 2) arc (-135:-45:1.414 and 2);
\draw[green!50!gray](-1.5,0.7)node {\SMALL $L_0$};
\draw[red](-1.5,-0.7)node {\SMALL $L_1$};
\draw[blue](1.5,-0.7)node {\SMALL $L_2$};
\fill (-3.5,0) circle (0.08) (-3.6,-0.4) node {\SMALL $p_1$};
\fill (3.5,0) circle (0.08) (3.6,-0.5) node {\SMALL $p_0$};
\fill (-1.5,0) circle (0.08) (1.5,0) circle (0.08);
\fill (0,0) circle (0.08) (0,-0.7) node {\SMALL $p_2$};
\end{scope}
\begin{scope}[shift={(0,2.5)}] 
\fill[orange!25!white](-3.5,0) arc (-135:-45:1.414 and 2) arc (-135:-45:1.06 and 2) arc (-135:-45:1.06 and 2) arc (-135:-45:1.414 and 2) arc (45:135:1.414 and 2) arc (45:135:2.12 and 3) arc (45:135:1.414 and 2);
\draw[green!50!gray,thick](-3.5,0) arc (135:45:1.414 and 2) arc (135:45:2.12 and 3) arc (135:45:1.414 and 2);
\draw[red,thick](-3.5,0) arc (-135:-45:1.414 and 2) arc (-135:-45:1.06 and 2);
\draw[blue,thick] (-0.53,0.32) circle  (0.05);
\draw[blue,thick,double](0,0) arc (45:80:1.06 and 1.15);
\draw[blue,thick](0,0) arc (-135:-45:1.06 and 2) arc (-135:-45:1.414 and 2);
\draw[green!50!gray](-1.5,0.7)node {\SMALL $L_0$};
\draw[red](-1.5,-0.7)node {\SMALL $L_1$};
\draw[blue](1.5,-0.7)node {\SMALL $L_2$};
\fill (-3.5,0) circle (0.08) (-3.6,-0.4) node {\SMALL $p_1$};
\fill (3.5,0) circle (0.08) (3.6,-0.5) node {\SMALL $p_0$};
\fill (-1.5,0) circle (0.08) (1.5,0) circle (0.08);
\fill (0,0) circle (0.08) (0,-0.7) node {\SMALL $p_2$};
\end{scope}
\begin{scope}[shift={(6,-2.5)}] 
\fill[orange!25!white](-3.5,0) arc (-135:-45:1.414 and 2) arc (-135:-45:1.06 and 2) arc (-135:-45:1.06 and 2) arc (-135:-45:1.414 and 2) arc (45:135:1.414 and 2) arc (45:135:2.12 and 3) arc (45:135:1.414 and 2);
\draw[green!50!gray,thick](-3.5,0) arc (135:45:1.414 and 2) arc (135:45:2.12 and 3) arc (135:45:1.414 and 2);
\draw[red,thick](-3.5,0) arc (-135:-45:1.414 and 2) arc (-135:-45:1.06 and 2);
\draw[red,thick] (0.53,0.32) circle  (0.05);
\draw[red,thick,double](0,0) arc (135:100:1.06 and 1.15);
\draw[blue,thick](0,0) arc (-135:-45:1.06 and 2) arc (-135:-45:1.414 and 2);
\draw[green!50!gray](-1.5,0.7)node {\SMALL $L_0$};
\draw[red](-1.5,-0.7)node {\SMALL $L_1$};
\draw[blue](1.5,-0.7)node {\SMALL $L_2$};
\fill (-3.5,0) circle (0.08) (-3.6,-0.4) node {\SMALL $p_1$};
\fill (3.5,0) circle (0.08) (3.6,-0.5) node {\SMALL $p_0$};
\fill (-1.5,0) circle (0.08) (1.5,0) circle (0.08);
\fill (0,0) circle (0.08) (0,-0.7) node {\SMALL $p_2$};
\end{scope}
\begin{scope}[shift={(9,2.5)}] 
\fill[gray!15!white](-3.5,0) arc (-135:-45:1.414 and 2) arc (-135:-45:1.06 and 2) arc (-135:-45:1.06 and 2) arc (-135:-45:1.414 and 2) arc (45:135:1.414 and 2) arc (45:135:2.12 and 3) arc (45:135:1.414 and 2);
\fill[orange!25!white](0,0) arc (-135:-45:1.06 and 2) arc (-135:-45:1.414 and 2) arc (45:75:1.414 and 2) arc (-47.5:-90:2) arc (45:135:1.06 and 1.15);
\draw[green!50!gray,thick](-3.5,0) arc (135:45:1.414 and 2) arc (135:45:2.12 and 3) arc (135:45:1.414 and 2);
\draw[red,thick](-3.5,0) arc (-135:-45:1.414 and 2) arc (-135:-45:1.06 and 2);
\draw[red,thick,double](0,0) arc (135:45:1.06 and 1.15) arc (-90:-48:2);
\draw[blue,thick](0,0) arc (-135:-45:1.06 and 2) arc (-135:-45:1.414 and 2);
\draw[green!50!gray](-1.5,0.7)node {\SMALL $L_0$};
\draw[red](-1.5,-0.7)node {\SMALL $L_1$};
\draw[blue](1.5,-0.7)node {\SMALL $L_2$};
\fill (-3.5,0) circle (0.08) (-3.6,-0.4) node {\SMALL $p_1$};
\fill (3.5,0) circle (0.08) (3.6,-0.5) node {\SMALL $p_0$};
\fill (-1.5,0) circle (0.08) (1.5,0) circle (0.08);
\fill (0,0) circle (0.08) (0,-0.7) node {\SMALL $p_2$};
\fill (2.82,0.53) circle (0.08);
\end{scope}
\draw[thick](-9,0.5)--(-6,0.5) arc (90:60:4);
\draw[thick](-9,-0.5)--(-6,-0.5) arc (-90:-60:4);
\draw[thick](-4,0)--(9,0);
\draw[thick](-4,0) arc (0:-60:0.8);
\fill (-4.4,-0.7) circle (0.08);
\fill (-4,0) circle (0.08);
\fill (4,0) circle (0.08);
\fill (-9,-0.5) circle (0.08);
\fill (-9,0.5) circle (0.08);
\fill (9,0) circle (0.08);
\draw[->](-9,1.35)--(-9,0.75);
\draw[->](-10,-1.5)--(-9.1,-0.75);
\draw[->](-3.5,-1.45)--(-4.2,-0.9);
\draw[->](-0.5,1.35)--(-0.5,0.25);
\draw[->](3.5,-1.35)--(2,-0.25);
\draw[->](9,1.35)--(9,0.25);
\end{tikzpicture}
\caption{A one-dimensional family of propagating discs with a concave
corner}\label{fig:ainfty}
\end{figure}
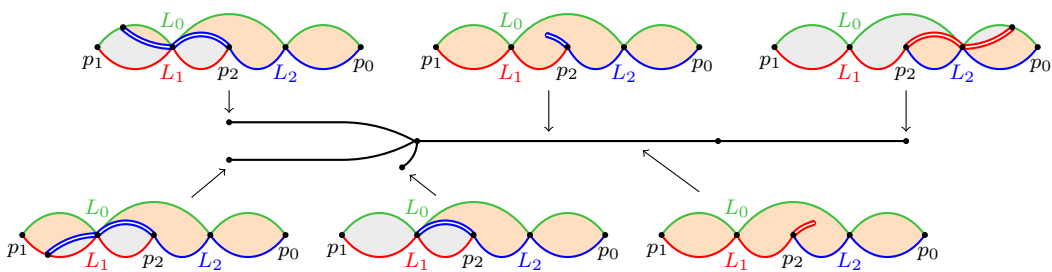

We first consider the case of a one-parameter family of propagating perturbed holomorphic discs with
$k+1$ marked points mapping to Floer generators $p_0,p_1,\dots,p_k$, and whose
image has a concave corner at one of the marked points, say 
$p_i\in \cX(L_{i-1},L_i)$. The simplest stratum in this family consists of
configurations where the slit does not propagate; in this case,
one of the components of $u(S)$ maps to an immersed polygon with 
a concave corner at $p_i$, and all the other components are rigid.
Such configurations deform by moving the boundary branch point along either
$L_{i-1}$ or $L_i$ to create a slit in the polygon, which extends from the concave 
corner along either Lagrangian as depicted in the central part of Figure \ref{fig:ainfty}.
In usual Lagrangian Floer theory on Riemann surfaces, each component of a
1-dimensional moduli space is an interval, whose ends are reached when
the slit extends all the way across and eventually hits the
boundary of the concave polygon, breaking it into a pair of smaller convex polygons.
These broken configurations contribute to the coefficient of $p_0$
in the left-hand side of \eqref{eq:ainfty}, and the $A_\infty$-relation
expresses the fact that they arise in cancelling pairs. (The
area and holonomy weights of broken configurations match those of the
unbroken configuration of which they are extremal deformations, hence they
are equal at both ends of the interval).

In our setting, as the slit extends across the concave polygon, 
it may hit a node through which the disc $u:S\to M$ propagates,
rather than the boundary of the disc.  When this happens, the
deformation naturally extends further (into different moduli spaces), as one can allow the slit to
propagate
into the next component of $u(S)$, and so on until it eventually hits the
boundary of the propagating disc. However, if the map $u$ locally covers more than once
the component of $M$ into which the slit is being extended,
there may be more than one way in which it can be slit along the
appropriate Lagrangian. This is illustrated on the left side of
Figure \ref{fig:ainfty}, where the left-most component of $S$ (a strip with
boundary on $L_0$ and $L_1$) is assumed to enter the left-side node with input degree $k_{in}=2$, 
so that there are two different ways in which this holomorphic disc can be
slit along $L_2$. 
Thus, the union of moduli spaces we consider is not an interval,
but rather a {\em tree} which may fork into several branches each
time the slit travels through a node.
An additional contribution to the boundary of the union of moduli spaces
can also arise when, rather than continuing through the node, the
slit stops at the node and breaking occurs through a constant component at
the node; see the bottom center diagram in Figure \ref{fig:ainfty}. (For
clarity this configuration is depicted in the figure as lying at the end of 
a short edge of the tree, but in fact there is no 1-dimensional stratum 
associated to this configuration.)

We claim that there is still a cancellation between the ends of the moduli space 
where $u$ breaks into a pair of propagating discs by extending a slit along either $L_{i-1}$ or
$L_i$. Area and holonomy weights behave just as in the usual case, so the key
ingredient is Lemma \ref{l:tot_mult_invariance}, which asserts an equality
between the total propagation multiplicities of the various configurations
that can arise before and after the slit extends through a node (up to any
broken configurations).  This in turn implies that
the sum of the propagation multiplicities of all the broken configurations
which arise as the slit extends towards one direction (for instance the three ends at the left of
Figure \ref{fig:ainfty}) is equal to the propagation multiplicity of the
initial disc $u$ -- and hence, arguing similarly when
the slit extends in the other direction, to the sum of the propagation multiplicities
of all the configurations at the other end of the moduli space (for instance
the single end at the right of Figure \ref{fig:ainfty}).

The argument is the same in the case of a repeated corner. The main
difference is that, if $p_{i}$ and $p_{i+1}$ coincide, then a slit can only
develop along $L_i$; when this slit shrinks all the
way to a point, a
constant disc contributing to $\mu^2(p_{i+1},p_i)$ breaks off from the rest
of the configuration. Lemma \ref{l:tot_mult_invariance} now implies that the
propagation multiplicity of this broken configuration with no slit is equal to the sum of
the propagation multiplicities of all the broken configurations where the
slit extends all the way across the propagating disc.

Thus, the argument for propagating discs whose image has a concave
or repeated corner has been reduced to Lemma \ref{l:tot_mult_invariance}, which we now
prove.

\subsubsection{Proof of Lemma \ref{l:tot_mult_invariance}}\label{sss:proof_lemma_tot}
The proof of Lemma \ref{l:tot_mult_invariance} is essentially combinatorial
in nature, as we need to compare the propagation multiplicities of configurations where
the slit stops just before a node vs.\ those where the slit has extended across
the node. There are three different cases, which need to be handled
separately: (1) the slit propagates ``backwards''
through a node, i.e.\ we extend it towards a component of $S$ that lies further away from the output
marked point $z_0$,
as in Figure \ref{fig:ainfty} left; (2) it propagates
``forward'' through a node, i.e.\ we extend it towards a component that lies
closer to the output marked point, as in Figure \ref{fig:ainfty}
right; and (3) we extend the slit through a node that
connects to a constant disc component at the output. In each case, the
desired equality of propagation multiplicities reduces to a purely combinatorial identity (Lemmas
\ref{l:ainfty_combi_backwards}--\ref{l:ainfty_combi_outputbounce}).

\subsubsection*{Case 1: the slit propagates backwards through a node}
Consider a node of $S$, mapping to a node $p_v\in M$, where the output 
of a component $D_{in}$ of $S$ mapping to $M_{e_{in}}$, with
boundaries on $L_{j_1}$ and $L_{j_2}$ ($j_1<j_2$), is attached to an
input end of a component $D_{out}$ mapping to $M_{e_{out}}$. Denote by $k_{in}\geq 1$ and
$k_{out}\geq 0$ the degrees of $u$ in the two strip-like ends near the node.
Assume that a slit is being extended along $L_i$ from the component $D_{out}$
backwards through the node and into $D_{in}$. Since the slit comes in from
$D_{out}$, necessarily either $i<j_1$ or $i>j_2$. When $i>j_2$ as pictured
on Figure \ref{fig:ainfty} left (resp.\ $i<j_1$), once extended into $D_{in}$ the slit breaks the local picture 
near $p_v$ into two propagating discs:
\begin{itemize}
\item one with boundary on $L_{j_2}$
and $L_i$ (resp.\ $L_i$ and $L_{j_1}$), propagating {\em backwards} through the node from $M_{e_{out}}$
to $M_{e_{in}}$, with input degree $1\leq a\leq k_{out}$ in $M_{e_{out}}$
determined by the position of the incoming slit within $D_{out}$, and arbitrary output degree $0\leq
b\leq k_{in}-1$ in $M_{e_{in}}$ (there are $k_{in}-1$ choices for how to
extend the slit into $D_{in}$);
\item the other with boundary on $L_{j_1}$ and $L_i$ (resp.\ $L_i$ and $L_{j_2}$), propagating forward
from $M_{e_{in}}$ to $M_{e_{out}}$, with input degree $k_{in}-b$ in
$M_{e_{in}}$ and output degree $k_{out}-a$ in $M_{e_{out}}$.
\end{itemize}
When $a=k_{out}$, another possibility (corresponding to the bottom center diagram of
Figure~\ref{fig:ainfty}) is that the slit ends at
the node $p_v$ and breaks the configuration into:
\begin{itemize}
\item an incoming propagating disc
(involving all the components of $u$ that lie on the $D_{in}$ side of the
node, plus one of the two pieces delimited by the slit on the $D_{out}$
side; in gray on Figure \ref{fig:ainfty} bottom center) 
that comes into the node from both $D_{in}$ and $D_{out}$ and ends with a 
constant component at $p_v$, and 
\item an outgoing propagating disc (the remaining portions of the curve on
the $D_{out}$ side) which
has an input at $p_v$ with boundary on $L_{j_1}$ and $L_i$ (resp.\ $L_i$ and
$L_{j_2}$), with local
degree $k_{out}-a=0$ as required by Definition \ref{def:proptree} for inputs at nodes.
\end{itemize}

Recall that the propagation coefficient at the node $v$ for the initial configuration (with
local degrees $k_{in}$ and $k_{out}$),
$C^{v;e_{in},e_{out}}_{k_{in},k_{out}}$, is defined to be the coefficient of
$t_{e_{out}/v}^{k_{out}}$ in the expansion of $t_{e_{in}/v}^{-k_{in}}$ as a
power series in $t_{e_{out}/v}$; whereas the product of the propagation
coefficients for the two nodes after inserting the slit as described above
is $C^{v;e_{out},e_{in}}_{a,b} C^{v;e_{in},e_{out}}_{k_{in}-b,k_{out}-a}$.
Meanwhile, in the case of a broken configuration involving a constant component at $p_v$ 
(for $a=k_{out}$), the contribution of the nodes of the constant component
to the propagation multiplicity is $K^{v;e_{in},e_{out}}_{k_{in},k_{out}}$,
the coefficient of the constant term in the expression of 
$t_{e_{in}/v}^{-k_{in}}t_{e_{out}/v}^{-k_{out}}$ as a linear combination of
negative powers of $t_{e_{in}/v}$ and $t_{e_{out}/v}$. Thus, 
the invariance of the total propagation multiplicities under
passing the slit through the node follows from:

\begin{lemma} \label{l:ainfty_combi_backwards}
Given $v,e_1,e_2$, and integers $k_{1}\geq 1$ and $k_{2}\geq a\geq 1$,
\begin{equation} \label{eq:ainfty_combi_backwards}
\sum_{b=0}^{k_{1}-1} C^{v;e_{2},e_{1}}_{a,b} C^{v;e_{1},e_{2}}_{k_{1}-b,k_{2}-a}
+\delta_{a,k_{2}} K^{v;e_{1},e_{2}}_{k_{1},k_{2}} =
C^{v;e_{1},e_{2}}_{k_{1},k_{2}}.
\end{equation}
\end{lemma}

\proof
Denote $t_1=t_{e_1/v}$ and $t_2=t_{e_2/v}$. 
The rational function $t_1^{-k_1}t_2^{-a}$ has a partial fraction
decomposition into a finite linear combination of
$1,t_1^{-1},\dots,t_1^{-k_1},t_2^{-1},\dots,t_2^{-a}$, so that we can write
$t_1^{-k_1}t_2^{-a}=K^{v;e_1,e_2}_{k_1,a}+P_1(t_1^{-1})+P_2(t_2^{-1})$, where
$P_1(t_1^{-1})$ is a polynomial 
in $t_1^{-1}$ without constant term (the polar part at $x_1$), and $P_2(t_2^{-1})$ is a
polynomial in $t_2^{-1}$ without constant term (the polar part at $x_2$).
On the other hand, near $x_1$ we have the power series expansion
$t_2^{-a}=\sum\limits_{b\geq 0} C^{v;e_2,e_1}_{a,b} t_1^b$, which yields the
Laurent series expression
$$t_1^{-k_1} t_2^{-a}=\sum_{b\geq 0} C^{v;e_2,e_1}_{a,b} t_1^{b-k_1}.$$
Comparing the polar parts at $x_1$ (i.e., using the fact that $P_2$ expands
near $x_1$ as a power series in $t_1$ without negative powers), we conclude that
$$P_1(t_1^{-1})=\sum_{b=0}^{k_1-1} C^{v;e_2,e_1}_{a,b} t_1^{b-k_1}.$$
This, in turn, yields a Laurent series expression for $t_1^{-k_1}t_2^{-a}$
near $x_2$, using the fact that each monomial in $P_1$ has a power series
expansion $t_1^{b-k_1}=\sum\limits_{c\geq 0}C^{v;e_1,e_2}_{k_1-b,c}\,t_2^c$:
$$t_1^{-k_1}t_2^{-a}=P_2(t_2^{-1})+K^{v;e_1,e_2}_{k_1,a}+
\sum_{b=0}^{k_1-1}\sum_{c\geq 0}
C^{v;e_2,e_1}_{a,b} C^{v;e_1,e_2}_{k_1-b,c}\, t_2^c.$$
On the other hand, starting from $t_1^{-k_1}=\sum\limits_{d\geq 0}
C_{k_1,d}^{v;e_1,e_2}\,t_2^d$ we also have
$$t_1^{-k_1} t_2^{-a}=\sum_{d\geq 0} C^{v;e_1,e_2}_{k_1,d} t_2^{d-a}.$$
Comparing the coefficients of $t_2^{k_2-a}$ in these two
expressions immediately gives \eqref{eq:ainfty_combi_backwards}.
\endproof

\subsubsection*{Case 2: the slit propagates forward through a node}
Consider again a node of $S$, mapping to a node $p_v\in M$, where the output 
of a component $D_{in}$ of $S$ mapping to $M_{e_{in}}$, with
boundaries on $L_{j_1}$ and $L_{j_2}$ ($j_1<j_2$), is attached to an
input end of a component $D_{out}$ mapping to $M_{e_{out}}$. Assume furthermore
that the restriction of $u$ to $D_{out}$ is not a constant map, and denote by $k_{in}\geq 1$ and
$k_{out}\geq 0$ the degrees of $u$ in the two strip-like ends near the node.
Assume that a slit is being extended along $L_i$ from the component $D_{in}$
forward through the node and into $D_{out}$. Since the slit comes in from
$D_{in}$, necessarily $j_1<i<j_2$ (see Figure \ref{fig:ainfty} right), and
once extended into $D_{out}$ the slit breaks the local picture 
near $p_v$ into two propagating discs, both going forward through the node
from $M_{e_{in}}$ to $M_{e_{out}}$, one of them with input degree $1\leq
a\leq k_{in}-1$ (determined by the position of the slit in $D_{in}$) and
output degree $0\leq b\leq k_{out}$ (which can be chosen freely, there are
$k_{out}+1$ choices for how to extend the slit into $D_{out}$), while the
other has input degree $k_{in}-a$ and output degree $k_{out}-b$. The
invariance of total propagation multiplicities then reduces to:

\begin{lemma} Given $v,e_1,e_2$ and integers $1\leq a\leq k_1-1$ and
$k_2\geq 0$, \begin{equation}\label{eq:ainfty_combi_forward}
\sum_{b=0}^{k_2} C_{a,b}^{v;e_1,e_2} C_{k_1-a,k_2-b}^{v;e_1,e_2}=C_{k_1,k_2}^{v;e_1,e_2}.  
\end{equation}
\end{lemma}

\proof This is immediate from expressing $t_1^{-k_1}$ as the product of
$t_1^{-a}=\sum_{b\geq 0} C_{a,b}^{v;e_1,e_2}\, t_2^b$ and
$t_1^{a-k_1}=\sum_{d\geq 0} C_{k_1-a,d}^{v;e_1,e_2}\, t_2^d$, and taking the
coefficient of $t_2^{k_2}$ in the resulting power series.\endproof

\subsubsection*{Case 3: the slit propagates through a constant output component}

\begin{figure}[t]
\begin{tikzpicture}[scale = 0.5]
\begin{scope}[shift={(-4,0)}] 
\fill[orange!25!white](-2,-0.828) arc (-90:-45:2.828) arc (-135:-90:2.828)--(2,0.828) arc (90:135:2.828) arc (45:90:2.828);
\draw[green!50!gray,thick](-2,0.828) arc (90:45:2.828);
\draw[red,thick](-2,-0.828) arc (-90:-45:2.828) arc (-135:-90:2.828);
\draw[blue,thick](0,0) arc (135:90:2.828);
\fill (1.02,-0.32) circle (0.08);
\draw[thick,double](2,-0.414) arc (-90:-110:2.828 and 1.414);
\draw[green!50!gray](-1.5,1.2) node {\SMALL $L_0$};
\draw[red](-0.4,-1.1)node {\SMALL $L_j$};
\draw[blue](1.5,1.2)node {\SMALL $L_k$};
\draw(2.5,-0.4)node {\SMALL $L_i$};
\fill (0,0) circle (0.12) (0,0.45) node {\SMALL $p_0$};
\end{scope}
\draw[->,thick](-0.5,0)--(1.5,0);
\begin{scope}[shift={(4.5,0)}] 
\fill[orange!25!white](-2,-0.828) arc (-90:-45:2.828) arc (-135:-90:2.828)--(2,0.828) arc (90:135:2.828) arc (45:90:2.828);
\fill[gray!15!white](-2,-0.828) arc (-90:-45:2.828) arc (-135:-90:2.828)--(2,-0.414) arc (-90:-135:2.828 and 1.414) arc (-45:-90:2.828 and 1.414);
\draw[green!50!gray,thick](-2,0.828) arc (90:45:2.828);
\draw[red,thick](-2,-0.828) arc (-90:-45:2.828) arc (-135:-90:2.828);
\draw[blue,thick](0,0) arc (135:90:2.828);
\draw[thick,double](-2,-0.414) arc (-90:-45:2.828 and 1.414) arc (-135:-90:2.828 and 1.414);
\draw[green!50!gray](-1.5,1.2) node {\SMALL $L_0$};
\draw[red](-0.4,-1.1)node {\SMALL $L_j$};
\draw[blue](1.5,1.2)node {\SMALL $L_k$};
\draw(2.5,-0.4)node {\SMALL $L_i$};
\fill (0,0) circle (0.12) (0,0.45) node {\SMALL $p_0$};
\end{scope}
\end{tikzpicture}
\caption{Extending a slit through a constant output component}\label{fig:slitoutput}
\end{figure}
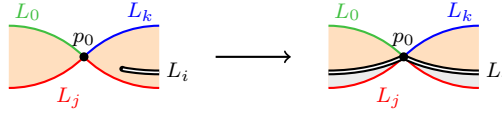

Next we consider the case where a slit is extended along $L_i$ 
into a constant output component $D_{out}$ (a constant disc with 
two inputs, carrying the output marked point $z_0\in S$ and mapping to a node $p_0=p_v\in M$). 
Denote by $D_1$ and $D_2$ the two components of $S$ adjacent to $D_{out}$,
$M_{e_1}$ and $M_{e_2}$ the components of $M$ into which they map,  and
$k_1,k_2\geq 1$ the degrees of the restrictions of $u$ to their output strip-like ends.
A slit which extends along $L_i$ within the component $D_2$ and reaches the constant
output component can be extended further back into $D_1$, as shown in
Figure \ref{fig:slitoutput}. This decomposes the local picture near $p_v$
into two propagating discs:
\begin{itemize}
\item one with boundary on $L_i$ and $L_j$, which propagates
through the node from $D_2$ towards $D_1$, with input degree $1\leq a\leq
k_2-1$ in $M_{e_2}$ as
determined by the position of the slit in $D_2$, and output degree $0\leq
b\leq k_1-1$ in $M_{e_1}$ (there are $k_1$ possible choices);
\item the other with incoming ends with boundaries on $L_0$ and $L_i$ on one
hand and $L_i$ and $L_k$ on the other hand, of degrees $k_1-b$ and $k_2-a$
in $M_{e_1}$ and $M_{e_2}$ respectively, ending at a constant component at
$p_0=p_v$.
\end{itemize}
The invariance of the sum of all propagation multiplicities now follows from:
\begin{lemma} \label{l:ainfty_combi_outputbounce} Given $v,e_1,e_2$ and integers $k_1\geq 1$ and $1\leq a\leq k_2-1$,
\begin{equation}\label{eq:ainfty_combi_outputbounce}
\sum_{b=0}^{k_1-1} C_{a,b}^{v;e_2,e_1} K_{k_1-b,k_2-a}^{v;e_1,e_2}=
K_{k_1,k_2}^{v;e_1,e_2}.
\end{equation}
\end{lemma}
\proof
As in the proof of Lemma \ref{l:ainfty_combi_backwards}, setting $t_1=t_{e_1/v}$ and 
$t_2=t_{e_2/v}$, we start with the partial fraction decomposition
$t_1^{-k_1}t_2^{-a}=K_{k_1,a}^{v;e_1,e_2}+P_1(t_1^{-1})+P_2(t_2^{-1})$, and
recall that $$P_1(t_1^{-1})=\sum_{b=0}^{k_1-1}
C_{a,b}^{v;e_2,e_1} t_1^{b-k_1}.$$
Multiplying by $t_2^{a-k_2}$, we obtain 
\begin{equation}\label{eq:ainfty_combi_lemma}
t_1^{-k_1}t_2^{-k_2}=\Bigl(K_{k_1,a}^{v;e_1,e_2}+P_2(t_2^{-1})\Bigr) t_2^{a-k_2}+
\sum_{b=0}^{k_1-1} C_{a,b}^{v;e_2,e_1} t_1^{b-k_1}\,t_2^{a-k_2}.
\end{equation}
This expression can in turn be decomposed into partial fractions and
expressed as a linear combination of
$1,t_1^{-1},\dots,t_1^{-k_1},t_2^{-1},\dots,t_2^{-k_2}$; we are interested
in the constant term of this decomposition. The first part of the right-hand side of
\eqref{eq:ainfty_combi_lemma} only involves negative powers of $t_2$, so it
does not contribute to the constant term. Meanwhile,
the constant term in the partial fraction decomposition of $t_1^{b-k_1}\,t_2^{a-k_2}$ is
$K_{k_1-b,k_2-a}^{v;e_1,e_2}$; therefore, the constant term in the
partial fraction decomposition of the right-hand side  of
\eqref{eq:ainfty_combi_lemma} is
$$\sum_{b=0}^{k_1-1} C_{a,b}^{v;e_2,e_1}\,K_{b-k_1,a-k_2}^{v;e_1,e_2},$$
which is exactly the left-hand side of \eqref{eq:ainfty_combi_outputbounce}.
On the other hand, the constant term in the partial fraction decomposition
of $t_1^{-k_1}t_2^{-k_2}$ (the left-hand side of
\eqref{eq:ainfty_combi_lemma}) is, by definition, equal to $K_{k_1,k_2}^{v;e_1,e_2}$.
The lemma then follows by comparing these two quantities.
\endproof

\subsubsection{Bifurcated propagating discs}\label{sss:bifurcated}
Besides propagating discs with a concave corner, there is another type of
configuration which gives rise to 1-dimensional moduli spaces of propagating
discs: ``bifurcated'' discs in which, near one of the nodes $p_v$ of $M$, 
$S$ has {\em two} incoming components $D_1,D_2$ and one outgoing component
$D_{out}$, each of which maps to a different component of
$M$ ($M_{e_1}$, $M_{e_2}$, $M_{e_{out}}$, where $e_1,e_2,e_{out}$ are the
three edges meeting at $v$). 

\begin{figure}[t]
\begin{tikzpicture}[scale = 0.65]
\begin{scope}[shift={(0,0)}] 
\draw(-2,2)node {\SMALL (c)};
\fill[orange!25!white](0,0) arc (120:90:4) -- (2,-0.536) arc (-90:-120:4);
\fill[gray!15!white](0,0) arc (-120:-150:4) -- (-0.536,2) arc (30:0:4);
\fill[gray!15!white](0,0) arc (120:150:4) -- (-0.536,-2) arc (-30:0:4);
\draw[green!50!gray,thick](-0.536,2) arc (30:0:4) arc (120:90:4);
\draw[red,thick](-1.464,1.464) arc (-150:-120:4) arc (120:150:4);
\draw[blue,thick](-0.536,-2) arc (-30:0:4) arc (-120:-90:4);
\draw[green!50!gray](0.4,0.65) node {\small $L_i$};
\draw[red](-0.7,0)node {\small $L_j$};
\draw[blue](0.4,-0.65)node {\small $L_k$};
\fill (0,0) circle (0.08);
\draw(-0.7,1.3)node {\SMALL $D_1$};
\draw(-0.7,-1.3)node {\SMALL $D_2$};
\draw(1.3,0)node {\SMALL $D_{out}$};
\end{scope}
\begin{scope}[shift={(-5,0)}] 
\draw(-2,2)node {\SMALL (b)};
\fill[orange!25!white](0,0) arc (120:90:4) -- (2,-0.536) arc (-90:-120:4);
\fill[gray!15!white](0,0) arc (-120:-150:4) -- (-1,1.732) -- (0,0);
\fill[orange!25!white] (0,0) -- (-1,1.732) -- (-0.536,2) arc (30:0:4);
\fill[gray!15!white](0,0) arc (120:150:4) -- (-0.536,-2) arc (-30:0:4);
\fill[blue](-0.95,1.645) circle (0.06);
\draw[blue,thick,double] (0,0) -- (-0.95,1.645);
\draw[green!50!gray,thick](-0.536,2) arc (30:0:4) arc (120:90:4);
\draw[red,thick](-1.464,1.464) arc (-150:-120:4) arc (120:150:4);
\draw[blue,thick](-0.536,-2) arc (-30:0:4) arc (-120:-90:4);
\draw[green!50!gray](0.4,0.65) node {\small $L_i$};
\draw[red](-0.7,0)node {\small $L_j$};
\draw[blue](0.4,-0.65)node {\small $L_k$};
\fill (0,0) circle (0.08);
\draw(-1,1.2)node {\tiny $D_1$};
\draw(-0.7,-1.3)node {\SMALL $D_2$};
\draw(1.3,0)node {\SMALL $D_{out}$};
\end{scope}
\begin{scope}[shift={(-10,0)}] 
\draw(-2,2)node {\SMALL (a)};
\fill[orange!25!white](0,0) arc (120:90:4) -- (2,-0.536) arc (-90:-120:4);
\fill[gray!15!white](0,0) arc (120:150:4) -- (-1,-1.732) -- (0,0);
\fill[orange!25!white] (0,0) -- (-1,-1.732) -- (-0.536,-2) arc (-30:0:4);
\fill[gray!15!white](0,0) arc (-120:-150:4) -- (-0.536,2) arc (30:0:4);
\fill[green!50!gray](-0.95,-1.645) circle (0.06);
\draw[green!50!gray,thick,double] (0,0) -- (-0.95,-1.645);
\draw[green!50!gray,thick](-0.536,2) arc (30:0:4) arc (120:90:4);
\draw[red,thick](-1.464,1.464) arc (-150:-120:4) arc (120:150:4);
\draw[blue,thick](-0.536,-2) arc (-30:0:4) arc (-120:-90:4);
\draw[green!50!gray](0.4,0.65) node {\small $L_i$};
\draw[red](-0.7,0)node {\small $L_j$};
\draw[blue](0.4,-0.65)node {\small $L_k$};
\fill (0,0) circle (0.08);
\draw(-0.55,-1.6)node {\tiny $D_{2}$};
\draw(-0.7,1.3)node {\SMALL $D_1$};
\draw(1.3,0)node {\SMALL $D_{out}$};
\end{scope}
\begin{scope}[shift={(5,0)}] 
\draw(-2,2)node {\SMALL (d)};
\fill[orange!25!white](0,0) arc (120:90:4) -- (2,0) -- (0,0);
\fill[gray!15!white](0,0) arc (-120:-90:4) -- (2,0) -- (0,0);
\fill[orange!25!white](0,0) arc (-120:-150:4) -- (-0.536,2) arc (30:0:4);
\fill[gray!15!white](0,0) arc (120:150:4) -- (-0.536,-2) arc (-30:0:4);
\fill[red](1.9,0) circle (0.06);
\draw[red,thick,double] (0,0) -- (1.9,0);
\draw[green!50!gray,thick](-0.536,2) arc (30:0:4) arc (120:90:4);
\draw[red,thick](-1.464,1.464) arc (-150:-120:4) arc (120:150:4);
\draw[blue,thick](-0.536,-2) arc (-30:0:4) arc (-120:-90:4);
\draw[green!50!gray](0.4,0.65) node {\small $L_i$};
\draw[red](-0.7,0)node {\small $L_j$};
\draw[blue](0.4,-0.65)node {\small $L_k$};
\fill (0,0) circle (0.08);
\draw(-0.7,1.3)node {\SMALL $D_1$};
\draw(-0.7,-1.3)node {\SMALL $D_2$};
\draw(1.6,-0.3)node {\tiny $D_{out}$};
\end{scope}
\begin{scope}[shift={(4,-4)}]
\draw[thick,dotted](-10,0.5)--(-9,0.5);
\draw[thick,dotted](-13,-0.5)--(-12,-0.5);
\draw[thick,dotted](1,0)--(2,0);
\draw[thick](-9,0.5)--(-6,0.5) arc (90:60:4);
\draw[thick](-12,-0.5)--(-6,-0.5) arc (-90:-60:4);
\draw[thick](-4,0)--(1,0);
\draw[thick](-4,0) arc (0:60:0.8);
\fill (-4.4,0.7) circle (0.08);
\fill (-4,0) circle (0.08);
\draw[->](-9,1.65)--(-9,0.75);
\draw[->](-13,2)--(-12,-0.25);
\draw[->](-4.4,1.65)--(-4.4,0.9);
\draw[->](0.5,1.65)--(0.5,0.25);
\end{scope}
\end{tikzpicture}
\caption{Decomposing a bifurcated disc into a pair of propagating discs}\label{fig:ainfty_bifurcated}
\end{figure}
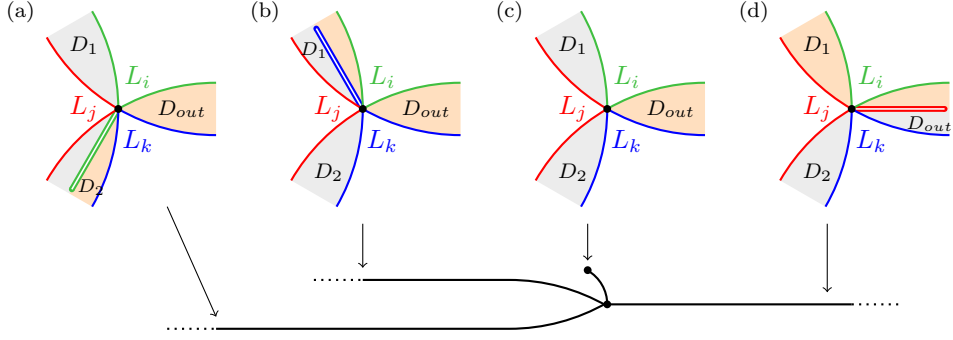

If the outgoing component near the bifurcated node does not surject locally
onto a neighborhood of the node in $M_{e_{out}}$ (i.e., the output degree is $k_{out}=0$), 
then such a bifurcated disc can be realized immediately as a
broken configuration of two rigid propagating discs, one including $D_1$ and
$D_2$ and ending at a constant component at $p_v$, and the other starting
with an input at $p_v$ and including $D_{out}$ (see Figure
\ref{fig:ainfty_bifurcated}\,(c)). In general (regardless of the value of $k_{out}$), 
this configuration can also deform by growing a slit into any one of the three
components $D_1$, $D_2$, or $D_{out}$, which has the effect of locally
breaking the bifurcated configuration into a pair of honest propagating
discs. Thus, the moduli space of propagating discs extends
into three types of directions, corresponding to the three
ways in which a slit can be created and extended into $S$; see Figure
\ref{fig:ainfty_bifurcated}\,(a)(b)(d). (For each of these
there may be multiple possibilities if the degree of that component of $u$
is greater than 1). As the slit expands into the appropriate component of
$S$, it will eventually either hit the boundary of the domain or pass through
other nodes and extend into other components, as in the case of discs with concave corners.
This part of the story works exactly as in the previous section and is
handled by Lemma \ref{l:tot_mult_invariance}; the new ingredient, rather, is 
the cancellation that needs to occur between the combinatorial propagation multiplicities 
associated to the various ways of creating a slit and locally decomposing a bifurcated node into a pair of rigid
propagating discs.

Denote by $k_1\geq 1$, $k_2\geq 1$ and $k_{out}\geq 0$ the degrees of $D_1$, $D_2$ and
$D_{out}$ near the bifurcated node. As noted above, if $k_{out}=0$ then
there is a broken configuration in which one of the two rigid propagating
discs contains $D_1$ and $D_2$ and ends with a constant component at $p_v$
(Figure \ref{fig:ainfty_bifurcated}\,(c)); the nodes adjacent to the constant component 
contribute $K^{v;e_1,e_2}_{k_1,k_2}$ to the propagation multiplicity of this broken configuration.
Meanwhile, for each $0\leq b\leq k_1-1$ there are
deformations in which a slit grows into the $D_1$ component, decomposing the
local picture into a propagating disc consisting of $D_2$ (incoming into
$p_v$) attached to part of $D_1$ (outgoing with degree $b$) on one hand, and
a propagating disc consisting of the remaining portion of $D_1$ (incoming
into $p_v$ with degree $k_1-b$) attached to $D_{out}$ (Figure \ref{fig:ainfty_bifurcated}\,(b)).
Similarly, there are configurations with a slit in the $D_2$ component,
where one propagating disc consists of $D_1$ (incoming into $p_v$) attached
to part of $D_2$ (outgoing with degree $0\leq a\leq k_2-1$), and the other
consists of the rest of $D_2$ (incoming into $p_v$ with degree $k_2-a$)
attached to $D_{out}$ (Figure \ref{fig:ainfty_bifurcated}\,(a)).
The last case is when the slit lies in $D_{out}$; one propagating disc
consists of $D_1$ attached to part of $D_{out}$ (outgoing with degree $0\leq
c\leq k_{out}$) and the other consists of $D_2$ attached to the rest of
$D_{out}$ (outgoing with degree $k_{out}-c$) (Figure
\ref{fig:ainfty_bifurcated}\,(d)). Comparing the sum of the propagation
multiplicities of the configurations with a slit in one of the input discs
to those with a slit in the output component $D_{out}$ then amounts to
checking the following identity:

\begin{lemma}\label{l:ainfty_tripod}
Given a vertex $v$ of $G$ with adjacent edges $e_1,e_2,e_3$, and
integers $k_1,k_2\geq 1$ and $k_3\geq 0$,
\begin{equation}\label{eq:ainfty_tripod}
\sum_{a=0}^{k_2-1} C_{k_1,a}^{v;e_1,e_2} C_{k_2-a,k_3}^{v;e_2,e_3}
+\sum_{b=0}^{k_1-1} C_{k_2,b}^{v;e_2,e_1} C_{k_1-b,k_3}^{v;e_1,e_3}
+\delta_{k_3,0}K_{k_1,k_2}^{v;e_1,e_2}=
\sum_{c=0}^{k_3} C_{k_1,c}^{v;e_1,e_3} C_{k_2,k_3-c}^{v;e_2,e_3}.
\end{equation}
\end{lemma}

\proof The equality follows from comparing two ways of expressing $t_1^{-k_1}t_2^{-k_2}$ as a
power series in $t_3$. On one hand, we can start from
$t_1^{-k_1}=\sum\limits_{c\geq 0} C_{k_1,c}^{v;e_1,e_3} t_3^c$ and
$t_2^{-k_2}=\sum\limits_{d\geq 0} C_{k_2,d}^{v;e_2,e_3} t_3^d$. Multiplying
these two expressions, we arrive at a power series in $t_3$ in which the
coefficient of $t_3^{k_3}$ is exactly the right-hand side of
\eqref{eq:ainfty_tripod}. On the other hand, we can proceed as 
in the proof of Lemma \ref{l:ainfty_combi_backwards} to obtain the partial
fraction decomposition
\begin{equation}\label{eq:partialproduct}
t_1^{-k_1}t_2^{-k_2}=K_{k_1,k_2}^{v;e_1,e_2}+P_1(t_1^{-1})+P_2(t_2^{-1})=
K_{k_1,k_2}^{v;e_1,e_2}+\sum_{b=0}^{k_1-1} C_{k_2,b}^{v;e_2,e_1}\,t_1^{b-k_1}
+\sum_{a=0}^{k_2-1} C_{k_1,a}^{v;e_1,e_2}\,t_2^{a-k_2}.
\end{equation}
Substituting $t_1^{b-k_1}=\sum\limits_{d\geq 0} C_{k_1-b,d}^{v;e_1,e_3}\, t_3^d$
and $t_2^{a-k_2}=\sum\limits_{d\geq 0} C_{k_2-a,d}^{v;e_2,e_3}\, t_3^d$, we
arrive at a power series in $t_3$ in which the coefficient of $t_3^{k_3}$ is
the left-hand side of \eqref{eq:ainfty_tripod}.
\endproof

This completes the case by case analysis and the proof of Theorem
\ref{thm:ainfty}.

\subsection{Infinite Hamiltonian perturbations}\label{ss:FMH}

We now describe a version of the Fukaya category of $M$ which can be
expressed in terms of local pieces. This construction involves large (in a certain sense, ``infinite'') 
Hamiltonian perturbations and is very similar to H. Lee's thesis
\cite{Lee}. While Lee decomposes a Riemann surface into pairs of pants, we
will instead use neighborhoods of the vertices
(i.e., {\em mirrors} of pairs of pants) as building blocks.

For each half-edge $e/v\in E,$ we choose an identification of $M_e$ with $[0,4]\times S^1/\sim,$ 
where $\{0\}\times S^1$ (resp.\ $\{4\}\times S^1$) is identified with $v$
(resp.\ $v'$), in such a way that the symplectic form is $\frac{A_e}{4}d\tau\wedge d\psi,$ 
where $\tau$ and $\psi$ are the coordinates on $[0,4]$ and $S^1=\R/\Z.$ 

We assume that near $\tau=1$ and near $\tau=3$ (and in fact, over the whole
support of the further perturbations we introduce below) the Hamiltonian $h$ used in the definition 
of the category $\cF(M)$ can be expressed as a function of the $\tau$ coordinate only.
Choose a sequence of $C^\infty$ functions $f_n:[0,4]\to\R$, constant away
from $\tau=1$ and $\tau=3$, and converging to a continuous
function $f:[0,4]\to\R$, such that:

\begin{itemize}
\item[(i)] $f=f_n=0$ near $\tau=0$ and $\tau=4$, and $f$ and $f_n$ are constant
near $\tau=2$; $f_n=f$ on $[0,1-\frac1n]\cup [3+\frac1n,4]$, and $f_n-f$ is
constant on $[1+\frac1n,3-\frac1n]$.

\item[(ii)] $f''_n\leq 0$ on $[0,1)\cup (3,4]$ and $f''_n\geq 0$ on $(1,3)$
(hence the same holds for $f''$);

\item[(iii)] $f'_n(1)=-n$, $f'_n(3)=n$, $\lim\limits_{\tau\to 1}f'(\tau)=-\infty,$ and $\lim\limits_{\tau\to 3}f'(\tau)=+\infty.$
\end{itemize}

\noindent 
We consider Hamiltonian perturbations $H_n=\varepsilon h+\frac{A_e}{4}f_n(\tau)$. The
assumption that $h$ only depends on $\tau$ over the support of $f'_n$ ensures
that the Hamiltonian flows generated by $h$ and $f_n$ commute, and that
the time 1 flow of $H_n$ differs from that of $\varepsilon h$ by a rotation
$\psi\mapsto \psi+f'_n(\tau)$.  The category $\cF(M;H)$ is defined using
$H_n$ instead of $\varepsilon h$ as Hamiltonian perturbation for Floer
complexes, and taking $n\to \infty$ in a manner we discuss below.

We impose some additional conditions on the objects of $\cF(M;H)$.
For v.b. type Lagrangians, we will assume that the coordinate $\tau$ is strictly 
monotonic on every component (so that the Lagrangian only passes once
through the ``necks'' at $\tau=1$ and $\tau=3$); we note that every v.b.
type object of $\cF(M)$ is isomorphic to an object which satisfies this
condition. We also assume that the generators of $CF^*(L,L';\varepsilon h)$ all lie away from $\tau=1$
and $\tau=3$. Point type Lagrangians aren't necessary for our argument, but can
be allowed as long as they are disjoint from the
circles at $\tau=1$ and $\tau=3$; this excludes objects which are supported at the
boundary of the pieces of our decomposition.

Given a pair of objects of v.b.~type $(L,\cE),(L',\cE'),$ we consider the Floer complexes 
$CF^*((L,\cE),(L',\cE');H_n)$ whose generators are indexed by the set
$\cX(L,L';H_n)$ of time 1 trajectories of the
Hamiltonian vector field of $H_n$ starting at $L$ and ending at $L'$, or
equivalently, intersections of $\phi^1_{H_n}(L)$ with $L'$. 
For any value of $n$, we can use $H_n$ instead of $\varepsilon h$ in the
construction of Section \ref{s:Amodel}, and arrive at an $A_\infty$-category
$\cF(M;H_n)$ which is quasi-equivalent to $\cF(M)$. However, due to the lack of {\em a
priori} bound on the degrees of propagating discs with given inputs,
H.\ Lee's argument \cite{Lee} does not allow us to conclude that
the $A_\infty$-operations $\mu^k_{H_n}$ can be understood from local
considerations for any finite value of $n$, even if we restrict ourselves to
a finite set of objects.

To address this, we define $CF^*((L,\cE),(L',\cE');H)$ to be a completion of
the countably infinite dimensional $K$-vector space whose generators correspond to
(morphisms between fibers of the local systems at the end points of) time 1
trajectories of the Hamiltonian vector field of $H=\varepsilon h+\frac{A_e}{4}f(\tau)$
in the complement of the circles $\tau=1$ and $\tau=3$ where the flow is not defined.
Namely, $CF^*(L,L';H)$ consists of formal sums of elements
$\rho_p\in\Hom(\cE_{p},\cE'_{p'})$ such that $\|\rho_p\|\to 0$ (i.e., $\val(\rho_p)\to+\infty$).

The Floer complexes
$CF^*((L,\cE),(L',\cE');H_n)$ stabilize as $n\to \infty$, in the following
sense. For $\delta>0$, let $$\cN_{\delta}=\bigcup_{e\in E(G)}
\cN_{e,\delta},\ \text{where}\ \cN_{e,\delta}=([1-\delta,1+\delta]\cup
[3-\delta,3+\delta])\times S^1\subset M_e.$$
Then the generators of
$CF^*((L,\cE),(L',\cE');H_n)$ which lie outside of $\cN_{1/N}$ remain
exactly the same for all $n\geq N$, and so we can think of
$CF^*((L,\cE),(L',\cE');H)$ as the completion of the naive limit of these Floer complexes.
(Because our counts of discs are weighted by symplectic area rather than by topological
energy, we can directly identify Floer generators with each other for large
values of $n$, without the action rescaling discussed in \cite{AuSm}).

Considering the effect of the rotations $\psi\mapsto \psi+f'_n(\tau)$
induced by the perturbations, we see that, under mild assumptions on the geometry of $L$ and $L'$ near $\tau=1$ and
$\tau=3$, the set $\cX(L,L';H_n)$ (resp.\ $\cX(L,L';H)$) differs from $\cX(L,L';\varepsilon h)$
by adding:
\begin{itemize}
\item $n$ (resp.\ infinitely many) degree 1 generators in $(0,1)\times S^1;$

\item $2n$ (resp.\ ``twice'' infinitely many) degree 0 generators in $(1,3)\times S^1;$

\item $n$ (resp.\ infinitely many) degree $1$ generators in $(3,4)\times S^1.$ 
\end{itemize}

Fix Lagrangians $L_0,\dots,L_k$ of v.b.~type and input generators
$p_i\in \cX(L_{i-1},L_i;H_n)$. 
Consider a component $u_e:D_e\to M_e$ of a propagating perturbed
holomorphic disc for the Hamiltonian $H_n$ which maps to $M_e$, and assume
that $([0,1)\times S^1/\sim) \subset M_e$ contains part of the
image of $u_e$, but not its output. Then the lift of $u_e$ to the universal cover of
$M_e-\{p_v,p_{v'}\}$ has a maximum ``width'' along the $\psi$ direction which is
determined by the inputs of $u_e$ and, for those inputs which map to
the node $p_v$ at $\tau=0$, the local degree of $u_e$ in the strip-like end
near the node. However, the perturbation $H_n$ prevents any 
portion of $u_e$ of width less than $n$ from crossing $\tau=1$ in the
increasing $\tau$ direction from input to output. Therefore, if $n$ is
sufficiently large compared to the sum of the local degrees of $u_e$ at the inputs
which map to $p_v$, we arrive at a contradiction, and the output of $u_e$
must also lie at $\tau<1$; see \cite[Lemma 3.5]{Lee} (see also
\cite[Proposition 5.5]{AuSm}). 

We arrive at the following conclusion. For each vertex $v$ of $G,$ we denote by $P_v$ the union of subsets $([0,3]\times S^1/\sim)\subset M_e$ for each half-edge $e/v.$
For each edge $e$, denote by $N_e\subset M_e\subset M$ the subset $[1,3]\times S^1$.
(Thus, when $e$ is the only edge connecting $v$ to $v'$, $N_e=P_v\cap P_{v'}$).
Then:

\begin{prop}\label{prop:hlee}
Given any propagating perturbed holomorphic disc $u:S\to M$ for the
Hamiltonian $H_n$, with boundary on $L_0,\dots,L_k$ and inputs $p_i\in
\cX(L_{i-1},L_i;H_n)$, one of the following holds:
\begin{itemize}
\item the image of $u$ is entirely contained inside $P_v$ for some $v\in V(G)$,
and the output marked point lies outside of $N_e$ for all $e/v$;
\item the image of $u$ is entirely contained inside $N_e$ for some $e\in
E(G)$;
\item at least one of the input generators $p_i$ lies within $\cN_{k/n}$;
\item the disc $u$ propagates through a node of $M$ with output
degree $k_{out}\geq n/k$.
\end{itemize}
\end{prop}

In the last case, propagation with output degree $\geq n/k$ implies that the
symplectic area of the disc is bounded below by a constant
multiple of $n$. Therefore, we have: 

\begin{prop} \label{prop:lowarealocal}
For a given collection of input generators $p_i\in
\cX(L_{i-1},L_i;H)$ and a constant $A>0$, there exists $N=N(A)$ such that,
for $n\geq N$, any propagating perturbed holomorphic disc with inputs
$p_1,\dots,p_k$ and with area $\leq A$ lies entirely within a single piece
$P_v$ (or $N_e$), and its output lies outside of $\cN_{1/N}$.
Moreover, the moduli spaces of such discs are in bijection with each other for all
$n\geq N$.\end{prop}

\proof The first part of the statement is immediate from Proposition \ref{prop:hlee}, 
since for $n$ sufficiently large the area bound precludes propagation with output degree $\geq n/k$.
Moreover, the bound on propagation degrees implies a bound on the ``width'' of
each component of the propagating disc along the $\psi$ coordinate, and
hence for the output as well, whereas 
the generators near $\tau=1$ and $\tau=3$ correspond to trajectories of the
Hamiltonian flow which wrap more and more around the $S^1$ direction.
Finally, the existence of a bijection between the moduli spaces
of propagating discs for different values of $n\geq N$ is immediate for
discs which do not cross $\tau=1$; for those which cross $\tau=1$ (necessarily in the decreasing
$\tau$ direction from input to output), recasting
solutions to the perturbed Cauchy-Riemann equation as polygons with boundary on
the images of  $L_i$ under the Hamiltonian flow makes it clear that
increasing the value of $n$ simply deforms these polygons by widening
the strip-like portions that cross the neck at $\tau=1$.
(See also \cite[Section 3]{Lee} and \cite[Section 5]{AuSm} for related
arguments.)
\endproof

This allows us to define $A_\infty$-operations in $\F(M;H)$ as the naive
limits of the operations using Hamiltonians $H_n$: given $p_i\in
\cX(L_{i-1},L_i;H)$ and unitary $\rho_i\in \hom(\cE_{i-1|p_i},\cE_{i|p'_i})$, we
define
\begin{equation}\label{eq:mu_k_H}
\mu^k_H(\rho_k,\dots,\rho_1)=\lim_{n\to \infty} \mu^k_{H_n}(\rho_k,\dots,\rho_1),
\end{equation}
i.e.\ the element of $CF^*((L_0,\cE_0),(L_k,\cE_k);H)$ which agrees mod
$T^A$ with $\mu^k_{H_n}(\rho_k,\dots,\rho_1)$ for all $n> N(A)$.
We then extend this definition to finite sums of generators by linearity and then
to arbitrary inputs in the completed morphism spaces by continuity.

Concretely, $\mu^k_H(\rho_k,\dots,\rho_1)$ can be understood as a weighted
count of propagating discs in which the Hamiltonian perturbations are chosen
to be large enough relative to the given inputs and to the local degrees $k_{out}$
at the nodes of $S$; by Proposition \ref{prop:hlee} these discs remain
within a single $P_v$, and so the disc can only propagate through one node of $M.$

One small technical comment is in order: in the above construction we have
defined $A_\infty$-operations using the same Hamiltonian $H_n$ for the inputs and output 
of $\mu^k_{H_n}$, which means for $k\geq 2$ the perturbed Cauchy-Riemann equation involves
a 1-form $\beta$ that is not closed (for compact $M$ this is not a problem, since $H_n$ is
bounded; in the wrapped setting one should instead appeal to Abouzaid's rescaling trick on the noncompact components of $M$).
However one could also have used as in \cite{Lee} and \cite{AuSm} a closed 1-form in the Cauchy-Riemann equation and have
$\mu^k$ take values in a Floer complex with the Hamiltonian perturbation $kH_n$,
whose geometric behavior is essentially the same as that of $H_{kn}$.  The
details of the construction of the limit for $n\to \infty$ are then
different (and potentially more involved if one introduces a ``telescope'' model
for the chain-level limit of complexes for different Hamiltonians), but even
then it is possible under mild geometric assumptions on the Lagrangians
$L_i$ to rephrase the construction in terms of a (completed) naive limit.

We note the following consequence of Proposition
\ref{prop:hlee}, which we will use in Section \ref{s:equivalence}:

\begin{prop}\label{prop:ainfty_ideal}
For each $v\in V(G)$, the (completed) span of the generators of the Floer complexes
which lie outside of $P_v$ form an $A_{\infty}\!$-ideal in $\cF(M;H).$ 
We denote by $\cF(P_v;H)$ the quotient of $\cF(M;H)$ by this $A_{\infty}\!$-ideal.
Similarly, for each edge $e$ the span of the generators which lie outside of $N_e$ 
form an $A_{\infty}\!$-ideal in $\cF(M;H).$ We denote by $\cF(N_e;H)$ the quotient 
of $\cF(M;H)$ by this $A_{\infty}\!$-ideal. 
\end{prop}

\subsection{Continuation $A_\infty$-homomorphisms}\label{ss:continuation}
We end this section with the construction of $A_\infty$-homomorphisms from $\cF(M)$ to $\cF(M;H)$
via continuation maps in Lagrangian Floer theory (see e.g.\ \cite{SeBook});
because our comparison argument relies on a different approach (see Section \ref{ss:HPLFloer}), 
we skip some of the details involved in the construction of the higher terms. 

We construct an $A_\infty$-homomorphism $\mathfrak{K}_n:\cF(M)\to \cF(M;H_n)$, whose $k$-th order
term
$$\mathfrak{K}_n^k:\bigotimes_{i=1}^k CF^*((L_{i-1},\cE_{i-1}),(L_i,\cE_i);\varepsilon h)\to 
CF^{*+1-k}((L_0,\cE_0),(L_k,\cE_k);H_n)$$
counts perturbed propagating holomorphic discs with $k$ inputs,
for a Hamiltonian perturbation which interpolates between
$\varepsilon h$ at the inputs and $H_n$ at the output. 

The first order map $\mathfrak{K}_n^1$
is the easiest one to describe. Fix a smooth family of Hamiltonians
$H_\sigma$, $\sigma\in \R_{\geq 0}$ such that $H_\sigma=\varepsilon h$ for
$\sigma=0$ and $H_\sigma=H_n$ for $\sigma=n$; also fix a smooth
nonincreasing function $\sigma:\R\to \R$ such that $\sigma=n$ on
$(-\infty,-1)$ and $\sigma=0$ on $(1,\infty)$. The domain of a propagating disc
with a single input is a linear chain of discs with two marked points each
(i.e., strips $\R\times [0,1]$), $S=D_1\cup\dots\cup D_\ell$ (with $D_1$ carrying the input marked point
$z_1$ and $D_\ell$ carrying the output $z_0$). We then consider maps $u:S\to M$ in
which one of the components $D_j$ solves the usual Floer continuation equation
$$(du-X_{H_{\sigma(s)}}\,dt)^{0,1}=0$$
with Hamiltonian $\varepsilon h$ at the input $(s\to +\infty)$ and $H_n$ at
the output $(s\to -\infty)$, while the components $D_1,\dots,D_{j-1}$ 
(resp.\ $D_{j+1},\dots,D_\ell$) which precede (resp.\ follow) it along the
way from the input to the output are perturbed holomorphic strips for the Hamiltonian
$\varepsilon h$ (resp.\ $H_n$). Counting such perturbed propagating discs
(for all possible choices of the component of $S$ where continuation takes
place) which are rigid  (i.e., belong to moduli spaces of solutions with expected dimension 0, or
equivalently, the input and output generators have the same degree), with
signs and weights as in the definition of the $A_\infty$-operations,
yields the map $\mathfrak{K}_n^1:CF^*((L_0,\cE_0),(L_1,\cE_1);\varepsilon h)\to
CF^*((L_0,\cE_0),(L_1,\cE_1);H_n)$, which is easily checked to be a chain
map by considering one-dimensional moduli spaces.

\begin{remark}\label{rmk:continuation}
Although the definition allows the change of Hamiltonian
to happen in any component of the propagating disc, the components of a regular 
rigid continuation trajectory are themselves rigid; this implies that the continuation
must actually take place in the input component $D_1$, resp.\ the output
component $D_\ell$, if the input and output 
are degree 1, resp.\ degree 0 generators of the respective Floer complexes.

In fact, in our setting, continuation
trajectories starting at a degree 1 generator in the interior of $M_e$ are
necessarily constant. Therefore, $\mathfrak{K}_n^1$ is the naive inclusion
on $CF^1$, while for degree 0 generators it differs from
the naive inclusion (constant trajectories) by counts of propagating
perturbed holomorphic strips in which the output component is a continuation trajectory
from $\varepsilon h$ to $H_n$ in the usual sense and all other components
are perturbed holomorphic strips for the Hamiltonian $\varepsilon h$.
\end{remark}

Moreover, the same arguments as in the previous section show that
continuation trajectories of bounded symplectic area (or energy --- the two
are equivalent because $f_n$ and $f$ are uniformly bounded), hence bounded
propagation degrees through the nodes of $M$, must stabilize as $n\to
\infty$, i.e.\ the moduli spaces are the same for all sufficiently large
values of $n$.  This allows us to define
$\mathfrak{K}^1:CF^*((L_0,\cE_0),(L_1,\cE_1);\varepsilon h)\to
CF^*((L_0,\cE_0),(L_1,\cE_1);H)$ by $\mathfrak{K}^1(\rho)=\lim_{n\to \infty}
\mathfrak{K}^1_n(\rho)$.  Taking the limit $n\to \infty$ in the identity
$\mu^1_{H_n}\circ \mathfrak{K}^1_n=\mathfrak{K}^1_n\circ \mu^1$ shows that $\mathfrak{K}^1$ is also
a chain map.

It follows from Remark \ref{rmk:continuation} that, for degree 0 Floer
generators, $\mathfrak{K}^1$ counts propagating discs in which the output component
of $S$ is a continuation trajectory from $\varepsilon h$ to $H$ (i.e.\ $H_n$
for sufficiently large $n$), while the other components are solutions to
Floer's equation for the Hamiltonian $\varepsilon h$; whereas for
degree 1 generators continuation happens at the input.%
\medskip

The higher order terms of the $A_\infty$-homomorphisms $\mathfrak{K}_n$ involve the
choice, for each stable nodal domain $S=\bigsqcup D_j/\sim$ (and continuously
and consistently over the moduli space of these), of a one-parameter
family of Hamiltonian perturbation data, such that at one end of the family
the Hamiltonian is $\varepsilon h$ everywhere except in the strip-like
end near the output marked point $z_0$ where the continuation to $H_n$ takes
place, and at the other end of the family the Hamiltonian is $H_n$ everywhere except in the
strip-like ends near the input marked points $z_1,\dots,z_k$.
One way to achieve this is to choose for each $S$ a 
smooth function $s:S-\{z_0,\dots,z_k\}\to \R$ such that 
$\lim_{z\to z_0} s(z)=-\infty$ at the output marked point, $\lim_{z\to z_i}
s(z)=+\infty$ at the input marked points, and on each component of $S,$
$s$ decreases monotonically from the inputs to the output. This
choice should be made continuously
over the moduli space of stable nodal discs and consistently with respect to
degenerations of the domain. We then consider solutions of the Floer
continuation equation involving the Hamiltonians $H_{\sigma\left(s(z)-s_0\right)}$, where
the parameter $s_0\in \R$ is allowed to vary and determines the level set of
$s$ near which the Hamiltonian perturbation switches from $\varepsilon h$ to
$H_n$.

Since $H_n=\varepsilon h$ near the nodes of $M$ (and we can ensure that the
same holds for all $H_\sigma$), the details of the behavior of the continuation
perturbation as $s_0$ varies through the value of $s$ at a node of $S$ are
not particularly important. What does require more care is the case where
some components of $S$ are unstable (strips), and the most obvious 
constructions fail to account for domain automorphisms if continuation
proceeds simultaneously across several unstable components of $S$.
Conceptually the simplest approach is to stabilize the domain by adding a boundary marked
point to each unstable component of $S$, where we require the
$\tau$-coordinate of the appropriate component $M_e$ to take a prescribed
value. (Alternatively, by considering the structure of rigid continuation
configurations as in Remark \ref{rmk:continuation} one can exclude a number
of potential cases and devise an ad hoc definition for the remaining ones).

As in the case of the linear term, observing that contributions to
$\mathfrak{K}_n^k$ from propagating discs whose area is below a fixed threshold
stabilize for sufficiently large $n$, we can take the limit as $n\to\infty$
and set $\mathfrak{K}^k(\rho_k,\dots,\rho_1)=\lim_{n\to \infty}
\mathfrak{K}^k_n(\rho_k,\dots,\rho_1)$.

We claim that the $A_\infty$-functor $\mathfrak{K}:\cF(M)\to \cF(M;H)$ is a quasi-equivalence.
The usual method to establish such a result is to construct another $A_\infty$-functor in the opposite direction by considering
Floer-theoretic continuation maps with the roles of $H$ and $\varepsilon h$
reversed, and show that it is a quasi-inverse to $\mathfrak{K}$ by a homotopy
argument. We expect that this can be done in our setting,
but it is easier to proceed differently. Namely, it suffices to show that
the linear terms of the $A_\infty$-functor $\mathfrak{K}$ are quasi-isomorphisms of chain
complexes; this will follow from the argument in Section \ref{ss:HPLFloer}
where we show that $\mathfrak{K}^1$ coincides with a purely
algebraic construction based on the homological perturbation lemma.

\section{The B-model: generalized Tate curves from combinatorial data} \label{s:Bmodel}

\subsection{Generalized Tate curve in terms of formal schemes}

Given combinatorial data as in Definition \ref{def:combidata},
the following is a particular case of Mumford's construction (actually, its version over the universal power series ring). We first take the $\Z$-scheme $X^0,$ which is obtained as a union of $X^0_{v},$ where we identify $x_{e/v}$ and $x_{e/v'}$ for $v\ne v'.$ The resulting nodal points are denoted by $x_{e}\in X^0.$

Let us choose the following affine open subsets $U^0_{e},W^0_{v}\subset X^0.$ For $v\in V,$ the subset $W^0_{v}$ is ($X^0_{v}$ minus nodal points). For $e\in E$ we take $v,$ $v'$ to be the endpoints of $e,$ and define $U^0_{e}$ to be $X^0_{v}\cup X^0_{v'}$ minus nodal points other than $x_{e}.$ We have isomorphisms
$$W_{v}^0\cong\Spec\Z[t^{\pm 1},(1-t)^{-1}],$$
$$U_{e}^0\cong \Spec \Z[t_{e/v},(1-t_{e/v})^{-1},t_{e/v'},(1-t_{e/v'})^{-1}]/(t_{e/v}t_{e/v'}).$$
The first of these isomorphisms of course depends on a choice of coordinate $t$ on $X_{v}^0$ taking values $0,1,\infty$ at the marked points. 

We now define the formal scheme $\mX$ over $\Z[[q_{e},e\in E]].$ Its reduction modulo all $q_{e}$ will be exactly $X^0.$ We first take the affine formal schemes
$\cU_{e},\cW_{v},$ given by
$$\cW_{v}:=\Spf \cO(W_{v}^0)[[q_{f},f\in E]];$$
$$\cU_{e}:=\Spf\Z[T_{e/v},(1-T_{e/v})^{-1},T_{e/v'},(1-T_{e/v'})^{-1}][[q_{f},f\in E]]/(T_{e/v}T_{e/v'}-q_{e}).$$
It is easy to see that for $e/v,$ $e/v',$ we have a natural isomorphism
$$\widehat{\cO(\cU_{e})[T_{e/v}^{-1}]}\xto{\sim}\cO(\cW_{v}),\quad T_{e/v}\mapsto t_{e/v},\, T_{e/v'}\mapsto \frac{q_{e}}{t_{e/v}}.$$ 

This allows us to glue together all $\cU_{e}$ in the obvious way, and this way we obtain our formal scheme $\mX.$ It is easy to see from Grothendieck algebraization theorem that there is a unique (up to canonical isomorphism) algebraic curve $X$ over $\Z[[\{q_{e}\}]]$ such that the reduction of $X$ mod $q_{e}$ is identified with $X^0,$ and the formal neighborhood of $X^0$ at $X$ is identified with $\mX.$

However, the algebraization is essentially impossible to write down explicitly, and we don't need that since the categories of coherent sheaves and of perfect complexes are naturally obtained from the formal scheme. That is, we have $\Coh(X)\simeq\Coh(\mX),$ $\Perf(X)\simeq \Perf(\mX).$

\begin{remark}\label{remark:inverting_q} Although in general punctured formal schemes (objects like $\mX-X^0$) are not easy to deal with, here they are not too much different
from usual schemes. Namely, if we want to invert some collection $q_{e_1},\dots,q_{e_l}$ (for example, all $q_{e}$'s), then we simply take a ringed space $\mX'$ with the same underlying topological space, and define the sections on affine subsets by 
$$\cO_{\mX'}(\cU)=\widehat{\cO_{\mX}(\cU)[(q_{e_1}\dots q_{e_l})^{-1}]}.$$
Then we will have $\Coh(\mX')=\Coh(\mX)/(q_{e_1}\dots q_{e_l}\text{-torsion}),$ and similarly for $\Perf(\mX').$
\end{remark}

From now on, we denote by $K$ the Novikov field $\hhat{\mk[T^{\R}]},$ where $\mk$ is some field of coefficients. As above, we denote by $A_e\in\R_{>0}$ the symplectic areas of $2$-spheres on the A side. Taking continuous homomorphism $$\Z[[\{q_e, e\in E\}]]\to K,\quad q_e\mapsto T^{A_e}$$ 
(or some other element of valuation $A_e$ if we allow a bulk deformation of the A-model),
we get the extension of scalars $X_K$ of $X.$ The B side will be the curve $X_K.$

\subsection{The Schottky groupoid}

We now give the description of the curve $X_K$ in terms of rigid analytic geometry. To avoid confusion, we put $Y_v:=X^0_v\times_{\Z} K\cong\PP^1_K,$ and keep the notation $t_{e/v}$ for the chosen projective coordinates.

We denote by $\pi_1(G)$ the fundamental groupoid of the graph $G.$ We define the functor $g:\pi_1(G)\to \Sch/K$ by sending each $v\in V$ to $Y_v,$ and for each edge $e$ connecting $v$ and $v'$ we send the morphism $e:v\to v'$ to the map $g_{e/v}:X^0_v\to X^0_{v'},$ given by $t_{e/v'}(g_{e/v}(x))=\frac{q_e}{t_{e/v}(x)}.$

Fixing a vertex $v_0\in V,$ we get the {\it Schottky group} $\Gamma_{G,v_0}:=\pi_1(G,v_0),$ which acts faithfully on $Y_{v_0}.$ The group $\Gamma_{G,v_0}$ is free on $g=g(X_K)$ generators, and its non-identity elements are acting by hyperbolic transformations of $Y_v.$ 

If we now consider each $Y_v$ as a rigid analytic space, then we define $F_v\subset Y_v$ to be the set of limit points of the $\pi_1(G,v)$-action ($F_v$ is naturally a Cantor
set; see e.g.\ \cite[Chapter I]{Gerritzen-VdPut}).
Then the curve $X_K$ is identified, as a rigid analytic space, with the quotient of the collection $\{Y_v-F_v\}_{v\in V}$ by $\pi_1(G)$ (the same as the quotient $(Y_v-F_v)/\Gamma_{G,v}$ for any $v\in V$).

For each vertex $v\in V,$ and any real numbers $1>s_{e/v}>|q_e|,$ for each half-edge $e/v,$ we define the open analytic subset $U_{v,\{t_{e/v}\}}\subset X_K$ as the image of $$\{1\geq |t_{e/v}|\geq s_{e/v}\text{ for }e/v\}\subset Y_v-F_v.$$ Also, for any half-edge $e/v,$ and any $1>s_1\geq s_2>|q_e|$ we define $U_{e/v,s_1,s_2}$ as the image of $$\{s_1\geq |t_{e/v}|\geq s_2\}\subset Y_v-F_v.$$ Clearly, if the edge $e$ connects $v$ and $v'$ then $U_{e/v,s_1,s_2}=U_{e/v',\frac{|q_e|}{s_2},\frac{|q_e|}{s_1}}.$ For two distinct $v,v'\in V,$ and collections $\{s_{e/v}\},\{s_{e/v'}\},$ we have
$$U_{v,\{s_{e/v}\}}\cap U_{v',\{s_{e/v'}\}}=\bigsqcup_{e}U_{e/v,\frac{|q_e|}{s_{e/v'}},s_{e/v}},$$
where the union is over the edges connecting $v$ and $v',$ and we put $U_{e/v,s_1,s_2}=\emptyset$ if $s_1<s_2.$

We will mostly use the following open affinoid subsets:
$$U_v:=U_{v,\{|q_e|^{\frac34}\}_{e/v}},\quad U_e:=U_{e/v,|q_e|^{\frac14},|q_e|^{\frac34}}.$$





\section{Construction of the equivalence}\label{s:equivalence}

\subsection{The assignment of vector bundles to objects of $\cF(M)$} \label{ss:assignbundles}

Recall that a v.b. type object of $\cF(M)$ is a pair $(L,\cE),$ where $L$ is a graph with vertices in $V(G)$ and edges going in each of $M_e,$ and $\cE$ is a local system of free finitely generated $\cO_K$-modules on $L.$ We fix such a Lagrangian graph $L_0,$ so that the pair $(L_0,\cO_K)$ will correspond to the structure sheaf $\cO_X.$

Now, for any object $(L,\cE)\in\cF(M)$ and each edge $e\in E(G)$ connecting $v,v'\in V(G)$ we have the following invariants:

\begin{itemize}
\item $r_e(L)=r_e(L_0,L)\in \Z,$ the rotation number of $L$ with respect to $L_0$ in $M_e$ in the negative direction. 
The sum $\smash{\sum\limits_{e\in E(G)}}r_e(L)$ will be the slope of the corresponding vector bundle.

\item $S_{e/v}(L)=S_{e/v}(L_0,L),$ the signed area bounded by $L_0$ and $L$ on the universal cover of 
$M_e\setminus\{p_v,p_{v'}\},$ where we take the lifts which are close to each other when we approach $p_{v'}.$ 
We have $S_{e/v'}(L)+S_{e/v}(L)=r_{e}(L)A_e.$

\item the monodromy $R_{\cE,e/v}:\cE_v\to\cE_{v'}.$
\end{itemize}

We define the vector bundle $\Phi(L,\cE)$ on $X$ as follows. First, its pullbacks to $Y_v-F_v$ are given by $\cE_v\otimes_{\cO_K}\cO_{Y_v-F_v}.$ Then, we need to describe the action of the groupoid $\pi_1(G).$ 
For each edge $e$ considered as a morphism from $v$ to $v'$ in $\pi_1(G),$ the corresponding isomorphism $$u_{e/v}:\cE_v\otimes\cO_{Y_v-F_v}\to g_{e/v}^*(\cE_{v'}\otimes\cO_{Y_{v'}-F_{v'}})$$ is given by 
\begin{equation}\label{eq:gluingdata} u_{e/v}=R_{\cE,e/v}\otimes T^{-S_{e/v}(L)}t_{e/v}^{r_e(L)}.\end{equation}
If the A-model is bulk deformed, this formula should be
corrected by the exponential of the integral of $\mathfrak{b}$ over the area
bounded by $L_0$ and $L$ on the universal cover of $M_e\setminus \{p_v,p_{v'}\}$.

By \cite{Fa}, the vector bundle $\Phi(L,\cE)$ is semistable of slope $\sum\limits_{e\in E(G)}r_e(L).$

\subsection{Abstract Homological Perturbation Lemma (HPL) for complexes}
\label{ss:HPL}

We recall the following abstract setup for homological perturbation, for which we refer to \cite{CL}. For simplicity the base field will be the Novikov field $K.$ 

Let $(\cK,d_{\cK})$ and $(\cL,d_{\cL})$ be complexes. Suppose that we are given  maps $i,p,h,$ where $i:\cL\to\cK$ and $p:\cK\to \cL$ are morphisms of complexes, and $h:\cK\to\cK$ is a map of (cohomological) degree $-1$ such that $pi=1_{\cL},$ $1_{\cK}-ip=dh+hd,$ $h^2=0,$ $ph=0,$ $hi=0.$ 

Now let us take a perturbation $\delta$ of the differential $d_{\cK},$ satisfying the Maurer-Cartan equation $[d_{\cK},\delta]+\delta^2=0.$ Hence, $\tilde{d}_{\cK}=d_{\cK}+\delta$ is a differential on $\cK.$ Assume that the endomorphism $(\id_{\cK}+h\delta)$ of $\cK$ is invertible (hence, so is $(\id_{\cK}+\delta h)$). Then there are natural perturbations for $d_{\cL},$ $i,$ $p$ and $h,$ so that all of the relations continue to hold:
$$\tilde{d}_{\cL}=d_{\cL}+p\delta(\id+h\delta)^{-1}i,\quad \tilde{i}=(\id+h\delta)^{-1}i,\quad \tilde{p}=p(\id+\delta h)^{-1},\quad \tilde{h}=(\id+h\delta)^{-1}h.$$

In particular, $\tilde{i}$ and $\tilde{p}$ are quasi-isomorphisms of complexes with perturbed differentials. 

\begin{remark}\label{remark:top_nilpotent}Suppose that (the graded components of) $\cK$ and $\cL$ are Banach vector spaces over $K,$ and the maps $i,p,h,\delta$ are continuous. Then the assumption that $(\id+h\delta)$ is invertible would follow from the assumption that $h\delta:\cK\to\cK$ is locally topologically nilpotent, i.e. for any homogeneous $x\in\cK$ we have $\lim\limits_{n\to\infty}(h\delta)^n(x)=0.$ Indeed, in this case we have $$(\id+h\delta)^{-1}(x)=\sum\limits_{n=0}^{\infty}(-1)^n(h\delta)^n(x),\quad (\id+\delta h)^{-1}(x)=\sum\limits_{n=0}^{\infty}(-1)^n(\delta h)^n(x).$$ Hence, the formulas for $\tilde{d}_{\cL},\tilde{i},\tilde{p},\tilde{h}$ can also be expressed as infinite sums.\end{remark}

Note that for $(\cK,\cL,i,p,h)$ and $(\cK',\cL'\,i',p',h')$ as above one can define
their tensor product to be $(\cK\otimes\cK',\cL\otimes\cL',i\otimes i', p\otimes p',h''),$ where $h''=h\otimes\id+ip\otimes h'.$   

Now suppose that $(\cK,\cL,i,p,h)$ is as above and $\mu_{\cK}=(\mu_{\cK}^1,\mu_{\cK}^2,\dots)$ is an $A_{\infty}$-structure on $\cK.$ Then we get the data $(T(\cK[1]),T(\cL[1]),i',p',h')$ as above, using the formulas for the tensor product. We get a coderivation $\delta:T(\cK[1])\to T(\cK[1])$ of degree $1,$ with components $\delta^1=\mu^1_{\cK}-d_{\cK},\delta^2=\mu^2_{\cK},\dots.$ Assuming that $(\id_{\cK}+h\delta^1)$ is invertible, we easily see that the same holds for $(\id_{T(\cK[1])}+h'\delta).$ Applying the above formulas, we get the deformed differential $\tilde{d}_{T(\cL[1])},$ which is in fact a coderivation, hence it gives an $A_{\infty}$-structure $\mu_{\cL}$ on $\cL.$ The deformed morphisms $\tilde{i'}:T(\cL[1])\to T(\cK[1]),$ $\tilde{p'}:T(\cK[1])\to T(\cL[1])$ are in fact morphisms of DG coalgebras, hence they give morphisms of $A_{\infty}$-algebras $\tilde{i}:(\cL,\mu_{\cL})\to (\cK,\mu_{\cK}),$ $\tilde{p}:(\cK,\mu_{\cK})\to (\cL,\mu_{\cL}),$ which are quasi-isomorphisms. For details, see \cite[Section 3.3]{CL}.

Note that if we are in the setup of Remark \ref{remark:top_nilpotent}, then the expressions of $\mu_{\cL},\tilde{i},\tilde{p}$ as infinite sums are actually the standard summations over trees. The summations for the components $\mu_{\cL}^n,$ $\tilde{i}_n,$ $\tilde{p}_n$ would be finite if $\delta^1=0.$

The same construction applies also to $A_{\infty}$-categories. We will use it
in Section \ref{ss:HPLFloer} to argue that
the variant $\F(M;H)$ of our A-model construction involving
``infinite'' Hamiltonian perturbations is quasi-isomorphic to $\F(M)$. 
We will also use it to obtain expressions for the theta functions corresponding to the generators of the Floer complexes (provided that they are concentrated in degree zero). 

\subsection{Infinite Hamiltonian perturbations and \v{C}ech complexes}

Recall from Proposition \ref{prop:ainfty_ideal} that, denoting by $P_v$ the
union of the subsets $([0,3]\times S^1/\sim)\subset M_e$ for each half-edge
$e/v$, and by $N_e$ the subset $[1,3]\times S^1\subset M_e$, the generators
of the Floer complexes which lie outside of $P_v$ (resp.\ $N_e$) span (after
completion) an $A_\infty$-ideal in $\cF(M;H)$, and we denote by $\cF(P_v;H)$
(resp.\ $\cF(N_e;H)$) the quotient of $\cF(M;H)$ by this $A_\infty$-ideal.

These quotients come with $A_\infty$-functors $\cF(M;H)\to \cF(P_v;H)$ and $\cF(P_v;H)\to
\cF(N_e;H)$, which are surjective on morphisms and have vanishing higher order
terms. Hence, the naive chain level limit embeds fully faithfully into the homotopy limit $$\lim(\prod\limits_{v\in V(G)}\cF(P_v;H)\prarr\prod\limits_{e\in E(G)}\cF(N_e;H))\hto \holim(\prod\limits_{v\in V(G)}\cF(P_v;H)\prarr\prod\limits_{e\in E(G)}\cF(N_e;H)),$$ 
and the Fukaya category $\cF(M;H)$ embeds fully faithfully into the naive limit:
$$\cF(M;H)\hto \lim(\prod\limits_{v\in V(G)}\cF(P_v;H)\prarr\prod\limits_{e\in E(G)}\cF(N_e;H)).$$ 
We will show in Section \ref{ss:localfunctors} that there are natural equivalences
$$\Perf(\cF(P_v;H))\simeq \Perf(U_v),\quad \Perf(\cF(N_e;H)\simeq \Perf(U_e),$$ under which the functors $\cF(P_v;H)\to \cF(N_e;H)$ correspond to the restriction functors $\Perf(U_v)\to \Perf(U_e)$. 
Thus, we get a fully faithful embedding
$$\cF(M;H)\to\holim(\prod\limits_{v\in V(G)}\Perf(U_v)\prarr\prod\limits_{e\in E(G)}\Perf(U_e))\simeq \Perf(X_K),$$
which induces a fully faithful functor $\Psi:\Perf(\cF(M;H))\to \Perf(X_K)$. 
But on the level of isomorphism classes of objects this functor sends $(L,\cE)$ exactly to $\Phi(L,\cE)$. 
Since the vector bundles of the form $\Phi(L,\cE)$ generate the category $\Perf(X_K),$ we conclude that $\Psi$ is an equivalence.          

\subsection{HPL for Hamiltonian perturbations}
\label{ss:HPLFloer}

We now show that HPL provides a quasi-equivalence between the
$A_{\infty}$-categories $\cF(M)$ and $\cF(M;H)$; we
also explain how this can be viewed as an algebraic counterpart to
the continuation functor $\mathfrak{K}:\cF(M)\to \cF(M;H)$ described in
Section \ref{ss:continuation}.

As in Section \ref{ss:FMH}, for each edge $e$ connecting vertices $v$ and
$v'$ we identify $M_e$ with $[0,4]\times (\R/\Z)/\!\sim$ with coordinates
$(\tau,\psi)$ (with $\tau=0$ corresponding to the node $p_v$ and $\tau=4$ to
$p_{v'}$). Consider two v.b.-type Lagrangians $L,L'$. Without loss
of generality we assume that the generators of $CF^*(L,L';\varepsilon h)$ lie away from
the support of the perturbations $f_n$ and $f,$
and that $\cX(L,L';H_n)$ (resp.\ $\cX(L,L';H)$) differs from
$\cX(L,L';\varepsilon h)$ by adding, in each component $M_e$:
\begin{itemize}
\item $n$ (resp.\ infinitely many) degree 1 generators
$q_{e/v,1},q_{e/v,2},\dots$ (in increasing order of $\tau$ coordinates) in $(0,1)\times S^1;$

\item $2n$ (resp.\ ``twice'' infinitely many) degree 0 generators
$\dots,p_{e/v,2},p_{e/v,1}$ (near $\tau=1$) and
$p_{e/v',1},p_{e/v',2},\dots$ (near $\tau=3$) in $(1,3)\times S^1;$

\item $n$ (resp.\ infinitely many) degree $1$ generators
$\dots,q_{e/v',2},q_{e/v',1}$ in $(3,4)\times S^1.$ 
\end{itemize}

Denote by $\mu^1_{nv}$ the ``naive'' (or ``low area'') part of the differential 
$\mu^1_H$ on $CF^*(L,L';H)$, only involving holomorphic
discs supported near $\tau=1$ or $\tau=3$ (without propagation) in a single
component $M_e$ of $M$. Thus,
$\mu^1_{nv}$ maps
$p_{e/v,k}$ to a multiple of $q_{e/v,k}$ for every half-edge $e/v$ and for all $k\geq
1$, and all other generators to zero. The areas of the discs
connecting $p_{e/v,k}$ to $q_{e/v,k}$ can be made arbitrarily small by
shrinking the support of the perturbations $f_n$ and $f$; this allows us to
assume that all other contributions to the Floer differential $\mu^1_H$ have larger
area than those which are recorded by $\mu^1_{nv}$.

Setting $\delta^1=\mu^1_H-\mu^1_{nv}$ (and $\delta^k=\mu^k_H$ for $k\geq 2$), we are
now in the setup of abstract HPL. Namely, the natural inclusion 
$i:(CF^*(L,L'),0)\to (CF^*(L,L';H),\mu^1_{nv})$ is a map of complexes, and so is the 
projection $p:(CF^*(L,L';H),\mu^1_{nv})\to (CF^*(L,L'),0).$ 
Further, we choose the homotopy $h$ to be the map sending each new generator
of degree 1, $q_{e/v,k}$, to the corresponding degree zero generator
$p_{e/v,k}$, multiplied by the inverse of the coefficient that
arises in $\mu^1_{nv}$. 
Then the map $h\delta^1$ is locally topologically nilpotent, because the
symplectic areas of the
perturbed holomorphic discs which contribute to $\delta^1$ are larger than
those of the discs which contribute to $\mu^1_{nv}$. It follows that
$\id+h\delta$ is invertible (see Remark \ref{remark:top_nilpotent}), and we can apply HPL.

Applying this construction to the $A_\infty$-categories $\cF(M)$ and $\cF(M;H)$ (or rather,
to full subcategories whose objects satisfy the assumptions we have made above about the
absence of Floer generators near $\tau=1$ and $\tau=3$ and the behavior
of the Floer complexes under Hamiltonian perturbations), we arrive at the
existence of operations $\mu^k_{HPL}$ ($k\geq 1$) on the Floer complexes
$CF^*(L,L')$, given by the formulas in Section \ref{ss:HPL},
and $A_\infty$-functors $\tilde{i}$ and $\tilde{p}$ giving a
quasi-equivalence between this $A_\infty$-category and $\cF(M;H)$. 
\medskip

We now show that the operations $\mu^k_{HPL}$ obtained from
$\mu^k_H$ via Homological Perturbation theory are equal to the structure
maps $\mu^k$ of the Fukaya category $\cF(M)$, so that $\tilde{i}$ and
$\tilde{p}$ in fact yield a quasi-equivalence between $\cF(M)$ and $\cF(M;H)$.
We start with the differential, and recall that the HPL gives
\begin{equation}\label{eq:mu1_HPL}
\mu^1_{HPL}=p\delta^1(\id+h\delta^1)^{-1}i=\sum_{\ell=0}^\infty
(-1)^\ell p\delta^1(h\delta^1)^\ell i.
\end{equation}

Consider two v.b.-type Lagrangians $L,L'$ as above, and a propagating
holomorphic strip $u:S\to M$ contributing to the Floer differential on
$CF^*(L,L')$, connecting an input generator $p_1$ to an output generator
$p_0$ via a sequence of holomorphic strips contained successively in components
$M_{e_1},\dots,M_{e_\ell}$ (with $p_1\in M_{e_1}$ and $p_0\in M_{e_\ell}$), 
attached to each other via nodes $p_{v_1},\dots,p_{v_{\ell-1}}$. 
Since $u$ is rigid, its boundary travels along $L$ and $L'$ without
backtracking, and the $\tau$ coordinate varies monotonically along each
component. We orient each edge $e_j$ so that the strip travels in the increasing
$\tau$ direction along $M_{e_j}$ from input to output, i.e.\ $p_{v_j}$ lies at the $\tau=4$
end of $M_{e_j}$ and at the $\tau=0$ end of $M_{e_{j+1}}$.  Assume for
example that the $\tau$-coordinate of the input $p_1\in M_{e_1}$ is less than
1, and that the $\tau$-coordinate of the output $p_0\in M_{e_\ell}$ is
greater than 1, so that each component of $u$ passes through the circle
$\{1\}\times S^1\subset M_{e_i}$ (the other cases are similar). Denote by
$w_j\in \R_+$ the width of the $j$-th component of $u$ at $\tau=1$, i.e.\ the
difference in the values of the $\psi$ coordinate at $\tau=1$ on the two
boundaries of the lift of the strip to the universal cover of
$M_{e_j}-\{p_{v_{j-1}},p_{v_j}\}$, and let $k_j=\lceil w_j \rceil$. Then the Hamiltonian perturbation $H$ (or
$H_n$ for $n>\max(w_j)$) breaks each component of $u$ into a strip which
ends at the new degree 1 generator $q_{e_j/v_{j-1},k_j}$ before $\tau$
reaches 1, and one which starts from the new degree 0 generator
$p_{e_j/v_{j-1},k_j}$ just past $\tau=1$. Thus we can associate to $u$ a
sequence of $\ell+1$ perturbed propagating holomorphic strips contributing to
differential $\mu^1_H$ on $CF^*(L,L';H)$ (and hence to
$\delta^1=\mu^1_H-\mu^1_{nv}$), interspersed with $\ell$ low area
connecting trajectories between the pairs of generators $p_{e_j/v_{j-1},k_j}$
and $q_{e_j/v_{j-1},k_j}$. These are exactly the kinds of configurations
counted by the right-hand side of \eqref{eq:mu1_HPL}. Moreover, the
propagation multiplicity of $u$ is equal to the product of the propagation
multiplicities of the $\ell+1$ perturbed strips that it breaks into; its
area is the sum of the areas of these strips minus the sum of the areas of
the connecting trajectories between $p_{e_j/v_{j-1},k_j}$ and
$q_{e_j/v_{j-1},k_j}$, and similarly for holonomies. Finally, the sign
$(-1)^\ell$ is due to the overall sign contributions of the additional pairs of outputs at 
the new degree 1 generators $q_{e_j/v_{j-1},k_j}$ in the broken
configuration; each time the two trajectories
ending at $q_{e_j/v_{j-1},k_j}$ have opposite boundary orientations along
$L'$, so their signs differ by $-1$. It follows that $\mu^1_{HPL}=\mu^1$.

The argument for $\mu^k$, $k\geq 2$ is similar. Consider v.b.-type
Lagrangians $L_0,\dots,L_k$ which pairwise satisfy the simplifying
assumptions we have made about the behavior of the Floer complexes under
perturbation, and a rigid propagating holomorphic disc $u:S\to M$ with
boundary on $L_0,\dots,L_k$ which contributes to $\mu^k$. The intersection
of $u$ with a neighborhood $\cN_\delta$ of the circles $\{1\}\times S^1$ and
$\{3\}\times S^1$ in every component of $M$ is a union of strip-like
portions of the propagating disc. Among these, the strips which cross $\tau=1$
(resp.\ $\tau=3$) in the decreasing (resp.\ increasing) $\tau$ direction
are essentially unaffected by the Hamiltonian perturbations $H_n$, while
those which cross $\tau=1$ (resp.\ $\tau=3$) in the increasing (resp.\
decreasing) $\tau$ direction get broken up as described above as soon
as $n$ exceeds their width along the $\psi$ coordinate.
Thus, $u:S\to M$ gets broken into a union of perturbed
propagating discs contributing to the structure maps of $\cF(M;H)$, each of
them with inputs that are either inputs of $u$ (hence ``old'' generators from $\cX(L_{i-1},L_i;\varepsilon
h)$) or new degree $0$ generators $p_{e/v,k}$, and outputs that are either new degree 1
generators $q_{e/v,k}$ or the output of $u$. Because $h$ vanishes on all except new
degree 1 generators, which it maps to the corresponding new degree 0
generators, this type of configuration agrees exactly with the 
tree sum that appears in the HPL formula, and we conclude that
$\mu^k_{HPL}=\mu^k$.

This completes the proof that $\cF(M;H)$ is quasi-equivalent to $\cF(M)$
(via the $A_\infty$-functors $\tilde{i}$ and $\tilde{p}$ provided by HPL).

While not needed for our argument, it is also instructive to compare $\tilde{i}:\cF(M)\to
\cF(M;H)$ with the continuation functor $\mathfrak{K}$ described in \S
\ref{ss:continuation}. The HPL formula for the linear term is
$$\tilde{i}^1=(\id+h\delta^1)^{-1}i=\sum_{\ell=0}^\infty
(-1)^\ell (h\delta^1)^\ell i.$$
Since $h\delta^1$ vanishes on degree 1 generators, for $CF^1$ this
simplifies to the naive inclusion $i$. For degree 0 generators,
$\tilde{i}^1$ differs from the inclusion by counts of broken configurations consisting of
perturbed propagating holomorphic strips contributing to
differential $\mu^1_H$ (i.e., to $\delta^1=\mu^1_H-\mu^1_{nv}$), ending at
degree 1 generators $q_{e_j/v_{j-1},k}$, interspersed with (inverses of) low area
connecting trajectories between pairs of generators $q_{e_j/v_{j-1},k}$ and
$p_{e_j/v_{j-1},k}$. Arguing as above, such configurations correspond
almost exactly to propagating holomorphic discs for the Floer differential
$\mu^1$ (with Hamiltonian perturbation $\varepsilon h$), except for the 
component carrying the output, where the picture is different and can be
checked by explicit calculation to match the behavior of a Floer
continuation trajectory from the Hamiltonian perturbation $\varepsilon h$ to
the perturbation $H_n$ for $n$ sufficiently large. Comparing with the
description in Remark \ref{rmk:continuation}, we conclude that
$\tilde{i}^1=\mathfrak{K}^1$. This in turn implies that $\mathfrak{K}$ is a 
quasi-equivalence. We expect (but have not checked) that the higher terms
of the $A_\infty$-functors $\tilde{i}$ and $\mathfrak{K}$ can also be shown
to agree.

\subsection{The local functors}\label{ss:localfunctors}

We now describe the functors $\cF(P_v;H)\to \Perf(U_v)$ and $\cF(N_e;H)\to \Perf(U_e),$ which after gluing give the functor $\cF(M;H)\to\Perf(X_K).$ 
We start with $P_v.$ We send each v.b.-type object $(L,\cE)$ to the free sheaf $\cE_{v}\otimes_{\cO_K}\cO_{U_v}.$ 

Let $(L,\cE),(L',\cE')$ be two v.b.-type objects. Since the local
systems $\cE$ and $\cE'$ can be trivialized over $P_v$, we suppress them from the
notation and assume that we are dealing with trivial rank 1 local systems. We also assume for now that
the only element of $\cX(L,L';\varepsilon h)$ which lies inside $P_v$ is the node $p_v$ 
itself (this can always be achieved by a Hamiltonian isotopy), and the elements of $\cX(L,L';H)$ inside $P_v$ consist of the
generator $p_v$ together with infinitely many generators
$q_{e/v,0},q_{e/v,1},q_{e/v,2},\dots$ in degree 1 and $\dots,
p_{e/v,2},p_{e/v,1},p_{e/v,0},p_{e/v,-1},p_{e/v,-2},\dots$ in degree 0 in
$(0,3)\times S^1 \subset M_e$, for each half-edge $e/v$.
(We index the degree 0 generators so that $p_{e/v,k}$ lies near $\tau=1$ for $k\geq 0$, and
near $\tau=3$ for $k<0$).

The Floer differential maps $p_v$ to a linear combination of the three
degree 1 generators $q_{e/v,0}$ immediately adjacent to it along each of the
three edges, and each $p_{e/v,k}$ $(k\geq 0)$ to a multiple of the corresponding generator
$q_{e/v,k}$ (these do not involve propagation). It also maps
$p_{e/v,-k}$ $(k\geq 1)$ to 
$$\sum_{e'/v,\ e'\neq e}\, \sum_{\ell\geq 0} \,C^{v;e,e'}_{k,\ell}\,
T^{S_{e/v}(p_{e/v,-k})-S_{e'/v}(q_{e'/v,\ell})}\,q_{e'/v,\ell},$$
where the terms in the sum correspond to strips which propagate from $M_e$ to $M_{e'}$ through $p_v$ with
input degree $k$ and output degree $\ell$; here $C^{v;e,e'}_{k,\ell}$ is as
in Definition \ref{defi:propcoeff}, and $S_{e/v}(p_{e/v,-k})$ and
$-S_{e'/v}(q_{e'/v,\ell})$ are the areas of the two components. (As a matter
of convention we denote by $S_{e/v}(x)$ the signed area of a disc
connecting a Floer generator $x$ inside $M_e$ to $p_v$, so the signed
area of a disc from $p_v$ to $x$ is $-S_{e/v}(x)$.)

It follows from this that the Floer differential is surjective (even after
completion, as the construction of the Hamiltonians $H_n$ and $H$ ensures a
uniform bound on the areas of the trajectories connecting $p_{e/v,k}$ to 
$q_{e/v,k}$ independently of $k$), and the cohomology is concentrated in
degree zero, with generators
\begin{eqnarray}\label{eq:tildepv}
\tilde{p}_v&=&p_v\,+\,\sum_{e/v}\,
T^{-S_{e/v}(p_{e/v,0})}\,p_{e/v,0}\quad \text{and}\\
\label{eq:tildepev}
\tilde{p}_{e/v,-k}&=&p_{e/v,-k}\,+\sum_{e'/v,\ e'\neq e}\, \sum_{\ell\geq 0} \,C^{v;e,e'}_{k,\ell}\,
T^{S_{e/v}(p_{e/v,-k})-S_{e'/v}(p_{e'/v,\ell})}\,p_{e'/v,\ell},
\end{eqnarray}
where the exponents of $T$ correspond to the areas of trajectories between 
$p_v$ and the respective generators. The situation is similar for general
v.b.-type objects, after a suitable relabelling of the generators.

There is in fact a simple geometric model, which we denote by $\cF(P_v)$, 
where the Floer differential vanishes and morphism spaces are the
cohomologies of the morphism spaces in $\cF(P_v;H)$. Namely, we 
consider a Hamiltonian which behaves like $\varepsilon h$ in the interior of
$P_v$ and like $H$ near the boundary of $P_v$ (at $\tau=3$ in each of the three
components of $M$ which meet at $p_v$). It is still the case that
the generators outside of $P_v$ form an $A_\infty$-ideal, by the same argument
as in Section \ref{ss:FMH}; and the generators inside $P_v$ now consist of
$p_v$ and the $p_{e/v,-k}$ for all $e/v$ and $k\geq 1$, all in degree zero. 
Via either HPL or continuation maps, it can be seen that $\cF(P_v)$ is quasi-equivalent to
$\cF(P_v;H)$, with the linear term of the quasi-equivalence mapping $p_v$ to
$\tilde{p}_v$ and $p_{e/v,-k}$ to $\tilde{p}_{e/v,-k}$.

It is now apparent how to define the functor from $\cF(P_v)$ (resp.\
$\cF(P_v;H)$) to $\Perf(U_v)$ on morphism spaces (resp.\ closed degree 0
morphisms) between v.b.-type objects: we map
$p_v$ (resp.\ $\tilde{p}_v$) to the constant function 1 on $U_v$, and
$p_{e/v,-k}$ (resp.\ $\tilde{p}_{e/v,-k}$) to 
$$T^{S_{e/v}(p_{e/v,-k})}\, t_{e/v}^{-k}.$$
To prove that this is indeed a functor, we verify
that Floer products in $\cF(P_v)$ correspond to products of functions on
$U_v$: denoting by $\hat{p}_{e/v,-k}=T^{-S_{e/v}(p_{e/v,-k})}\,p_{e/v,-k}$
the Floer generators rescaled by appropriate area weights, and considering the
various types of propagating holomorphic discs in $P_v$ with inputs at two
given generators $p_{e_1/v,-k_1}$ and $p_{e_2/v,-k_2}$ lying on different
components ($e_1\neq e_2$), we have
$$\mu^2(\hat{p}_{e_1/v,-k_1},\hat{p}_{e_2/v,-k_2})=K^{v;e_1,e_2}_{k_1,k_2} p_v+
\sum_{b=0}^{k_1-1}C^{v;e_2,e_1}_{k_2,b} \hat{p}_{e_1/v,b-k_1} +
\sum_{a=0}^{k_2-1}C^{v;e_1,e_2}_{k_1,a} \hat{p}_{e_2/v,a-k_2},$$
which matches exactly the product formula in equation
\eqref{eq:partialproduct}. Meanwhile, for generators lying on the same
component the result is immediate since
$\mu^2(\hat{p}_{e/v,-k_1},\hat{p}_{e/v,-k_2})=\hat{p}_{e/v,-k_1-k_2}$.

(Defining the functor explicitly on the remaining part of the
morphism spaces in $\cF(P_v;H)$, if one wishes to do so, is best
accomplished by using homological perturbation theory to lift the strict
functor $\cF(P_v)\to \Perf(U_v)$ to an $A_\infty$-functor $\cF(P_v;H)\to
\Perf(U_v)$; however we will not need an explicit formula.)

Finally, verifying that the functor is full and faithful involves a
comparison of completions. Namely, morphisms in $\cF(P_v)$ are
infinite linear combinations of Floer generators such that the Novikov valuations of the
coefficients go to $+\infty$, whereas functions on the
open affinoid domain $U_v$ are linear combinations of the basis functions
$1$ and $t_{e/v}^{-k}$ for all $e/v$ and $k\geq 1$, such
that convergence holds whenever $|t_{e/v}|\geq |q_e|^{3/4}$ (i.e.,
$\val(t_{e/v})\leq \frac34 A_e$).
The fact that these two completions agree under our functor mapping
$p_{e/v,-k}$ to $T^{S_{e/v}(p_{e/v,-k})}\,t_{e/v}^{-k}$ follows directly
from the geometric fact that the area $S_{e/v}(p_{e/v,-k})$
of the degree $k$ disc connecting the generator $p_{e/v,-k}$ near
$\tau=3$ in $M_e$ to $p_v$ differs from $\frac34 k A_e$ by a bounded amount.
\medskip

The functor $\cF(N_e;H)\to \Perf(U_e)$ is constructed similarly, with all
v.b.-type objects mapped to free sheaves over $U_e$ and Floer generators
mapped to suitable multiples of powers of the coordinate $t_{e/v}$ (or
equivalently $t_{e/v'}$ for the other vertex). Viewing $N_e$ as a subset of
$P_v$, and considering a pair of v.b.-type Lagrangians which do not intersect in $P_v$
outside of the node $p_v$ as previously, their morphism space in
$\cF(N_e;H)$ is the completion of the span of the infinite sequence of generators
$p_{e/v,k}$, $k\in\Z$, all in degree zero, and we map each $p_{e/v,k}$ to
$T^{S_{e/v}(p_{e/v,k})}\,t_{e/v}^k$. The fact that the completions agree
under this functor follows again from the observation that
$S_{e/v}(p_{e/v,k})$ is close to $\frac34 |k| A_e$ for $k\ll 0$ (the generators
which lie near $\tau=3$), and to $-\frac14 k A_e$ for $k\gg 0$ (the generators near $\tau=1$).

By definition the restriction functor $\cF(P_v;H)\to \cF(N_e;H)$ maps
morphism spaces to each other simply by quotienting by all
the generators which lie outside of $N_e$; we denote this quotient map by $Q$. Composing with the
quasi-equivalence from $\cF(P_v)$ into $\cF(P_v;H)$ provided by HPL (or
continuation), we obtain a restriction functor $\cF(P_v)\to \cF(N_e;H)$.
In light of \eqref{eq:tildepv}--\eqref{eq:tildepev}, this maps
$p_v$ to $Q(\tilde{p}_v)=T^{-S_{e/v}(p_{e/v},0)} p_{e/v,0}$, 
$p_{e/v,-k}$ to $Q(\tilde{p}_{e/v,-k})=p_{e/v,-k}$ itself, and for $e'\neq e$, $p_{e'/v,-k}$ to 
$$Q(\tilde{p}_{e'/v,-k})=\sum_{\ell\geq 0} \,C^{v;e',e}_{k,\ell}\,
T^{S_{e'/v}(p_{e'/v,-k})-S_{e/v}(p_{e/v,\ell})}\,p_{e/v,\ell}.$$
These formulas are easily checked to agree with the restriction from $\Perf(U_v)$ to
$\Perf(U_e)$, using the fact that $t_{e'/v}^{-k}=\sum_{\ell=0}^\infty C_{k,\ell}^{v;e',e}\, t_{e/v}^\ell$.

\subsection{Theta functions}\label{ss:theta}
Now we show how the ingredients of the construction assemble to give a concrete
description of the mirror functor $\cF(M)\to \Perf(X_K)$, in the special case when the Floer complex $CF^*((L,\cE),(L',\cE'))$ is concentrated in degree zero,
by providing an explicit map 
\begin{equation}\label{eq:theta}\Phi_{L,L'}:CF^0((L,\cE),(L',\cE'))\to \Hom(\Phi(L,\cE),\Phi(L',\cE')).\end{equation}

We consider two objects $(L,\cE),$ $(L',\cE'),$ and an intersection point $x\in L\cap L'$ of degree zero, such that $x\in P_v.$ 
We explain how to associate to it a map 
\begin{equation}\label{eq:map_for_intersection_points}\Phi_{L,L',x}:\Hom_{\cO_K}(\cE_x,\cE'_x)\to \Hom_{\cO_{U_v}}(\Phi(L,\cE)_{\mid U_v},\Phi(L',\cE')_{\mid U_v}).\end{equation}

Take the half-edge $e/v$ in the graph $G$ such that $x\in M_e.$ 
We denote by $r_{e/v}(x)\in \Z$ the rotation number of
$\phi_H^1(L)$ relative to $L'$ in the negative direction along the path from
$p_v$ to $x$, and by $S_{e/v}(x)$ the signed area of a disc connecting $x$ to $p_v$ inside $M_e$, or equivalently, the region bounded by $L$ and $L'$ on
the universal cover of $(0,\tau(x))\times S^1\subset M_e$ (taking the lifts which approach each other as $\tau\to \tau(x)$).
In the case when $x=p_v,$ we have $r_{e/v}(p_v)=0$ and $S_{e/v}(p_v)=0$. 
If $e$ connects $v$ and $v',$ then
\begin{equation}\label{eq:relations_on_r_and_S}r_{e/v}(x)+r_{e/v'}(x)=r_{e}(L,L'),\quad S_{e/v}(x)-S_{e/v'}(x)+r_{e/v'}(x)A_e=S_{e/v}(L,L').\end{equation}

To each such $x$ we associate a monomial $T^{S_{e/v}(x)}t_{e/v}^{-r_{e/v}(x)},$ considered as a function on $U_v.$ Now, we define the map \eqref{eq:map_for_intersection_points} by the formula 
\begin{equation}\label{eq:monomial_for_a_point}\Phi_{L,L',x}(\varphi)=(R_{\cE',x,v}\varphi R_{\cE,v,x})\otimes T^{S_{e/v}(x)}t_{e/v}^{-r_{e/v}(x)}\in \Hom_{\cO_{U_v}}(\Phi(L,\cE)_{\mid U_v},\Phi(L',\cE')_{\mid U_v}),\end{equation} where $\varphi\in\Hom(\cE_x,\cE'_x),$ and $R_{\cE,v,x},$ $R_{\cE',x,v}$ denote the monodromies. Moreover, the morphisms $(R_{\cE',x,v}\varphi R_{\cE,v,x})\otimes T^{S_{e/v}(x)}t_{e/v}^{-r_{e/v}(x)}$ and $(R_{\cE',x,v'}\varphi R_{\cE,v',x})\otimes T^{S_{e/v'}(x)}t_{e/v'}^{-r_{e/v'}(x)}$
agree on $U_e,$ which follows from the gluing data \eqref{eq:gluingdata} for $\Phi(L,\cE),$ $\Phi(L',\cE')$ and
the observation that, using \eqref{eq:relations_on_r_and_S},
\begin{multline*} \left(T^{-S_{e/v}(L,L')}\, t_{e/v}^{r_e(L,L')}\right)\left(T^{S_{e/v}(x)}\,t_{e/v}^{-r_{e/v}(x)}\right) =
T^{S_{e/v}(x)-S_{e/v}(L,L')}\, t_{e/v}^{r_{e/v'}(x)}\\=T^{S_{e/v}(x)-S_{e/v}(L,L')+r_{e/v'}(x)A_e}
\,t_{e/v'}^{-r_{e/v'}(x)}=T^{S_{e/v'}(x)} \,t_{e/v'}^{-r_{e/v'}(x)}.\end{multline*}

Now we introduce some notation. Let us take any reduced path in $G,$ written as $\gamma=(v_0,e_1,\dots,e_n,v_n).$ It gives a map $g_{\gamma}:Y_{v_0}\to Y_{v_n},$ given by $$g_{\gamma}=g_{e_n/v_{n-1}}\circ\dots\circ g_{e_1/v_0}.$$ We denote by $u_{L,L',\gamma}:\Hom(\cE_{v_0},\cE'_{v_0}\otimes \cO_{Y_{v_0}-F_{v_0}})\to g_{\gamma}^*(\Hom(\cE_{v_n},\cE'_{v_n})\otimes \cO_{Y_{v_n}-F_{v_n}})$ the gluing morphism.

The morphism \eqref{eq:theta} is given by ``averaging'' the morphisms \eqref{eq:monomial_for_a_point}. Namely, for a half-edge $e_0/v_0,$ a point $x\in L\cap L'\cap (\mathrm{int}(M_{e_0})\sqcup\{v_0\}),$ and a morphism $\varphi:\cE_x\to\cE'_x,$ for any vertex $v\in V(G)$ we put
\begin{equation}\label{eq:averaging}\Phi_{L,L'}(\varphi)_{\mid U_v}=\sum\limits_{\gamma:v_0\to v}g_{\gamma*}(u_{L,L',\gamma}((R_{\cE',x,v_0}\varphi R_{\cE,v_0,x})\otimes T^{S_{e_0/v_0}(x)}t_{e_0/v_0}^{-r_{e_0/v_0}(x)})).\end{equation} This sum converges because of our assumption on the Floer complex $CF^{*}((L,\cE),(L',\cE'))$ to be concentrated in degree zero. The restrictions of $\Phi_{L,L'}(\varphi)$ to different $U_v$ agree on the intersections, so we get a well-defined morphism of vector bundles $\Phi(L,\cE)\to \Phi(L',\cE').$

\begin{example} Consider a graph $G$ with a single vertex
$v$, a single internal edge $e$ connecting $v$ to itself, and an external
edge $e'$, so that $M$ is the union of $M_e = \PP^1 / (0\sim
\infty)$ and $M_{e'}=\C$, glued together at the node of $M_e$.
Choose the combinatorial data for the vertex $v$ so that the coordinates
$t^{\pm 1}$  associated to the two ends of the edge $e$ are inverses of each other
and take the value 1 at the third puncture.
The mirror is then the punctured elliptic curve $X_K=(K^*-q_e^\Z)/q_e^\Z$, where
$q_e=T^{A_e}$. Let $L_0,L_1$ be v.b.-type objects such that $L_1$ rotates
once clockwise around $L_0$ as it travels from the node of $M_e$ back
to itself, i.e.\ $r_e(L_0,L_1)=1$, and $S_{e/v}(L_0,L_1)=\frac12 A_e$; we
equip both with the trivial local system. Following the construction in
\S \ref{ss:assignbundles}, our mirror functor associates to
$L_0$ the structure sheaf $\O_X$, and to $L_1$ a certain line bundle
obtained as a quotient of the trivial line bundle over $(K^*-q_e^\Z)$. 
For a path $\gamma_n$ from $v$ to itself which goes $n$ times around the edge $e$, the
coordinate change $g_{\gamma_n}$ is given by $t\mapsto q_e^{-n} t$, and the gluing
map $u_{L_1,\gamma_n}$ for the fibers of the line bundle $\Phi(L_1),$ or
equivalently $u_{L_0,L_1,\gamma_n}$ for $\Hom(\Phi(L_0),\Phi(L_1)),$ is given by 
multiplication by $q_e^{-n^2/2} t^{n}$. 
Let $x=p_v$ be the intersection of $L_0$ and $L_1$ which lies at the node, and
let $\varphi\in CF^0(L_0,L_1)$ be the morphism defined by this intersection
point and the identity map between the trivial local systems. Then
\eqref{eq:averaging} becomes
$$\Phi_{L_0,L_1}(\varphi)=\sum_{n\in\Z} g_{\gamma_n *}(u_{L_0,L_1,\gamma_n}(\mathrm{id}))=
\id\otimes \sum_{n\in\Z} 
q_e^{-n^2/2} \bigl(q_e^{n}t\bigr)^{n}=\id\otimes \sum_{n\in\Z} q_e^{n^2/2}\, t^n$$
(as a morphism from $\O_X$ to the line bundle $\Phi(L_1)$, expressed in the trivialization of 
$\Phi(L_1)$ on the $\Z$-cover of $X_K$ by $K^*-q_e^\Z$); 
this recovers the classical theta function for the elliptic curve $K^*/q_e^\Z$.
\end{example}

Now we explain how HPL provides the averaging in \eqref{eq:averaging}. We need to compute the map $$\Phi_{L,L'}\circ\tilde{i}^1:CF((L,\cE),(L',\cE'))\to \Hom(\Phi(L,\cE),\Phi(\cL',\cE')),$$ where $\tilde{i}^1=(\id+h\delta^1)^{-1}i,$ and $i, h$ and $\delta$ are as in Section \ref{ss:HPLFloer}. 
Take some $x\in \cX(L,L';\varepsilon h),$ $\varphi\in\Hom(\cE_x,\cE'_x).$ We have $\tilde{i}^1(\varphi)=\sum\limits_{n=0}^{\infty}(-h\delta^1)^n i(\varphi).$  

Now, the map $h\delta^1:CF^0((L,\cE),(L',\cE');H)\to CF^0((L,\cE),(L',\cE');H)$ is described explicitly as follows. The formula \eqref{eq:monomial_for_a_point} provides an identification $$CF^0((L,\cE),(L',\cE');H)\cong \bigoplus\limits_{v\in V(G)}\Hom_{\cO_K}(\cE_v,\cE'_v)\oplus \bigoplus\limits_{e\in E(G)}\Hom(\Phi(L,\cE)_{\mid U_e},\Phi(L',\cE')_{\mid U_e}).$$ Under this identification, we have
$$h\delta^1(\varphi)=\sum\limits_{e/v}(\varphi\otimes\cO_{U_v})_{\mid U_e}\quad\text{ for }\varphi\in \Hom_{\cO_K}(\cE_v,\cE'_v),$$
Further, for each edge $e:v\to v',$ for $\varphi\in\Hom_{\cO_K}(\cE_v,\cE'_v)$ and $n\leq r_e(L,L'),$ the propagation rule implies the following: $$h\delta^1((\varphi\cdot t_{e/v}^{-n})_{\mid U_e})=\sum\limits_{e'/v,e'\ne e}(\varphi\cdot t_{e/v}^{-n})_{\mid U_{e'}}\quad\text{ for }n\geq r_e(L,L'),$$
and
$$h\delta^1((\varphi\cdot t_{e/v}^{-n})_{\mid U_e})=\sum\limits_{e'/v,e'\ne e}(\varphi\cdot t_{e/v}^{-n})_{\mid U_{e'}}+\sum\limits_{e''/v',e''\ne e} (R_{\cE',e/v}\varphi R_{\cE,e/v}^{-1})\cdot (T^{-nS_e(L,L')}t_{e/v'}^{n-r_e(L,L')})_{\mid U_{e''}}$$ for $0< n< r_e(L,L').$ It follows that the map $\Phi_{L,L'}\circ\tilde{i}^1$  gives exactly the averaging  \eqref{eq:averaging}. It is important that $CF^{\bullet}(L,L')$ is concentrated in degree zero, hence all the rotation numbers $r_e(L,L')$ are strictly positive and we don't
``lose'' any monomials while propagating.

\section{Canonical map} \label{s:canonicalmap}

Recall that for a smooth projective curve $C$ over a field $\mk,$ of genus $g\geq 2,$ we have the canonical map $\can:C\to \PP(H^0(C,\omega_C)^*)\cong \PP(H^1(C,\cO_C)).$ On $\mk$-rational points it can be described as
$$p\mapsto \im(\Ext^1(\cO_p,\cO_C)\otimes\Ext^0(\cO_C,\cO_p)\to \Ext^1(\cO_C,\cO_C)=H^1(C,\cO_C)).$$ This map is a closed embedding unless $C$ is hyperelliptic in which case it is $2:1$ onto its image.

Note that even when $C$ is reduced singular of arithmetic genus $g\geq 2,$ we still have a map $\can:C^{sm}\to \PP(H^1(C,\cO_C)).$ Moreover, for any Zariski open (resp. analytic open) subset $U\subset C^{sm}$ and a regular (resp. analytic) vector field $\theta\in H^0(U,T_U),$ we have a regular (resp. analytic) map $\can_{\theta}:U\to H^1(C,\cO_C).$

We will compute this map in our situation for a general trivalent graph (say, without loops, although they can be allowed), both on the A-side and the B-side, and we will see that they match.

\subsection{Canonical map:\ analytic setup} \ Here for simplicity we choose some
non-Archimedean normed field $K,$ and take the extension of scalars $X_K$ from $\Z[[\{q_{e}\}]]$ to $K,$ where $q_{e}$ are sent to some elements of $\m_K.$ Also, take the Schottky group $\Gamma=\Gamma_{e_0/v_0}.$

Then in the framework of rigid analytic geometry $X_K$ is identified with a quotient $(\PP^1_K-F)/\Gamma,$ where $F$ is the set of limit points of the group $\Gamma$ (equivalently, $F$ is the closure of the set of fixed points of non-identity elements of $\Gamma$). Now take a rational function $\phi$ on $\PP^1_K,$ which is regular at each point of $F.$ Then the collection of principal parts of $\phi$ at its poles defines a class $[\phi]\in H^1(X_K,\cO_K).$

Let us compute this class. Note that $$H^1(X_K,\cO_{X_K})\cong H^1(\Gamma,K)=\Hom(\Gamma/[\Gamma,\Gamma],K).$$ 
Now let us choose some point $t_0\in \PP^1_K-F,$ such that $\phi$ is regular at each point of $\Gamma t_0.$ Then we have a well-defined analytic function $$f_{\phi}(t):=\sum\limits_{g\in\Gamma} (\phi(gt)-\phi(gt_0)),$$
which is $\Gamma$-invariant up to adding a constant.
The associated class $[\phi]\in H^1(X_K,\cO_{X_K})=H^1(\Gamma,K)$ is given by  the cocycle \begin{equation}c_{\phi}(\gamma)=f_{\phi}(t)-f_{\phi}(\gamma t)=\sum\limits_{g\in\Gamma}(\phi(g\gamma(t_0))-\phi(g(t_0))),\quad \gamma\in\Gamma.\end{equation}

This cocycle of course does not depend on the choice of $t_0.$ Moreover, if $\gamma\ne 1,$ and $y_{\gamma}^0,y_{\gamma}^{\infty}\in \PP^1_K$ are the fixed points of $\gamma,$ with $y_{\gamma}^0$ being the attractor, then we have
\begin{equation}c_{\phi}(\gamma)=\sum\limits_{\bar{g}\in\Gamma/\gamma^{\Z}}(\phi(g(y_{\gamma}^0))-\phi(g(y_{\gamma}^{\infty}))).\end{equation}

Now, if we have an analytic open subset $U\subset \PP^1_K-F,$ such that $U\cap g(U)=\emptyset$ for all $g\in\Gamma\setminus\{1\},$ then we have $U\cong \pr(U)\subset X_K,$ and choosing the vector field $t\frac{\partial}{\partial t}$ on $U,$ we get the lifted canonical map $\can_{t\frac{\partial}{\partial t}}:U\to H^1(\Gamma,K).$ By the above discussion, this map is 
given by 
\begin{equation}\label{eq:can_map_analytic}\can_{t\frac{\partial}{\partial t}}(s)=c_{\frac{s}{t-s}}\in H^1(\Gamma,K),\quad c_{\frac{s}{t-s}}(\gamma)=\sum\limits_{g\in\Gamma}(\frac{s}{g\gamma(t_0)-s}-\frac{s}{g(t_0)-s}).\end{equation}

We will see how this $1$-cocycle arises both in the formal scheme framework and in the Fukaya framework.



\subsection{Canonical map: formal scheme} 

Here by $\mX$ we denote either the formal scheme over $\Z[[\{q_{e}\}]]$ introduced above, or its extension of scalars to some (nicely behaved) topological ring $R$ (where $q_{e}$ are sent to some topologically nilpotent elements). We also fix some $e_0/v_0$ and the corresponding Schottky group $\Gamma=\Gamma_{e_0/v_0}.$

Recall the open subsets $\cU_{e},\cW_{v}\subset\mX.$ Note that each intersection $\cU_{e}\cap \cU_{e'}$ (for $e\ne e'$) is either empty, or of the form $\cW_{v},$ or of the form $\cW_{v}\sqcup\cW_{v'}.$ Thus, given a coherent sheaf $\cF$ on $\mX$ we can (quasi-isomorphically) modify the \v Cech complex of $\cF$ for the covering $\{\cU_{e}\},$ and take the following complex: $$\cK(\cF):=\{\bigoplus_{e\in E}\Gamma(\cU_{e},\cF)\xto{d}\bigoplus_{v\in V} \Gamma(\cW_{v},\cF)\tens{\Z}V_{v}\},$$
where $$V_{v}=(\Z\cdot e_{e_1/v}\oplus \Z\cdot e_{e_2/v}\oplus \Z\cdot e_{e_3/v})/\Z\cdot (e_{e_1/v}+e_{e_2/v}+e_{e_3/v}),$$ and
$$d(\{f_{e}\})_{v}=\sum\limits_{e'/v}f_{e'}e_{e'/v}.$$

It is not hard to check directly that we have a quasi-isomorphic subcomplex $\cK_{const}(\cO)\subset \cK(\cO),$ formed by constant local sections (on $\cU_{e}$ and $\cW_{v}$). We can write down explicitly the identification $H^1(\cK_{const}(\cO))\cong H^1(\Gamma,R).$ Namely, for $e/v,e'/v,$ denote by $\xi^{e,e'}_{v}:V_{v}\to\Z$ the functional $e_{e/v}^*-e_{e'/v}^*.$ Then an element $\{a_{v}\}\in \cK^1_{const}(\cO)$ defines a cocycle
\begin{equation}\label{eq:from_Cech_to_group_cohom}c_a\in H^1(\Gamma,R),\quad c_a(\gamma_P)=\xi_{v_1}^{e_1,e_2}(a_{v_1})+\dots+\xi_{v_n}^{e_n,e_1}(a_{v_n}),\end{equation}
for $P=(v_0,e_1,v_1,\dots,e_n,v_n=v_0).$  

Now let us take a rational function $\phi$ on $\PP^1_R$ which is regular at $0,1,\infty.$ By this we mean $\phi(t)=\frac{h_1(t)}{h_2(t)},$ where $h_2(t)$ is monic, $\deg(h_1)\leq\deg(h_2),$ and $h_2(0),h_2(1)\in R$ are invertible. Then we get a coherent sheaf $\cF_{h_2}\supset \cO,$ such that $\supp(\cF_{h_2}/\cO)\subset \cW_{v_0}$ and $\cF(\cW_{v_0})=\frac{1}{h_2(t_{e_0/v_0})}\cO(\cW_{v_0}).$  Then the ``principal parts'' of $\phi(t_{e_0/v_0})$ give a well-defined element of $H^1(\mX,\cO_{\mX}).$ Let us compute a representative of this class in $\cK^1_{const}(\cO).$

We first take the sections $f_{e}\in\Gamma(\cU_{e},\cF_{h_2}),$ given by
$$f_{e}=\sum\limits_{\substack{P=(v_0,e_1,\dots,e_n,v_n);\\
e/v_n,e\ne e_n}}(\phi(\gamma_P^{e_0,e}(T_{e/v_n}))-\phi(\gamma_P^{e_0,e}(0)))$$
(it can be checked directly that $f_{e}$ are well-defined), and then notice that $d(\{f_{e}\})\in \cK^1(\cF)$ is contained in $\cK^1_{const}(\cO)\subset \cK^1(\cO)\subset \cK^1(\cF).$

Thus, $d(f_{e})$ is our desired constant representative, which then gives a class in $H^1(\Gamma,R)$ by the formula \eqref{eq:from_Cech_to_group_cohom}. By straightforward combinatorial considerations one checks that the result actually agrees with \eqref{eq:can_map_analytic}.  Now taking $s\in R$ such that $s(1-s)$ is invertible, we see that the class 
$\can_{t_{e_0/v_0}\,\partial/\partial t_{e_0/v_0}}\in H^1(\mX,\cO_{\mX})$ is given again by the formula \eqref{eq:can_map_analytic}.

\begin{remark}To make sense of canonical map for $|s|<1$ we need to invert $q_{e_0}$ as described in Remark \ref{remark:inverting_q}; the computation works in exactly the same way.\end{remark}

\subsection{Canonical map: Fukaya category}

Here we take the singular symplectic manifold $M$ as above; recall that the symplectic areas are denoted by $A_{e},$ $e\in E.$ Again, we fix $e_0/v_0,$ and also take $v_0'\ne v_0,$ $e_0/v_0'.$ 

We take $L_0$ to be a v.b.-type Lagrangian with trivial rank one local system, 
corresponding to $\cO_{\mX}$ under mirror symmetry, and orient its $M_{e_0}$ component from $v_0'$ to $v_0$.

The Floer complex $\Hom(L_0,L_0)$ is just the complex computing the cohomology of the graph $G$ (with vertices being $v$ and edges being $e$). We denote by $p_{v},$ resp.\ $z_{e}$ its generators of degree 0, resp.~1,
corresponding to the points of $\cX(L_0,L_0)=L_0^+\cap L_0$. (Recall that $L_0^+=\phi^1_{\varepsilon h}(L_0)$ is a slight pushoff of $L_0$ in 
the counterclockwise direction near each vertex $p_v$, and intersects $L_0$ at the vertices and also once inside each component $M_e$).

Now, let $L_1$ be a point-type object, i.e.\ a circle on the $M_e$ component, placed between the points $z_{e_0}$ and $p_{v_0},$ oriented in such a way that $\Hom(L_0,L_1)$ is in degree zero, hence $\Hom(L_1,L_0)$ is in degree 1 (and we take the trivial local system on $L_1$ for simplicity). We put $y_1:=L_0^+\cap L_1,$ $y_2:=L_1\cap L_0.$ So, $y_1\in\Hom(L_0,L_1)$ and $y_2\in\Hom(L_1,L_0).$ We are interested in 
$$\mu^2(y_2,y_1)=\sum\limits_{e}a_{e}z_{e}\in\Hom^1(L_0,L_0).$$ Let us denote by $B$ the area of the half-sphere with the boundary $L_1,$ containing the node $v_0.$

Now we determine the constants $a_{e}$. First, for $e=e_0$, $L_0$, $L_1$ and $L_0^+$ bound a small thin triangle inside $M_{e_0}$ with vertices $y_1,y_2,z_{e_0}$; 
the corresponding perturbed disc has area zero since two of its edges lie on $L_0$, so its area weight is 1.
All the other holomorphic strips will propagate through the nodes, and to count them we introduce some notation.

Namely, for $e/v,e'/v,$ we denote by $C^{v;e,e'}_{k,l}\in\Z$ (where $k,l\geq 0$)  the constants such that $$\frac{1}{g_{v}^{e,e'}(t)^k}=\sum\limits_{l\geq 0}C^{v;e,e'}_{k,l}t^{l}.$$
For $k\geq 1$ these are exactly the propagation coefficients introduced in Section \ref{s:Amodel}; the constants $C^{v;e,e'}_{0,l}=\delta_{0,l}$ do not participate in the propagation rules but it is convenient to include them. We will also adopt the following notation: $$\delta_{v}^{e,e'}=\begin{cases}1 & \text{ for }e/v,\,e'/v,\,e\ne e';\\0 & \text{otherwise}\end{cases}.$$

Now, the perturbed propagating holomorphic strips contributing to $a_{e}$ (other than the already mentioned triangle) are divided into two types:

(I) the ones which first come to $v_0$ with some degree $k>0,$ then propagate along some path (in our graph), and finally arrive to the component $e$ with degree $0;$

(II) The same with $v_0'$ instead of $v_0.$

The contribution of the strips of type (I) is the following sum:
\begin{equation}\label{eq:strips_through_lambda_0}a_{e,v_0}=\sum\limits_{\substack{P=(v_0,e_1,v_1,\dots,e_n,v_n);\\
e/v_n,e\ne e_n,e_1\ne e_0}}\left(\sum\limits_{k,k_1,\dots,k_n>0}C^{v_0;e_0,e_1}_{k,k_1}C^{v_1;e_1,e_2}_{k_1,k_2}\dots C^{v_n;e_n,e}_{k_n,0}\,T^{\,kB+\smash{\sum\limits_{l=1}^n} k_i A_{e_i}}\right).\end{equation}
Now let us notice the following identity: for a reduced path $P$ as in \eqref{eq:strips_through_lambda_0}, and for $k>0$ we have
\begin{equation}\label{eq:identity_for_Schottky_groupoid}
\sum\limits_{k_1,\dots,k_n\geq 0}C^{v_0;e_0,e_1}_{k,k_1}C^{v_1;e_1,e_2}_{k_1,k_2}\dots C^{v_n;e_n,e}_{k_n,0}q_{e_1}^{k_1}\dots q_{e_n}^{k_n}=\gamma_P^{e_0,e}(0)^{-k}
\end{equation}
(note the non-strict inequalities for $k_i$).
Now, if $n>0,$ then let us denote by $P'$ the path $(v_0,e_1,\dots,e_{n-1},v_{n-1})$ (removing the last edge from $P$). Then from \eqref{eq:identity_for_Schottky_groupoid} we get
\begin{equation}\label{eq:identity_with_positive_k_i}
\sum\limits_{k_1,\dots,k_n>0}C^{v_0;e_0,e_1}_{k,k_1}C^{v_1;e_1,e_2}_{k_1,k_2}\dots C^{v_n;e_n,e}_{k_n,0}q_{e_1}^{k_1}\dots q_{e_n}^{k_n}=\begin{cases}\gamma_P^{e_0,e}(0)^{-k}-\gamma_{P'}^{e_0,e_{n}}(0)^{-k} & \text{for }n>0;\\
g_{v_0}^{e_0,e}(0)^{-k} & \text{for }n=0.\end{cases}.
\end{equation} 
Combining \eqref{eq:identity_with_positive_k_i} with \eqref{eq:strips_through_lambda_0}, and identifying $q_{e}$
 with $T^{A_{e}},$ we get
\begin{equation}\label{eq:strips_through_lambda_0_final}a_{e,v_0}=\sum\limits_{\substack{P=(v_0,e_1,v_1,\dots,e_n,v_n);\\
n>0,e/v_n,e\ne e_n,e_1\ne e_0}}\pm\left(\frac{T^B}{\gamma_P^{e_0,e}(0)-T^B}-\frac{T^B}{\gamma_{P'}^{e_0,e_n}(0)-T^B}\right)\pm\delta_{v_0}^{e_0,e}\frac{T^B}{g_{v_0}^{e_0,e}(0)-T^B}.\end{equation}
Now, the strips of type II are completely analogous. Taking into account the identity $$\frac{(\frac{q}{s})}{t-\frac{q}{s}}=-\left(\frac{s}{\frac{q}{t}-s}+1\right),$$ we get that the contribution of strips of type II equals
\begin{equation}\label{eq:strips_through_lambda_0'}a_{e,v_0'}=\sum\limits_{\substack{P=(v_0,e_0,v_0',e_1,v_1,\dots,e_n,v_{n+1});\\
n\geq 0,e/v_{n+1},e\ne e_n}}\pm\left(\frac{T^B}{\gamma_P^{e_0,e}(0)-T^B}-\frac{T^B}{\gamma_{P'}^{e_0,e_n}(0)-T^B}\right).\end{equation}
So, combining \eqref{eq:strips_through_lambda_0_final}, \eqref{eq:strips_through_lambda_0'}, and taking into account the small triangle in $M_{e_0}$, we get
\begin{multline}a_{e}=\sum\limits_{\substack{P=(v_0,e_1,v_1,\dots,e_n,v_n);\\
n>0,e/v_n,e\ne e_n}}\pm\left(\frac{T^B}{\gamma_P^{e_0,e}(0)-T^B}-\frac{T^B}{\gamma_{P'}^{e_0,e_n}(0)-T^B}\right)\\[-1em]
\pm\delta_{v_0}^{e_0,e}\frac{T^B}{g_{v_0}^{e_0,e}(0)-T^B}
\pm\delta_{e_0,e}.\end{multline}

This completes the calculation of $\mu^2(y_2,y_1)\in\Hom^1(L_0,L_0).$ To get the value of the corresponding class $c_{L_1}\in H^1(\Gamma,R)$ on an element $\gamma\in\Gamma,$ we simply need to sum up $\pm a_{e}$ along a path. The same combinatorics as in the previous subsection shows that
$$c_{L_1}=\can_{t\frac{\partial}{\partial t}}(T^B),$$
where the RHS is given by \eqref{eq:can_map_analytic}.
So, we see that the canonical map indeed allows one to identify the points in the annulus $\{1>|t_{e_0/v_0}|>|T^{A_{e_0}}|\}$ and the circles with $1$-dimensional local systems via $t_{e_0/v_0}=T^{B}\cdot(\text{monodromy}).$


\section{Epilogue: higher dimensions} \label{s:higherdim}

We expect that the constructions and results described in this paper
for curves and their mirrors admit higher dimensional generalizations;
the details are still tentative as of this writing, so much so that we
do not even formulate a precise conjecture.

On the B-side, we consider a rigid analytic space $X_K$ admitting a
{\em generalized pair-of-pants decomposition}, i.e.\ an open cover by
analytic subsets $U_v$ which are obtained from the $n$-dimensional
pair of pants $$\Pi_n=\bigl\{(z_0:z_1:\dots:z_{n+1})\in \PP^{n+1}\,\big|\,
{\textstyle \sum z_i=0},\ z_j\neq 0\ \forall j\bigr\}$$ 
(i.e., the complement of $n+2$ generic hyperplanes in $\PP^n$)
by imposing suitable inequalities on the valuations of the coordinates.
Such decompositions arise most commonly from maximal degenerations of
complex varieties to the tropical limit. The combinatorics of the
decomposition is then encoded by the {\em tropicalization} of $X_K$,
a polyhedral complex $\Sigma$ in which the vertices index the
pairs-of-pants $U_v$ in the decomposition of $X_K$, and the higher-dimensional
strata and their affine structures determine which subsets
$U_v$ overlap non-trivially and how the valuations of the coordinates 
transform under gluing maps. 

Each pair-of-pants $U_v$ has
$\binom{n+2}{n}$ {\em ends} $U_{\sigma/v}$, namely the subsets of $U_v$ 
where $n$ of the homogeneous coordinates
are smaller than the remaining two, and each of these ends corresponds
to one of the $\binom{n+2}{n}$ top-dimensional strata $\sigma\subset \Sigma$ which meet 
at $v$. Also, each top-dimensional stratum $\sigma$ of $\Sigma$ determines a subset
$U_\sigma\subset X_K$ along which the ends $U_{\sigma/v}$
of the pairs of pants $U_v$ for $v$ adjacent to $\sigma$ 
overlap. After choosing suitable
coordinates, $U_\sigma$ can be identified (in a valuation-preserving manner)
with the affinoid domain $\val^{-1}(\sigma)\subset (K^*)^n$ determined by the
affine structure on $\sigma$. This in turns yields $(K^*)^n$-valued
local coordinates $t_{\sigma/v}$ on each end $U_{\sigma/v}$ of each
pair of pants in the decomposition of $X_K$; by construction these coordinates
have the same valuations as ratios of homogeneous coordinates on the model
pair of pants $\Pi_n$ into which $U_v$ embeds, and they transform by monomial coordinate changes between the
different vertices $v$ adjacent to a given top-dimensional stratum $\sigma$.
(Note that a complete description of $X_K$ also involves gluings over strata of
all dimensions in $\Sigma$, which we omit from our discussion for
simplicity.)

On the A-side, we consider a stratified space $M$ formed by a union of 
toric K\"ahler manifolds $M_\sigma$ glued together along toric strata, with
local models along codimension $k$ strata given by
the product of $(\C^*)^{n-k}$ with the union $\amalg_k$ of all coordinate $k$-planes in
$\C^{k+2}$ (in particular there are $\binom{n+2}{n}$
top-dimensional strata meeting at each vertex).
For mirror symmetry purposes, the moment map images of the various
components of $M$ should match the strata of the polyhedral complex
$\Sigma=\mathrm{Trop}(X_K)$ and their affine structures.

One can then consider stratified Lagrangian submanifolds $L\subset M$
which, along each codimension $k$ stratum of $M$, are modelled on the product
of a smooth Lagrangian submanifold of $(\C^*)^{n-k}$ with the union of
the real positive loci in $\amalg_k$.  One may further require each
component $L_\sigma\subset M_\sigma$ of $L$ to be a section of the
moment map fibration $M_\sigma\to \sigma$, and equip $L$ with a unitary
rank 1 local system $\mathcal{E}$ (trivialized at the vertices). 
We expect that such
objects correspond to line bundles on $X_K$, constructed from 
trivial line bundles over each pair of pants $U_v$ by gluing them over
$U_\sigma$ via transition functions determined explicitly by the rotation 
numbers of $L$ with respect to the real positive locus $L_0$ inside each
stratum of $M$, the signed symplectic areas of triangular regions
bounded by $L$ and $L_0$, and the holonomy of the local system $\cE$,
as in Section \ref{ss:assignbundles}.

As in the case of curves, one expects the definition of morphism spaces in the
Fukaya category of $M$ to involve Hamiltonian perturbations whose flow 
rotates the asymptotic directions of the Lagrangians by a small positive
amount around each lower-dimensional stratum (as well as
wrapping near infinity when $M$ is non-compact). 
The structure maps of the Fukaya category should then
involve weighted counts of rigid configurations of (perturbed) holomorphic
discs which are allowed to propagate among the strata of $M$
through lower-dimensional strata. The simplest case, and the only one we
shall discuss here, involves a Floer trajectory
propagating from a top-dimensional stratum $M_{\sigma_{in}}$ to another
top-dimensional stratum $M_{\sigma_{out}}$ through a common vertex $v$
shared by the two strata. In this case, the local behavior of the
holomorphic disc near $v$ can be described by
associating to the incoming component 
a degree $k_{in}\in \Z_{>0}^n$, and to the outgoing component a
degree $k_{out}\in \Z_{\geq 0}^n$. We then expect that the propagation 
coefficient $C^{v;\sigma_{in},\sigma_{out}}_{k_{in},k_{out}}$ (the
local contribution to the multiplicity of the propagating Floer
trajectory) should be defined as
the coefficient of the monomial $t^{k_{out}}_{\sigma_{out}/v}$ in the
expansion of (the analytic continuation of) $t_{\sigma_{in}/v}^{-k_{in}}$ 
as a power series in terms of the coordinates $t_{\sigma_{out}/v}$ over
$U_{\sigma_{out}/v}$. 

Obviously, a lot of further work is needed to flesh out key 
details of this story and confirm that the proposed construction is 
sound and leads to a homological mirror symmetry statement; this work
is still in the early stages, but we
hope to have convinced the reader that the story developed in this paper
is likely to extend beyond the case of curves.


\begin{thebibliography}{AAEKO}
\bibitem[Ab]{AbGen}
   M.~Abouzaid, {\sl A geometric criterion for generating the Fukaya category}, 
   Publ.\ Math.\ IH\'ES {\bf 112} (2010), 191--240.
\bibitem[AA]{AA}
   M.~Abouzaid, D.~Auroux, 
   {\sl Homological mirror symmetry for hypersurfaces in $(\C^*)^n$}, 
   to appear in Geom.\ Topol., arXiv:2111.06543.
\bibitem[AAEKO]{AAEKO}
   M.~Abouzaid, D.~Auroux, A.~I.~Efimov, L.~Katzarkov, D.~Orlov,
   {\sl Homological mirror symmetry for punctured spheres}, 
   J.\ Amer.\ Math.\ Soc. {\bf 26} (2013), 1051--1083. 
\bibitem[AAK]{AAK}
   M.\ Abouzaid, D.\ Auroux, L.\ Katzarkov,
   {\sl Lagrangian fibrations on blowups of toric varieties and mirror symmetry
   for hypersurfaces},
   Publ.\ Math.\ IH\'ES {\bf 123} (2016), 199--282.
\bibitem[AG]{AG}
   M.~Abouzaid, S.~Ganatra, in preparation.
\bibitem[AS]{AS}
   M.~Abouzaid, P.~Seidel,
   {\sl Lefschetz fibration methods in wrapped Floer theory}, in preparation.
\bibitem[AV]{AV}
   M.\ Aganagic, C.\ Vafa,
   {\sl Mirror symmetry, D-branes and counting holomorphic discs},
   arXiv:hep-th/0012041.
\bibitem[AuSm]{AuSm}
   D.~Auroux, I.~Smith,
   {\sl Fukaya categories of surfaces, spherical objects and mapping class groups},
   Forum Math.\ Sigma {\bf 9} (2021), e26, 1--50.
\bibitem[Ca]{Cannizzo}
   C.~Cannizzo,
   {\sl Categorical mirror symmetry on cohomology for a complex genus 2 curve},
   Adv.\ Math. {\bf 375} (2020), 107392.
\bibitem[CLL]{ChanLauLeung}
   K.\ Chan, S.-C.\ Lau, N.\ C.\ Leung,
   {\sl SYZ mirror symmetry for toric Calabi-Yau manifolds},
   J.\ Differential Geom.\ {\bf 90} (2012), 177--250.
\bibitem[CL]{CL}
   J.~Chuang, A.~Lazarev, 
   {\sl On the perturbation algebra}, 
   Journal of Algebra {\bf 519} (2019), 130--148.
\bibitem[Cla]{Clarke}
   P.\ Clarke,
   {\sl Duality for toric Landau-Ginzburg models},
   arXiv:0803.0447.
\bibitem[Ef]{EfimovGenusG}
   A.~I.~Efimov, {\sl Homological mirror symmetry for curves of higher genus},
   Adv.\ Math.\ {\bf 230} (2012), 493--530.
\bibitem[ENS]{ENS}
   J.~B.~Etnyre, L.~L.~Ng, J.~M.~Sabloff,
   {\sl Invariants of Legendrian knots and coherent orientations},
   J.\ Symplectic Geom.\ {\bf 1} (2002), 321--367.
\bibitem[Fa]{Fa}
   G.~Faltings, 
   {\sl Semistable vector bundles on Mumford curves}, 
   Invent.\ Math.\ {\bf 74} (1983), 199--212.
\bibitem[GS]{GammageShende}
   B.~Gammage, V.~Shende,
   {\sl Mirror symmetry for very affine hypersurfaces},
   arXiv:1707.02959.
\bibitem[GvdP]{Gerritzen-VdPut}
   L.\ Gerritzen, M.\ van der Put,
   {\it Schottky groups and Mumford curves},
   Lect.\ Notes in Math.\ {\bf 817}, Springer, 1980.
\bibitem[GKR]{GKR}
   M.~Gross, L.~Katzarkov, H.~Ruddat, 
   {\sl Towards mirror symmetry for varieties of general type},
   Adv.\ Math.\ {\bf 308} (2017), 208--275.
\bibitem[Ha]{Hanlon}
   A.~Hanlon, {\sl Monodromy of monomially admissible Fukaya-Seidel categories mirror to toric
varieties}, Adv.\ Math.\ {\bf 350} (2019), 662--746.
\bibitem[Hi]{Hicks}
   J.~Hicks, {\sl Tropical Lagrangian hypersurfaces are unobstructed},
   J.\ Topology {\bf 13} (2020), 1409--1454.
\bibitem[HV]{HoriVafa}
   K.~Hori, C.~Vafa, {\sl Mirror symmetry},
   arXiv:hep-th/0002222.
\bibitem[Lee]{Lee}
   H.~Lee, 
   {\sl Homological mirror symmetry for open Riemann surfaces from pair-of-pants decompositions}, 
   arXiv:1608.04473.
\bibitem[LP]{LekiliPolishchuk}
   Y.~Lekili, A.~Polishchuk,
   {\sl Auslander orders over nodal stacky curves and partially wrapped Fukaya categories}, 
   J.\ Topology {\bf 11} (2018), 615--644.
\bibitem[Ma]{Matessi}
   D.~Matessi, {\sl Lagrangian pairs of pants}, Int.\ Math.\ Res.\ Not.\
   {\bf 2021}, 11306--11356.
\bibitem[Mi]{Mikhalkin}
   G.~Mikhalkin, {\sl Examples of tropical-to-Lagrangian correspondence},
   European J.\ Math.\ {\bf 5} (2019), 1033--1066.
\bibitem[Na]{Nadler}
   D.~Nadler,
   {\sl Mirror symmetry for the Landau-Ginzburg A-model $M=\C^n$, $W=z_1\dots z_n$},
   Duke Math.\ J.\ {\bf 168} (2019), 1--84.
\bibitem[Or]{Orlov}
   D.~Orlov, {\sl Triangulated categories of singularities and equivalences
   between Landau-Ginzburg models}, Sb.\ Math.\ {\bf 197} (2006), 1827--1840.
\bibitem[PZ]{PZ}
   A.~Polishchuk, E.~Zaslow,
   {\sl Categorical mirror symmetry: the elliptic curve},
   Adv.\ Theor.\ Math.\ Phys.\ {\bf 2} (1998), 443--470.
\bibitem[Ru]{Ruddat}
   H.~Ruddat,
   {\sl Perverse curves and mirror symmetry},
   J.\ Algebraic Geometry {\bf 26} (2017), 17--42.
\bibitem[Se1]{SeBook}
   P.~Seidel,
   {\it Fukaya categories and Picard-Lefschetz theory},
   Zurich Lect.\ in Adv.\ Math., European Math.\ Soc., Z\"urich, 2008.
\bibitem[Se2]{SeGenus2}
   P.~Seidel, {\sl Homological mirror symmetry for the genus two curve},  
   J. Algebraic Geometry {\bf 20} (2011), 727--769.
\bibitem[STZ]{STZ}
   N.~Sibilla, D.~Treumann, E.~Zaslow,
   {\sl Ribbon graphs and mirror symmetry},
   Selecta Math.\ {\bf 20} (2014), 979--1002.
\end{thebibliography}
\end{document}